%% file: sdof_ir.tex
\documentclass{article}

\usepackage[margin=1in]{geometry}
\usepackage{amsmath}
\usepackage{amssymb}
\usepackage{graphicx}
\usepackage{subcaption}
\usepackage{authblk}
\usepackage[firstinits=true, sorting=none, style=numeric-comp, backend=bibtex, maxnames=999]{biblatex}
\usepackage{hyperref}
\usepackage{cleveref}
\usepackage{xcolor}
\usepackage{multirow}
\usepackage{soul} \usepackage{comment} \usepackage{placeins} \usepackage{multirow}

\newcommand{\removeBib}[1]{\AtEveryBibitem{\clearfield{#1}} \AtEveryBibitem{\clearlist{#1}} \AtEveryBibitem{\clearname{#1}}} 

\removeBib{language}
\removeBib{month}
\removeBib{day}

\DeclareFieldFormat{titlecase}{\MakeTitleCase{#1}}

\newrobustcmd{\MakeTitleCase}[1]{\ifthenelse{\ifcurrentfield{booktitle}\OR\ifcurrentfield{booksubtitle}\OR\ifcurrentfield{maintitle}\OR\ifcurrentfield{mainsubtitle}\OR\ifcurrentfield{journaltitle}\OR\ifcurrentfield{journalsubtitle}\OR\ifcurrentfield{issuetitle}\OR\ifcurrentfield{issuesubtitle}\OR\ifentrytype{book}\OR\ifentrytype{mvbook}\OR\ifentrytype{bookinbook}\OR\ifentrytype{booklet}\OR\ifentrytype{suppbook}\OR\ifentrytype{collection}\OR\ifentrytype{mvcollection}\OR\ifentrytype{suppcollection}\OR\ifentrytype{manual}\OR\ifentrytype{periodical}\OR\ifentrytype{suppperiodical}\OR\ifentrytype{proceedings}\OR\ifentrytype{mvproceedings}\OR\ifentrytype{reference}\OR\ifentrytype{mvreference}\OR\ifentrytype{report}} {#1}
{\MakeSentenceCase{#1}}}

\providecommand{\keywords}[1]{\textbf{\textit{Keywords---}} #1}

\title{Tracking Superharmonic Resonances for Nonlinear Vibration of Conservative and Hysteretic Single Degree of Freedom Systems}
\author[1]{Justin H. Porter~\thanks{jp88@rice.edu}}
\author[1]{Matthew R. W. Brake~\thanks{brake@rice.edu}}

\affil[1]{Department of Mechanical Engineering, Rice University, Houston,
  TX 77005}

\begin{document}

\include{preamble.tex}

\pagebreak

\pagenumbering{arabic} 

\maketitle{}
\begin{abstract}

Many modern engineering structures exhibit nonlinear vibration. Characterizing such vibrations efficiently is critical to optimizing designs for reliability and performance. 
For linear systems, steady-state vibration occurs only at the forcing frequencies. 
However, nonlinearities (e.g., contact, friction, large deformation, etc.) can result in nonlinear vibration behavior including superharmonics - responses at integer multiples of the forcing frequency. 
When the forcing frequency is near an integer fraction of the natural frequency, superharmonic resonance occurs, and the magnitude of the superharmonics can exceed that of the fundamental harmonic that is externally forced. 
Characterizing such superharmonic resonances is critical to improving engineering designs. 
The present work extends the concept of phase resonance nonlinear modes (PRNM) to be applicable to general nonlinearities, and is demonstrated for eight different nonlinear forces. The considered forces include stiffening, softening, contact, damping, and frictional nonlinearities that have not been previously considered with PRNM.
The proposed variable phase resonance nonlinear modes (VPRNM) method can accurately track superharmonic resonances for hysteretic nonlinearities that exhibit amplitude dependent phase resonance conditions that cannot be captured by PRNM.
The proposed method allows for characterization of superharmonic resonances without constructing a full frequency response curve at every force level with the harmonic balance method. Thus, the present method allows for analysis of potential failures due to large amplitudes near the superharmonic resonance with reduced computational cost.
The consideration of single degree of freedom systems in the present paper provides insights into superharmonic resonances and a basis for understanding internal resonances for multiple degree of freedom systems.

\end{abstract}

\section{Introduction}
\label{sec:introduction}

The optimization of modern engineering structures to improve efficiency requires the consideration of nonlinearities such as those due to large deformation \cite{touzeNonlinearNormalModes2006, mignoletReviewIndirectNonintrusive2013} and friction in jointed connections \cite{brakeMechanicsJointedStructures2017, mathisReviewDampingModels2020}. 
These nonlinearities result in a variety of vibration phenomena not observed in linear systems, such as amplitude dependent frequency and damping properties. Additionally, harmonically forced nonlinear systems can exhibit superharmonic resonances where responses of an unforced higher harmonic at an integer multiple of the forcing frequency have amplitudes that can be on the order of or higher than the fundamental harmonic response. 
Superharmonic resonances can occur because nonlinear internal forces excite higher harmonics of the structure \cite{ferriFrequencyDomainSolutions1988, chenPeriodicResponseBlades1999, chenPERIODICFORCEDRESPONSE2000, krackHarmonicBalanceNonlinear2019}, and are not found in linear systems, which respond in steady-state at only the forcing frequency.
Following the definition of \cite{volvertPhaseResonanceNonlinear2021}, superharmonic resonances are defined as local maximums in the amplitude of a superharmonic component, which is a response at an integer multiple of the forcing frequency. These resonances are accompanied by phase shifts in the superharmonic that can be used to define a superharmonic phase resonance condition near the amplitude resonance \cite{volvertPhaseResonanceNonlinear2021}.

Extensive literature has documented superharmonic resonances for conservative nonlinearities including in experiments \cite{nayfehModalInteractionsDynamical1989}. 
In addition, more recent experiments have demonstrated superharmonic resonances for frictional nonlinearities \cite{claeysModalInteractionsDue2016, chenMeasurementIdentificationNonlinear2022, scheelChallengingExperimentalNonlinear2020}.
Given the significant amplitude of the higher harmonics, superharmonic resonances can have amplitudes that are easily twice that of single harmonic solutions over the same frequency region. Therefore, it is critical to understand the superharmonic resonances if a structure experiences forcing near an integer fraction of a resonant frequency (e.g., due to rotation in turbomachinery).
Furthermore, these responses occur away from primary resonances and could easily be missed by design processes focusing on characterizing primary resonances to avoid large vibration amplitudes.
Additionally, superharmonic resonances can sometimes result in an internal resonance of two modes vibrating at an integer ratio of frequencies \cite{nayfehEnergyTransferHighFrequency1995}.
In these cases, it is critical to characterize the internal resonance because it occurs in the resonant regime of the mode responding at the fundamental frequency and thus could correspond to a globally maximum response amplitude.

Superharmonic resonances and related internal resonances can be modeled with numerous techniques including the harmonic balance method (HBM) \cite{krackHarmonicBalanceNonlinear2019}, phase resonance nonlinear modes (PRNM) \cite{volvertPhaseResonanceNonlinear2021}, perturbation techniques \cite{nayfehNonlinearInteractionsAnalytical2000, nayfehNonlinearOscillations1995}, and invariant manifolds \cite{boivinNonlinearModalAnalysis1995, pesheckNonlinearModalAnalysis2001, jiangConstructionNonlinearNormal2005, liNonlinearAnalysisForced2022}. 
In addition, methods to characterize nonlinear normal modes can include significant higher harmonic components for some energy levels \cite{kerschenNonlinearNormalModes2009, peetersNonlinearNormalModes2009, krackNonlinearModalAnalysis2015, krackVibrationPredictionBladed2017}.
Perturbation techniques have been applied to numerous examples of superharmonic and internal resonances (e.g., see \cite{nayfehNonlinearInteractionsAnalytical2000}). 
Similarly, Melnikov-type analyses utilizing perturbations of solutions have derived criterion that are required for subharmonic and superharmonic resonances to occur \cite{cenedeseHowConservativeBackbone2020, guckenheimerAveragingPerturbationGeometric1983}.
However, perturbation and other analytical techniques can only be applied to weak nonlinearities. 
On the other hand, HBM has frequently been used to model nonlinear vibration systems under steady-state external forcing numerically and does not limit the strength of the nonlinearities. The inclusion of higher harmonics in HBM allows for the modeling of superharmonic resonances \cite{krackHarmonicBalanceNonlinear2019}.

While analytical methods exist for analyzing superharmonic and internal resonances, numerical approaches for efficiently understanding superharmonic resonances are still an open area of research. 
Recently, PRNM has been proposed as a method of tracking superharmonic, subharmonic, and ultra-subharmonic resonances in single degree of freedom (SDOF) and multiple degree of freedom (MDOF) systems \cite{volvertPhaseResonanceNonlinear2021}. However, this approach has only been applied for conservative cubic \cite{volvertPhaseResonanceNonlinear2021, volvertResonantPhaseLags2022, abeloosControlbasedMethodsIdentification2022a} and quadratic \cite{abeloosControlbasedMethodsIdentification2022a} nonlinearities. 
The advantage of PRNM compared to HBM is that one continuation can be conducted to calculate the superharmonic resonance responses across many force levels while evaluating the solution at a given force level for only a single forcing frequency. By contrast, HBM generally requires continuation over a range of frequencies for a discrete set of force levels, frequently requiring at least an order of magnitude more nonlinear solutions than PRNM to obtain a similar characterization. Therefore, PRNM provides an understanding of superhamonic resonances at significantly lower computational costs. 
PRNM has been demonstrated to achieve similar benefits for experimental tracking of superharmonic resonances under some assumptions about the form of the nonlinearity \cite{abeloosControlbasedMethodsIdentification2022a}.
Furthermore, PRNM presents an opportunity to address challenges noted in other work related to tracking resonance characteristics for structures across force levels in the presence of superharmonics \cite{krackVibrationPredictionBladed2017, krackMultiharmonicAnalysisDesign2012}.
The application of PRNM to more general nonlinearities (e.g., hysteretic nonlinearities) has not been evaluated.
Other recent work has demonstrated a new tracking method for response extrema across force levels, but that approach required sufficiently smooth nonlinearities and was more computationally expensive than phase resonance based approaches \cite{razeTrackingAmplitudeExtrema2024}.

Hysteretic nonlinearities are an important category of nonlinear forces that have not previously been considered with PRNM. 
Jointed connections are commonly modeled with hysteretic nonlinearities making hysteretic nonlinearities highly relevant for engineering applications \cite{brakeMechanicsJointedStructures2017, mathisReviewDampingModels2020}. 
Hysteretic nonlinearities display path dependencies (generally formulated as history variables) and thus the nonlinear forces cannot be evaluated based solely on the instantaneous displacement and velocity. 
This additional complexity generally precludes the use of analytical solutions. 
The resulting damping from hysteretic nonlinearities is nontrivial to remove, preventing the use of techniques developed for conservative systems.
Specific methods for modeling jointed connections, such as the Extended Periodic Motion Concept (EPMC) \cite{krackNonlinearModalAnalysis2015}, have been proposed to model assembled structures more efficiently compared to HBM. However, EPMC breaks down in the presence of superharmonics and internal resonances. Thus, a different approach is necessary when superharmonic resonances occur.

Recent experiments demonstrating superharmonic resonances for hysteretic systems \cite{claeysModalInteractionsDue2016, chenMeasurementIdentificationNonlinear2022, scheelChallengingExperimentalNonlinear2020} motivate the need to generalize superharmonic resonance tracking to allow for general nonlinearities and specifically hysteretic nonlinearities. 
Furthermore, several studies have numerically modeled structures with hysteretic nonlinearities that exhibit superharmonic resonances \cite{claeysModalInteractionsDue2016, wongSteadyStateOscillation1994, masianiResonantCoupledResponse2002, houBifurcationStabilityAnalysis2019, casiniMitigationStructuralVibrations2022, casiniRoleHystereticRestoring2022}.
Accurately characterizing superharmonic resonances for jointed structures is critical because the damping properties can change significantly with the amplitude and phase of the higher harmonics \cite{krackEfficacyFrictionDamping2016}.
Furthermore, superharmonic resonances have the potential to cause significant amplitudes away from primary resonances. 
Additionally, superharmonic resonances between two modes resulting in an internal resonances can significantly change the resonant response characteristics of a structure with friction \cite{krackEfficacyFrictionDamping2016}.

The present paper considers tracking superharmonic resonances for SDOF systems with eight different nonlinear forces including conservative stiffening\footnote{The term stiffening is adopted in the present work rather than hardening because it more precisely reflects the change in stiffness and to avoid confusion with how hardening refers to changes in material hardness independent of the stiffness in other related communities.} and softening nonlinearities, a unilateral spring, cubic damping, and two hysteretic nonlinearities. Of these eight nonlinear forces, only the stiffening and softening cubic nonlinearities have been previously discussed in the literature for PRNM.
\Cref{sec:system} introduces the vibrating system, the nonlinear forces, and provides an introduction to the superharmonic resonances behavior. 
Results from PRNM are presented in additional detail in \Cref{sec:modeling} and a decomposition of the nonlinear forces is derived as a basis for the present work.
\Cref{sec:apriori_phase} utilizes the force decomposition to provide an understanding of the superharmonic resonance behavior. Then, a new method for tracking superharmonic resonances is proposed in \Cref{sec:tracking_method}, and in \Cref{sec:results}, the
evolutions of superharmonic resonances with respect to varying external force amplitude are discussed. Finally, the major conclusions of the paper are presented (see \Cref{sec:concl}).
While only SDOF systems are analyzed in the present work, this represents an important step towards developing methods for MDOF systems with general nonlinear forces. Furthermore, the present work formulates the proposed method such that the approach could be applied to MDOF systems with minimal modifications, but such analyses and determinations of their accuracy are left to future work. The restriction to SDOF systems in the present work means that only superharmonic resonances and not related internal resonances are considered.

\section{System Description} \label{sec:system}

The present work considers the nonlinear vibration behavior of a single degree of freedom (SDOF) system with a nonlinear internal force $f_{nl}$ and external forcing with magnitude $ F $ and frequency $ \omega $. The equations of motion of the forced system are
\begin{equation} \label{eq:system}
	m \ddot{x} + c \dot{x} + kx + f_{nl}(x, \dot{x}) = F \cos(\omega t)
\end{equation}
where $ x $, $ \dot{x} $, and $ \ddot{x} $ are the displacement, velocity, and acceleration of the mass respectively. The system has mass $ m $, damping factor $ c $, and linear stiffness $ k $. 
The mass is fixed at $ m = 1.0 $ kg for all models and the linear stiffness $ k $ is varied such that the system including the linearized stiffness from the nonlinear force has a natural frequency of 1 rad/s.
A linear damping factor of $c = 0.01 $ kg/s is selected for all cases.
The present work considers eight different $ f_{nl} $ including conservative stiffening (Duffing and quintic stiffness), conservative softening (softening Duffing and Iwan inspired (II) softening), even (unilateral spring as formulated in \Cref{sec:unispring_force}), nonlinear damping (cubic damping), and hysteretic (Jenkins element and Iwan element) behavior. The full details of the nonlinear forces are presented in \Cref{sec:nlforces}. 
The hysteretic nonlinear forces require numerical treatment, which has not been previously addressed with PRNM \cite{volvertPhaseResonanceNonlinear2021,volvertResonantPhaseLags2022, abeloosControlbasedMethodsIdentification2022a}. Furthermore, of the considered nonlinear forces, only the Duffing nonlinearities have been analyzed previously with PRNM \cite{volvertPhaseResonanceNonlinear2021,volvertResonantPhaseLags2022, abeloosControlbasedMethodsIdentification2022a}.
Frequency response curves (FRCs) for the nonlinear system of \eqref{eq:system} are calculated using the HBM and continuation (see \Cref{sec:hbm} for more details and \Cref{sec:frcs_main} for example FRCs) \cite{krackHarmonicBalanceNonlinear2019}. The code for the present simulations is made available for reference \cite{porterTMDSimPy}.

\subsection{Nonlinear Forces}\label{sec:nlforces}

The present section summarizes the nonlinear forces used in this study. 
The form of the nonlinearities and the parameters are tabulated in Tables \ref{tab:nl_forces_refs} and \ref{tab:nl_params} respectively.
To generalize the results, values are nondimensionalized utilizing the mass $m$, linearized stiffness $k_{lin}$, and a reference displacement $x_{ref}$ with dimensional values presented in \Cref{tab:nl_forces_refs}. The linearized stiffness is defined as 
\begin{equation}
	k_{lin} = k + \dfrac{\partial f_{nl}}{\partial x} \bigg|_{x=0, \dot{x}=0}.
\end{equation}
The parameters for the nonlinear forces are nondimensionalized as
\begin{equation}
	\begin{split}
		\omega_0 &= \sqrt{k_{lin} / m}, \ \  
		\zeta_0 = c / (2 \sqrt{k_{lin} m}),
		\\
		\hat{k} &= k / k_{lin}, \ \ 
		\hat{k}_t = k_t / k_{lin}, \ \ 
		\hat{k}_{nl} = k_{nl} / k_{lin},
		\\
		\hat{\alpha} &= \alpha x_{ref}^2 / k_{lin}, \ \ 
		\hat{\eta} = \eta x_{ref}^4 / k_{lin} , \ \  
		\hat{\gamma} = \gamma (\omega_0 x_{ref})^3/ (k_{lin} x_{ref}),
		\\
		\hat{F}_s &= F_s / (k_{lin} x_{ref}).
	\end{split}
\end{equation}
The nondimensionalizations for plotting are 
\begin{equation}
	\begin{split}
		\hat{\omega} &= \omega / \omega_0  , \ \ \hat{t} = t / \omega_0,
		\\
		|\hat{X}| &= |X| / x_{ref}, \ \ 
 		\hat{x} = x / x_{ref}  ,
 		\\
		\hat{F} &= F / (k_{lin} x_{ref}), \ \ 
		\hat{f}_{ext} = F \cos(\omega t) / (k_{lin} x_{ref}), \ \ 
		\hat{f}_{nl} = f_{nl} / (k_{lin} x_{ref}).
	\end{split}
\end{equation}
Here, $|X|$ is the maximum value of $x(t)$ at any time over the cycle.

\FloatBarrier

\begin{table}[h!]
	\centering
	\caption{Table of nonlinear forces and reference quantities.}
	\label{tab:nl_forces_refs}
		\begin{tabular}{ccccc}
			\hline \hline
			Force & Form & $m$ [kg] & $k_{lin} $ [N/m] & $x_{ref}$ [m] 
			\\ \hline
			
			Stiffening Duffing & $ \alpha x^3 $ & 1  & 1  & 1
			\\
			Quintic Stiffness & $ \eta x^5 $ & 1  & 1  & 1 
			\\ \hline
			
			Softening Duffing & $ \alpha x^3 $ & 1  & 1  & 1 
			\\
			Conservative Softening II & \eqref{eq:numericalsoften} & 1  & 1  & $\phi_{max} = 1.6 $ from \eqref{eq:phi_max} 
			\\ \hline
			
			Unilateral Spring & $ \max(k_{nl} x , 0) $ & 1  & 1  & 1 
			\\ \hline

			Cubic Damping & $ \gamma \dot{x}^3 $ & 1  & 1  & 1 
			\\
			Jenkins Element & \eqref{eq:jenk_force} & 1  & 1  & $F_s/k_t = 0.8 $ 
			\\
			Iwan Element & \eqref{eq:iwan_distrib} \eqref{eq:slider_f_iwan} \eqref{eq:iwan_integral} & 1  & 1  & $\phi_{max} = 2.4 $ from \eqref{eq:phi_max} 
			
			\\ \hline\hline
		\end{tabular}
\end{table}

\begin{table}[h!]
	\centering
	\caption{Parameters for vibration systems and different nonlinear forces.}
	\label{tab:nl_params}
		\begin{tabular}{cccc}
			\hline \hline
			Force & $\zeta_0$ & $\hat{k}$ & Nondimensionalized Parameters 
			\\ \hline
			
			Stiffening Duffing & 0.005 & 1 & $ \hat{\alpha} = 1 $
			\\
			Quintic Stiffness & 0.005 & 1 & $ \hat{\eta} = 1 $ 
			\\ \hline
			
			Softening Duffing & 0.005 & 1 & $ \hat{\alpha} = -2.5 $e-4 
			\\
			Conservative Softening II$^*$ & 0.005 & 0.75 & $ \hat{k}_t = 0.25 $, $ \hat{F}_s = 0.125 $, $ \chi = \beta = 0 $
			\\ \hline
			
			Unilateral Spring & 0.005 & 0.75 & $ \hat{k}_{nl} = 0.5 $
			\\ \hline

			Cubic Damping & 0.005 & 1 & $ \hat{\gamma} = 0.03 $
			\\
			Jenkins Element$^*$ & 0.005 & 0.75 & $ \hat{k}_t = 0.25 $, $ \hat{F}_s = 0.25 $
			\\
			Iwan Element$^*$ & 0.005 & 0.75 & $ \hat{k}_t = 0.25 $, $ \hat{F}_s = 0.083333 $, $ \chi = -0.5 $, $\beta = 0 $

			\\ \hline\hline
			\multicolumn{4}{l}{\begin{tabular}[c]{@{}l@{}}*The Conservative Softening II, Jenkins Element, and Iwan Element models all use the same dimensional
			\\  
			value of $F_s = 0.2 $ [N], but differ in $x_{ref}$ and thus have different values of $\hat{F}_s$ listed here.\end{tabular}}
		\end{tabular}
\end{table}

\FloatBarrier

\subsubsection{Conservative Softening II Nonlinearity}

The conservative softening Iwan inspired (II) nonlinearity is based on the loading backbone of the 4-parameter Iwan model and has form of \cite{segalmanFourParameterIwanModel2005}
\begin{equation} \label{eq:numericalsoften}
	f_{nl} = 
	\begin{cases}
		k_t x - \left( \left(  \dfrac{k_t (\beta + \frac{\chi+1}{\chi+2}  ) }{ F_s(1+\beta) }  \right)^{1+\chi}
		\dfrac{k_t}{(1 + \beta)(\chi + 2) }  \right) |x|^{\chi + 2} \text{sgn}(x) &  |x| < \phi_{max}
		\\
		F_s \text{sgn}(x) &  |x| \geq \phi_{max}
	\end{cases}
	.
\end{equation}
Here, $ k_t $ is the initial stiffness of the nonlinear element, $ F_s $ is the saturation limit, $ \chi $ controls the shape of the softening part of the curve, and $ \beta $ defines the extent of a slope discontinuity at the point where the force reaches $ F_s $. In addition, $ \text{sgn}(\cdot) $ is the signum function and the parameter $ \phi_{max} $ is the displacement for the transition to a constant force and is calculated as
\begin{equation} \label{eq:phi_max}
	\phi_{max} = \dfrac{F_s (1 + \beta)}{k_t \left( \beta + \frac{\chi + 1}{\chi + 2} \right)}.
\end{equation}
This model is chosen since the 4-parameter Iwan model is popular for the modeling of bolted connections \cite{segalmanFourParameterIwanModel2005, mathisReviewDampingModels2020}. The conservative softening II nonlinearity retains some of the characteristics of the 4-parameter Iwan model while simplifying the force to be conservative.

\subsubsection{Unilateral Spring} \label{sec:unispring_force}

Next, the unilateral spring in \Cref{tab:nl_forces_refs} is an even nonlinearity as can be seen by reformulating the linear and nonlinear stiffness terms as
\begin{equation}
	kx + \max(k_{nl} x, 0) = \left[ k + \dfrac{k_{nl}}{2} \right] x + \dfrac{k_{nl}}{2} |x|.
\end{equation}
Furthermore, the derivative of the nonlinear force is not defined at zero displacement. For the nondimensionalization, the present work uses
\begin{equation}
	k_{lin} = k + \dfrac{k_{nl}}{2}.
\end{equation}

\subsubsection{Jenkins Element}

The Jenkins hysteretic nonlinearity is a stick slip element with initial stiffness of $ k_t $ and slip limit of $ F_s $ \cite{jenkinsAnalysisStressstrainRelationships1962}. The evolution of the element is calculated based on the previous values of displacement $ x_0 $ and force $ f_0 $ as
\begin{equation} \label{eq:jenk_force}
	f_{nl}(x) = \begin{cases}
		\underbrace{k_t (x - x_0) + f_0}_{f_{stuck}} & |f_{stuck}| < F_s 
		\\
		F_s\ \text{sgn}(f_{stuck}) & Otherwise
	\end{cases}
	.
\end{equation}

\subsubsection{Iwan Element}

Finally, the 4-parameter Iwan model is used for a second hysteretic nonlinearity. This element has a distribution of sliders with strengths $ \phi \in [0, \phi_{max}] $ of \cite{segalmanFourParameterIwanModel2005}
\begin{equation} \label{eq:iwan_distrib}
	\rho(\phi) = \dfrac{F_s (\chi + 1)}{\phi_{max}^{\chi + 2}  \left( \beta + \frac{\chi + 1}{\chi + 2} \right) } \phi^\chi 
		+ 
		\dfrac{F_s \beta}{\phi_{max}  \left( \beta + \frac{\chi + 1}{\chi + 2} \right) } 
	\delta(\phi-\phi_{max})
\end{equation}
where $ \delta(\cdot) $ is a Dirac delta function.
The contribution of each slider $ f_\phi $ is calculated similar to the Jenkins model as\footnote{Note that force contributed by slider $ \phi $ is $ k_t f_\phi $ since the tangential stiffness incorporated into the total force through $ \rho(\phi) $ and $ \phi_{max} $.}
\begin{equation} \label{eq:slider_f_iwan}
	f_\phi = \begin{cases}
		\underbrace{x - x_0 + f_{\phi,0}  }_{f_{\phi,stuck}} & |f_{\phi,stuck}| < \phi
		\\
		\phi \text{sgn}(f_{\phi,stuck}) & Otherwise
	\end{cases}
\end{equation}
with $ f_{\phi, 0} $ representing the value at the previous instant.
Then, the force is calculated as
\begin{equation}\label{eq:iwan_integral}
	f_{nl} = 
	\int_{0}^{\phi_{max}}   f_\phi	\rho(\phi) d\phi  .
\end{equation}
For the purposes of this work, this integral is discretized with 100 equally spaced sliders and the midpoint integration rule for the continuous portion of the distribution in \eqref{eq:iwan_distrib} plus a single slider with value $ \phi_{max} $ for the Dirac delta function in \eqref{eq:iwan_distrib}.

\FloatBarrier

\subsection{Example Frequency Response Curves} \label{sec:frcs_main}

\begin{figure}[!h]
	\centering
	\newcommand{\exwidth}{0.85}  
	\begin{subfigure}{\exwidth\linewidth}
		\centering
		\includegraphics[width=\linewidth]{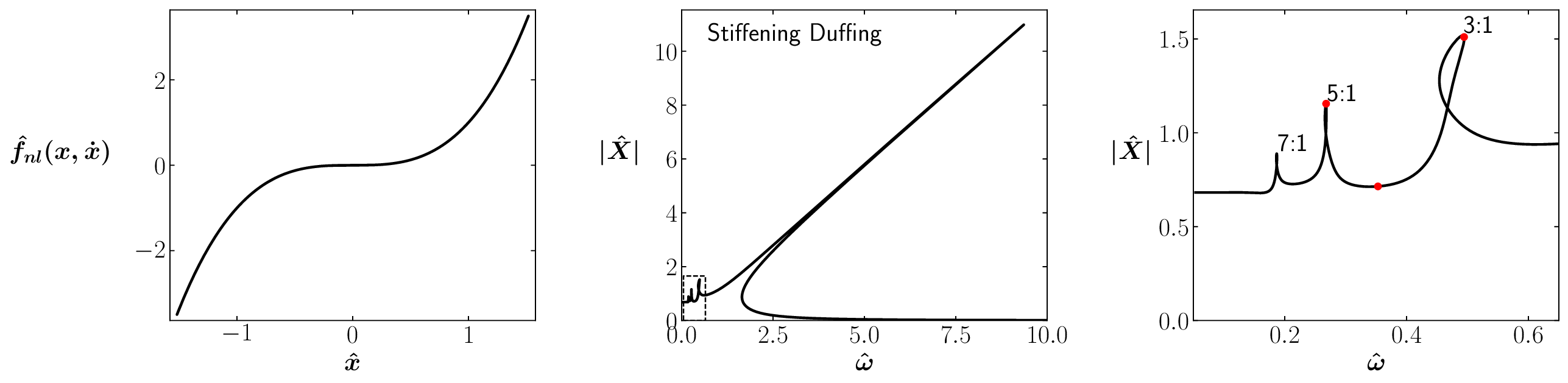}
		\caption{}\label{fig:example_FRC_duff}
	\end{subfigure}
												\begin{subfigure}{\exwidth\linewidth}
		\centering
		\includegraphics[width=\linewidth]{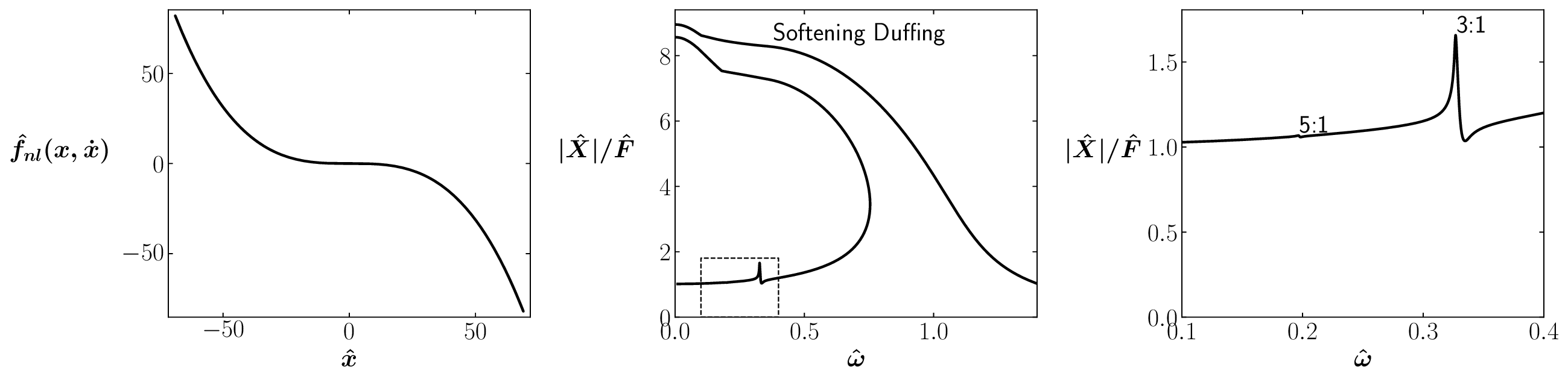}
		\caption{}\label{fig:example_FRC_softduff}
	\end{subfigure}
								\begin{subfigure}{\exwidth\linewidth}
		\centering
		\includegraphics[width=\linewidth]{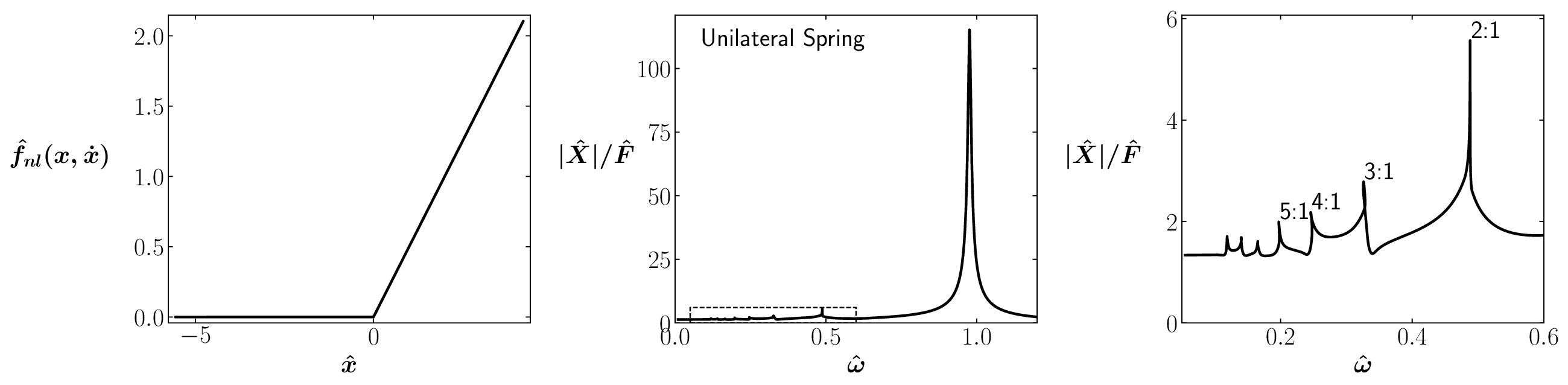}
		\caption{}\label{fig:example_FRC_unispring}
	\end{subfigure}
														\begin{subfigure}{\exwidth\linewidth}
		\centering
		\includegraphics[width=\linewidth]{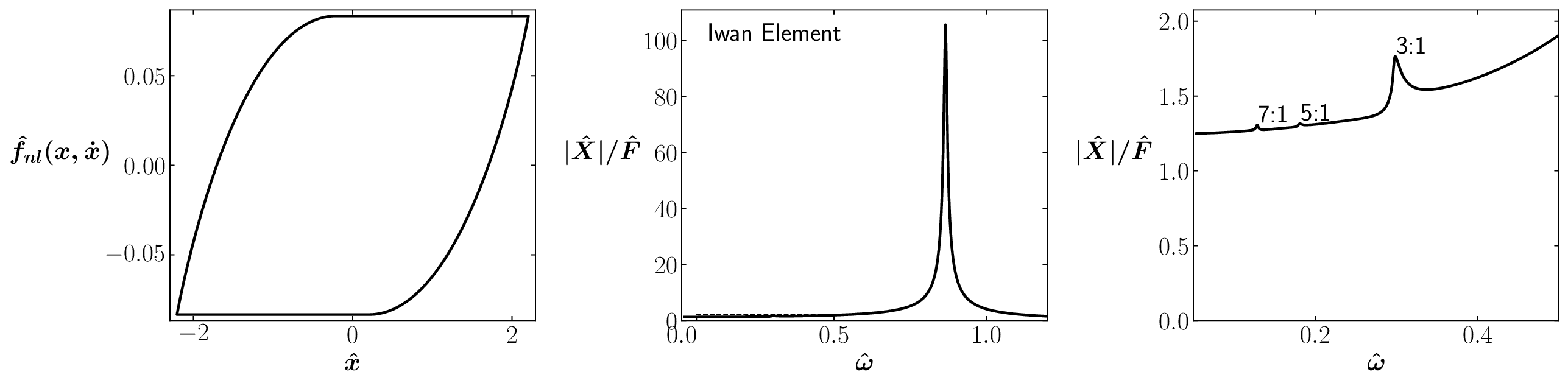}
		\caption{}\label{fig:example_FRC_iwan}
	\end{subfigure}
		\caption{Examples of FRCs illustrating superharmonic behavior calculated with harmonics 0 and 1 through 8. The left is the force displacement relationship for the lowest integer superharmonic resonance, the middle is the FRC over the full frequency range, and the right is a zoomed in version for the dashed box on the middle plot. The systems are 
		(a) stiffening Duffing with $ \hat{F} = 1.0 $, 
				(b) softening Duffing with $ \hat{F} = 8.0 $,
				(c) unilateral spring with $ \hat{F} = 1.0 $,
						and
		(d) Iwan element with $ \hat{F} = 1.25 $. The dots in (a) indicate the time series plots in \Cref{fig:example_time_series}.} 
	\label{fig:example_FRC}
\end{figure}

All of the nonlinear forces result in systems that exhibit superharmonic resonances. An $ n $:1 superharmonic resonance denotes when the $ n $th harmonic (with frequency $ n $ times the forcing frequency) responds at a local peak amplitude at a given forcing frequency. Several FRCs showing the diversity of superharmonic resonances are shown in \Cref{fig:example_FRC} (see \Cref{sec:appendix_frcs} for examples for the other nonlinear forces). 
From \Cref{fig:example_FRC}, it is clear that the qualitative behavior of the superharmonic resonances varies significantly with the different nonlinear forces. For instance in \Cref{fig:example_FRC_unispring}, the unilateral spring displays several more superharmonic resonances than the other nonlinear forces, and the superharmonic resonances at lower frequencies are not labeled since multiple harmonics simultaneously reach a peak value. Also of note, the 5:1 superharmonic resonance of the Iwan model in \Cref{fig:example_FRC_iwan} has a less prominent peak than the 7:1 superharmonic resonance, whereas the opposite is true for the stiffening Duffing nonlinearity (see \Cref{fig:example_FRC_duff}).

The superharmonic peaks shown for the SDOF systems are small relative to the primary resonance peaks. However, they do result in notable amplitudes compared to the response excluding higher harmonics at nearby frequencies (e.g., roughly double the amplitude for the Duffing oscillator at the 3:1 superharmonic resonance rather than at slightly higher or lower frequencies). 
Many structures are designed to avoid primary resonances (e.g., turbomachinery). However, if superharmonic resonances are neglected, the amplitude in some regions away from primary resonances could be significantly underpredicted, potentially resulting in structural failures.
Furthermore, the present work proposes a tracking method as a step towards generalizing to MDOF systems. 
For MDOF systems, it is possible to have two modes in resonance, with one at the fundamental frequency and one at an integer multiple of the fundamental frequency.
In that case, it is important to capture the superharmonic resonance behavior since it occurs in a high amplitude regime for the overall response. 
The SDOF systems considered here are a necessary step for developing methods that can be applied to MDOF systems.

\begin{figure}[!h]
	\centering
	\begin{subfigure}{0.6\linewidth}
		\centering
		\includegraphics[width=\linewidth]{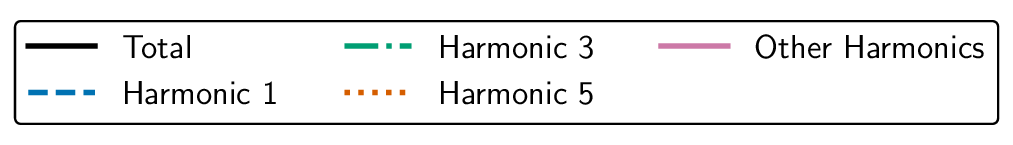}
			\end{subfigure}
	\\
	\begin{subfigure}{0.49\linewidth}
		\centering
		\includegraphics[width=\linewidth]{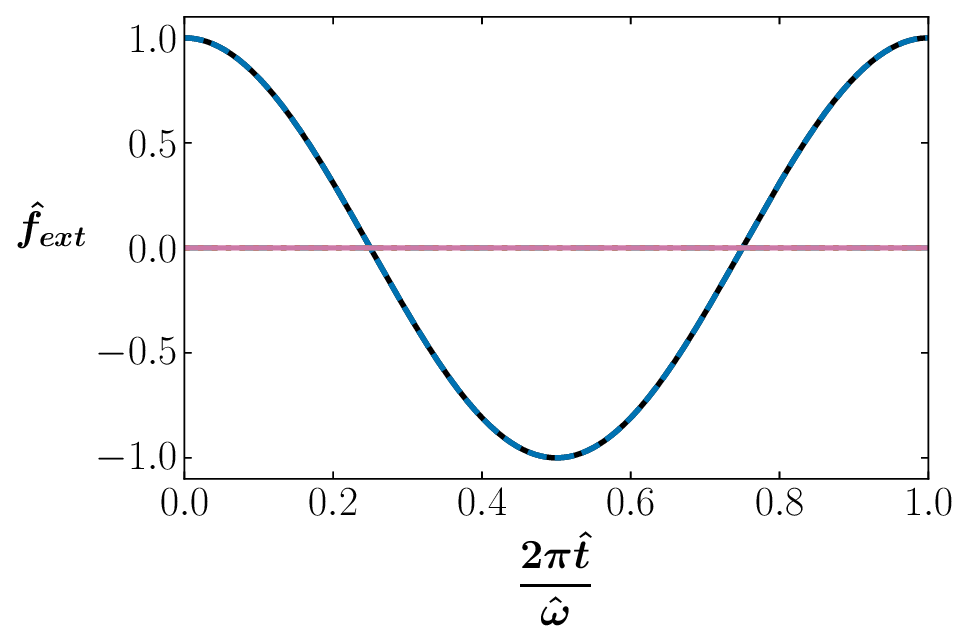}
		\caption{} \label{fig:time_force}
	\end{subfigure}
	\hfill
	\begin{subfigure}{0.49\linewidth}
		\centering
		\includegraphics[width=\linewidth]{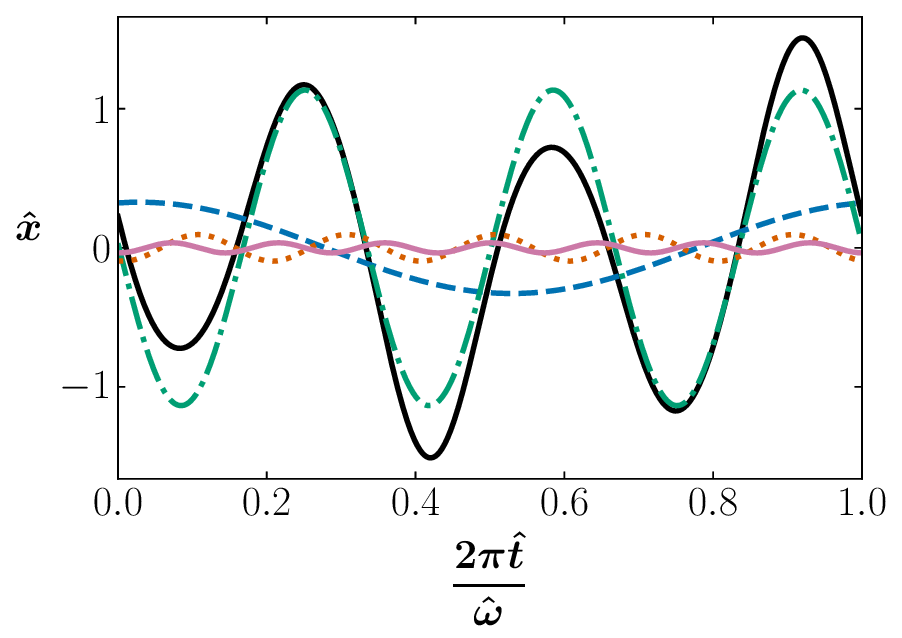}
		\caption{} \label{fig:time_3to1}
	\end{subfigure}
	\\
		\begin{subfigure}{0.49\linewidth}
		\centering
		\includegraphics[width=\linewidth]{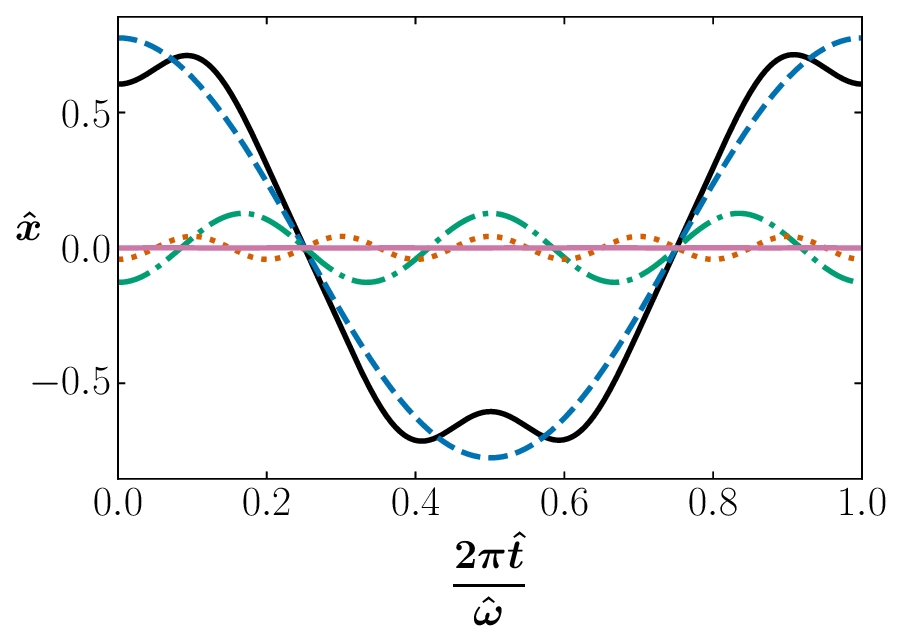}
		\caption{}  \label{fig:time_4to1}
	\end{subfigure}
	\hfill
	\begin{subfigure}{0.49\linewidth}
		\centering
		\includegraphics[width=\linewidth]{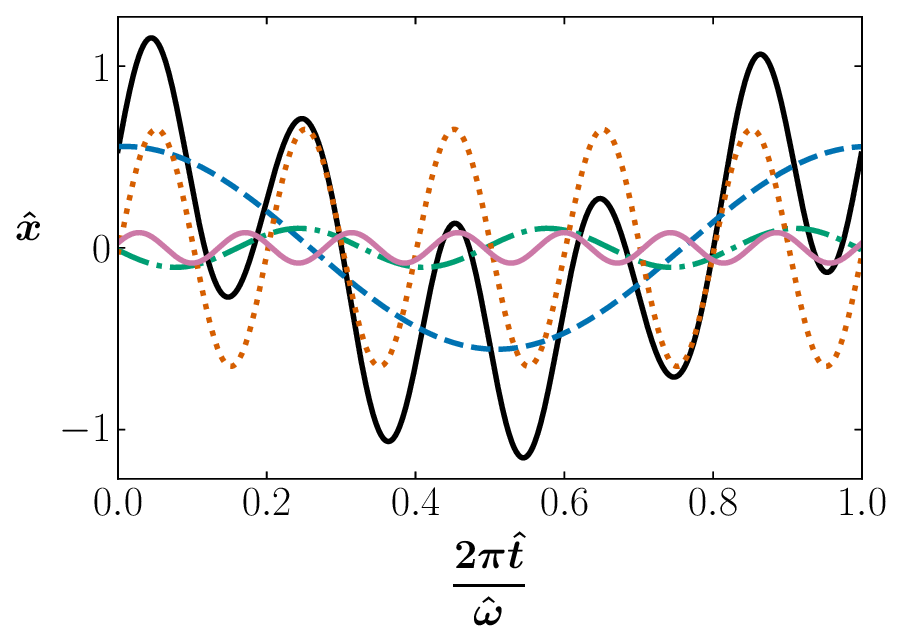}
		\caption{} \label{fig:time_5to1}
	\end{subfigure}
		\caption{Time series plots for the stiffening Duffing nonlinearity shown in \Cref{fig:example_FRC_duff} (corresponding to the red dots) (a) nondimensionalized applied external force $\hat{f}_{ext}$, (b) the 3:1 superharmonic resonance (the peak at $\hat{\omega} = 0.494$ in \Cref{fig:example_FRC_duff}) (c) response at $\hat{\omega} =0.350$ and (d) 5:1 superharmonic resonance (the peak at $\hat{\omega} = 0.268$ in \Cref{fig:example_FRC_duff})} 
	\label{fig:example_time_series}
\end{figure}

\Cref{fig:example_time_series} illustrates the time series of the response for three different responses for the stiffening Duffing nonlinearity shown in \Cref{fig:example_FRC_duff}. In all cases, the external forcing only contains the fundamental harmonic \Cref{fig:time_force} while higher harmonics are included in the responses. For the 3:1 and 5:1 cases (see \Cref{fig:time_3to1} and \Cref{fig:time_5to1} respectively), the third and fifth harmonics provide the largest components to the response resulting in the superharmonic resonance. The system does not show contributions from even harmonics and thus there is no 4:1 superharmonic resonance observed at the intermediate frequency of 0.35 rad/s\footnote{A 4:1 superharmonic resonance would be expected between the 3:1 and 5:1 superharmonic resonances near 0.35 rad/s and at approximately 0.25 times the natural frequency.} (see \Cref{fig:time_4to1}). In \Cref{fig:time_4to1}, the relative phase of the higher harmonics is such that the total response reaches a lower peak amplitude than just the response of the fundamental harmonic. Thus, the relative phase between the harmonics is critical to understanding if the presence of higher harmonics increases or decreases the total amplitude relative to a single harmonic case.

\begin{figure}[!h]
	\centering
	\begin{subfigure}{0.32\linewidth}
		\centering
		\includegraphics[width=\linewidth]{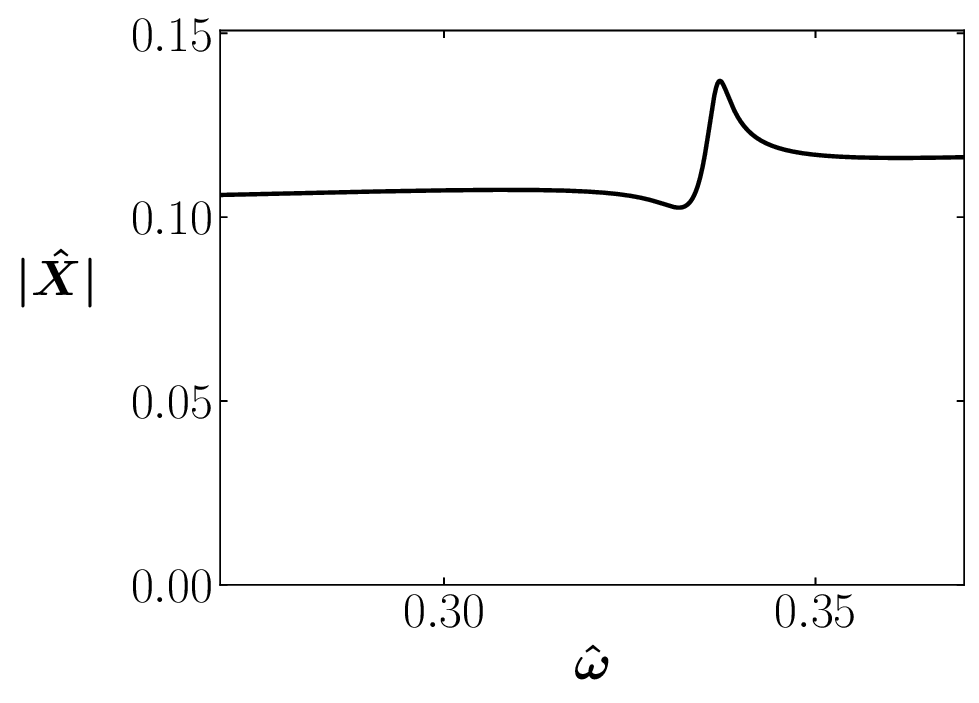}
		\caption{} 
	\end{subfigure}
	\hfill
	\begin{subfigure}{0.32\linewidth}
		\centering
		\includegraphics[width=\linewidth]{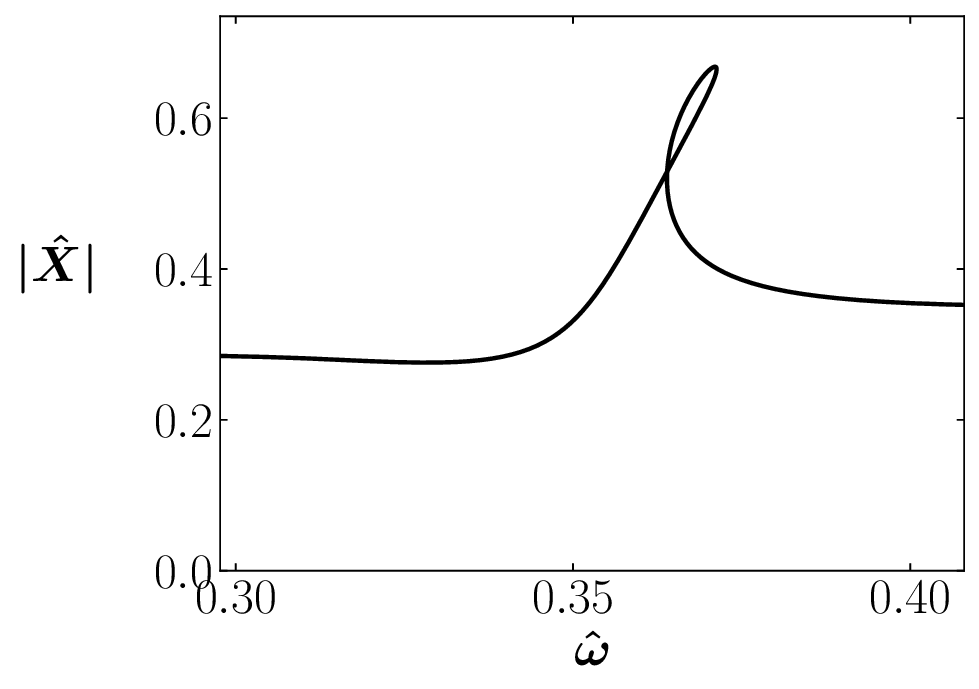}
		\caption{} 
	\end{subfigure}
	\hfill
	\begin{subfigure}{0.32\linewidth}
	\centering
	\includegraphics[width=\linewidth]{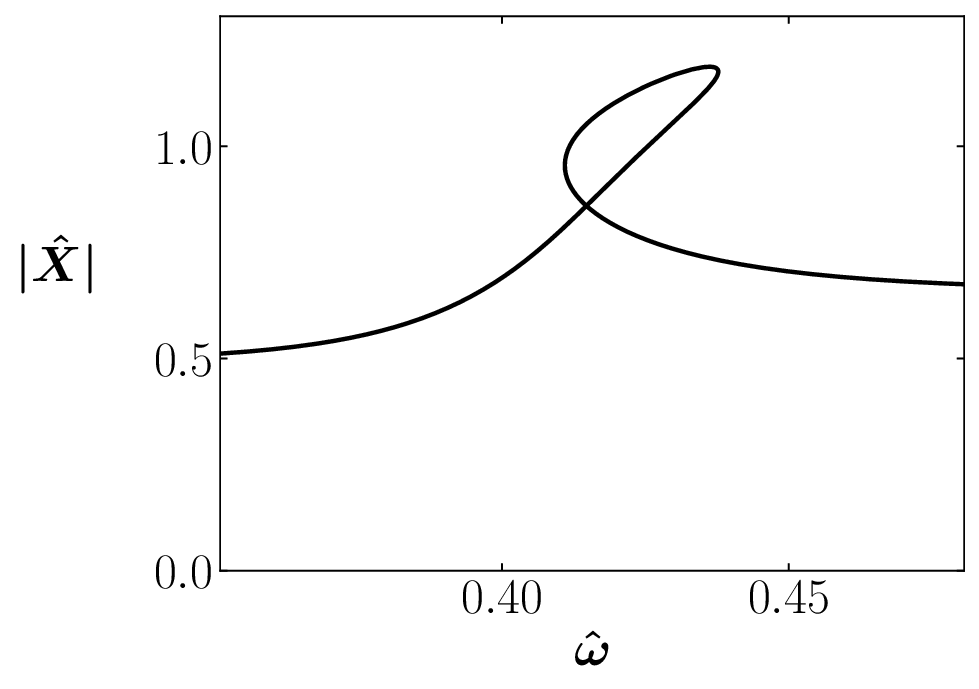}
	\caption{} 
	\end{subfigure}
	\\
		\begin{subfigure}{0.32\linewidth}
		\centering
		\includegraphics[width=\linewidth]{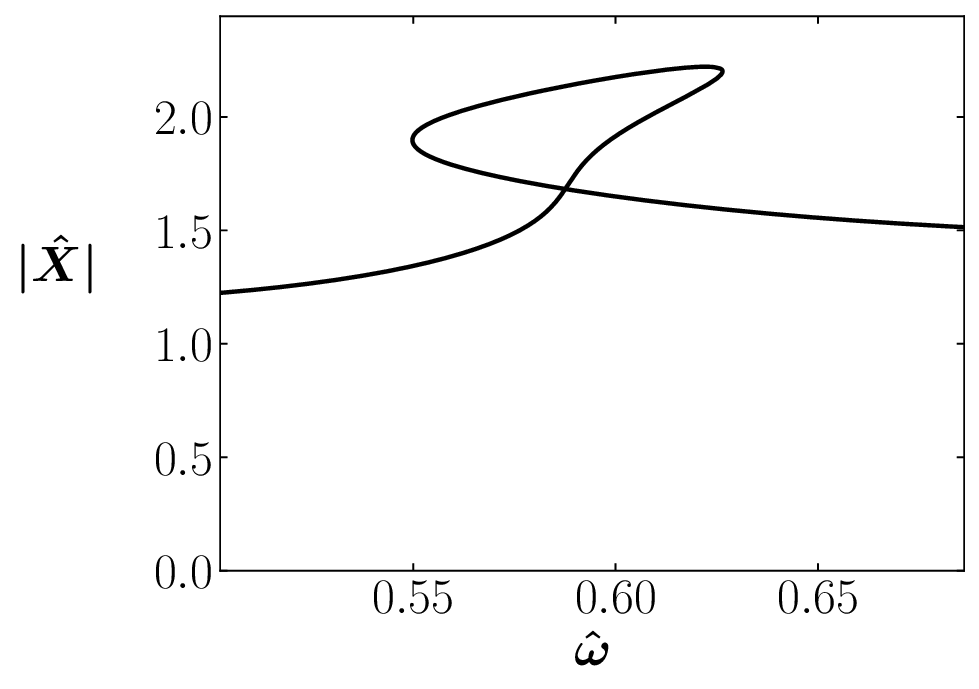}
		\caption{} 
	\end{subfigure}
	\hfill
	\begin{subfigure}{0.32\linewidth}
		\centering
		\includegraphics[width=\linewidth]{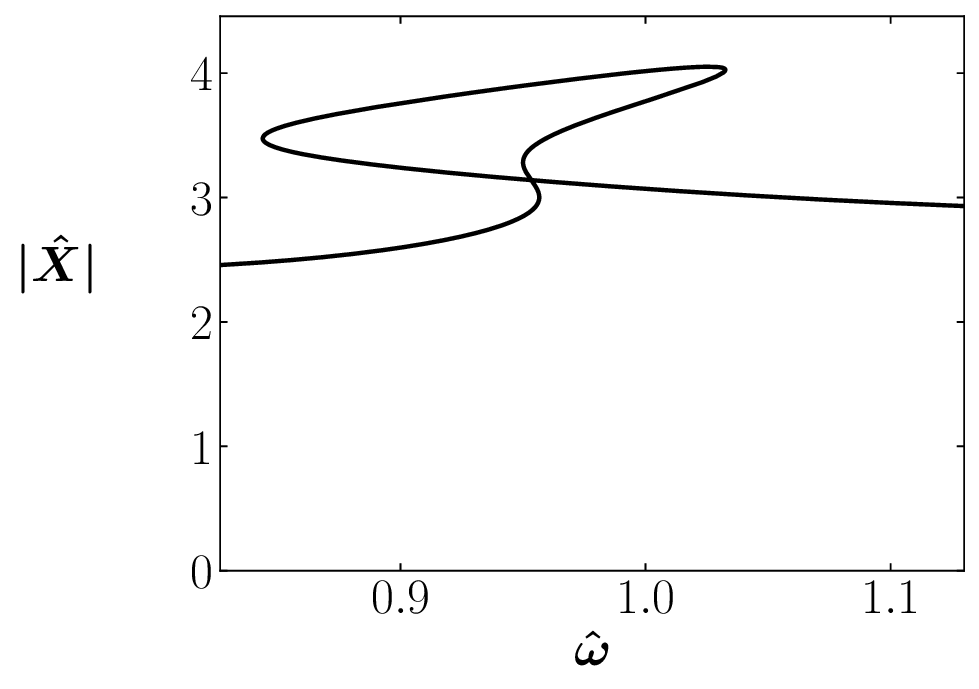}
		\caption{} 
	\end{subfigure}
	\hfill
	\begin{subfigure}{0.32\linewidth}
		\centering
		\includegraphics[width=\linewidth]{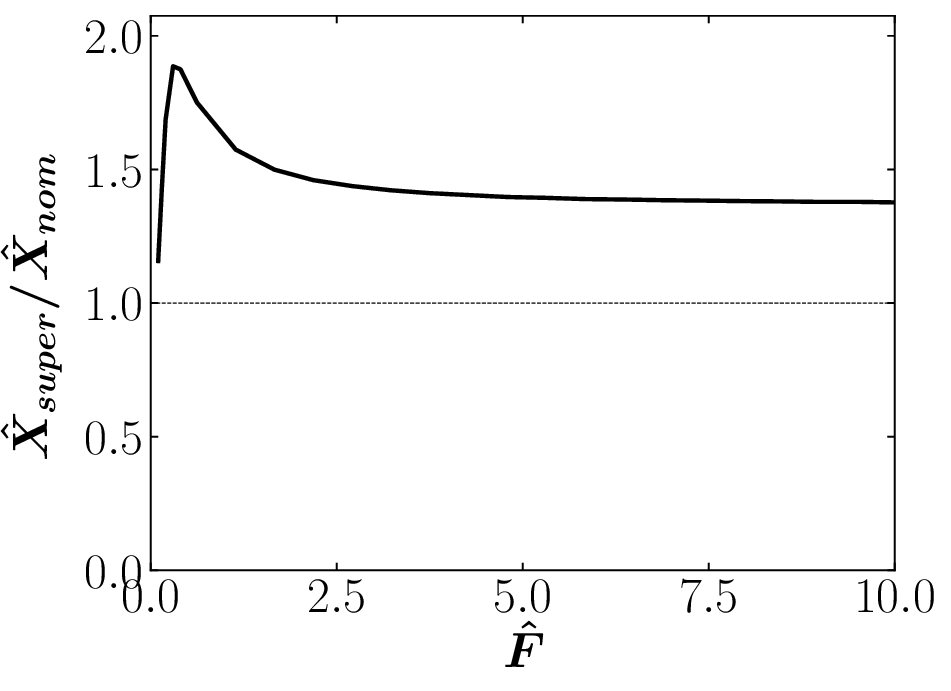}
		\caption{} 
	\end{subfigure}
		\caption{3:1 superharmonic resonances for stiffening Duffing case for nondimensional external force levels of (a) 0.10, (b) 0.30, (c) 0.62, (d) 2.18, and (e) 10.0. Figure (f) shows the magnitude of superharmonic resonance relative to a nominal response for a range of force levels. The nominal response amplitude $\hat{X}_{nom}$ for (f) is the nondimensional response amplitude at 1.1 times the frequency of the maximum third harmonic response, which is also the maximum frequency plotted in (a)-(e). $\hat{X}_{super}$ for (f) is the total nondimensional response amplitude at the peak amplitude of the third harmonic.} 
	\label{fig:example_duffing_superharmonic}
\end{figure}

\begin{figure}[!h]
	\centering
	\begin{subfigure}{0.32\linewidth}
		\centering
		\includegraphics[width=\linewidth]{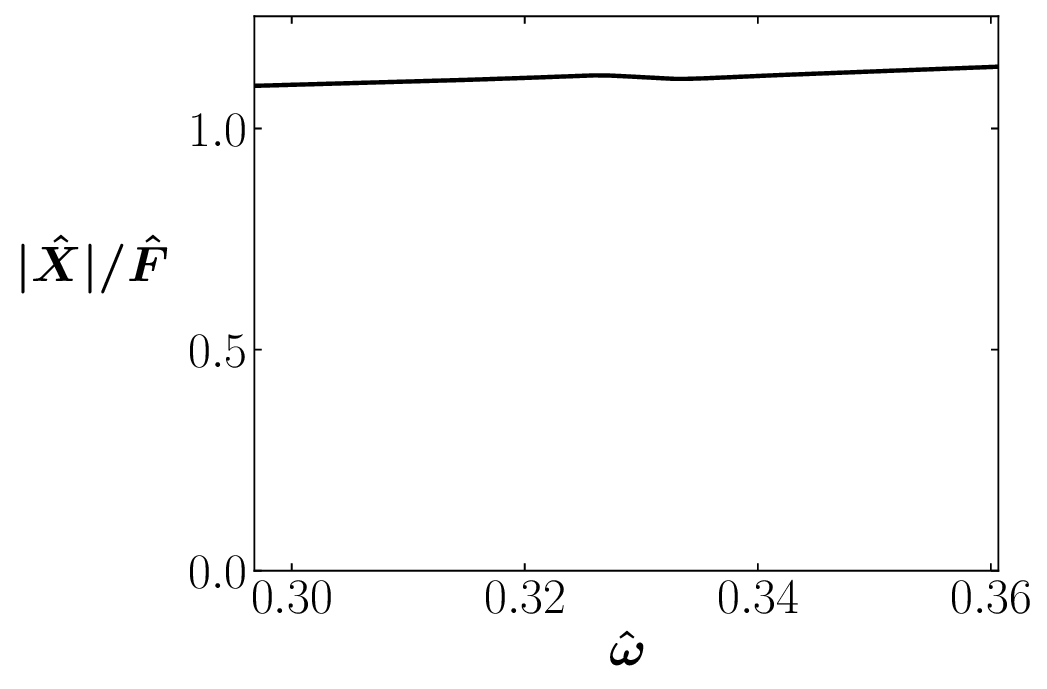}
		\caption{} 
	\end{subfigure}
	\hfill
	\begin{subfigure}{0.32\linewidth}
		\centering
		\includegraphics[width=\linewidth]{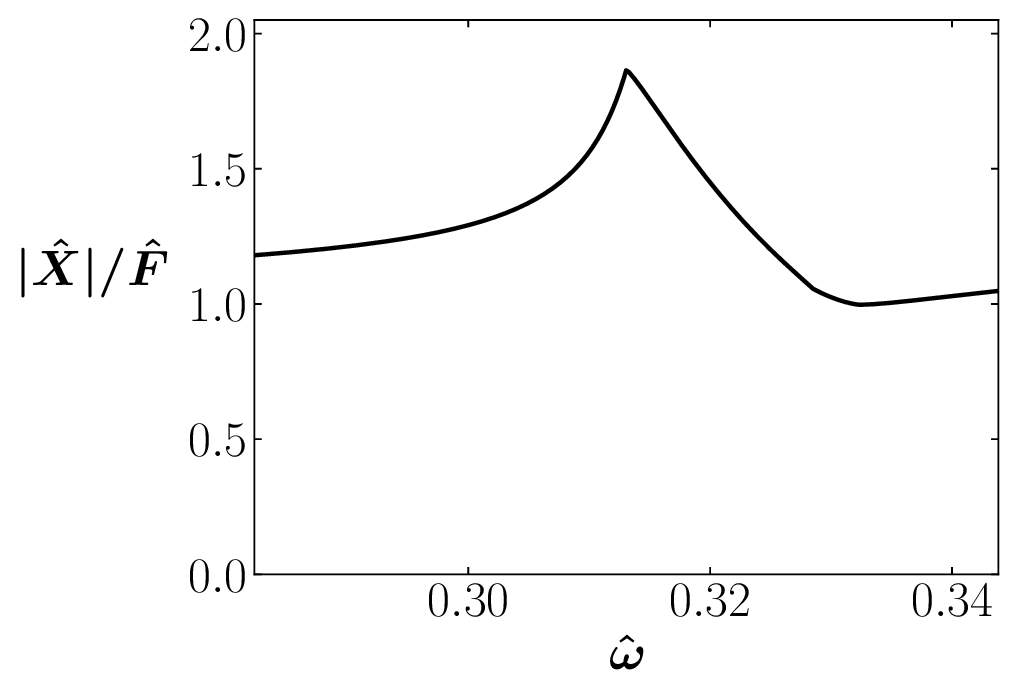}
		\caption{} 
	\end{subfigure}
	\hfill
	\begin{subfigure}{0.32\linewidth}
		\centering
		\includegraphics[width=\linewidth]{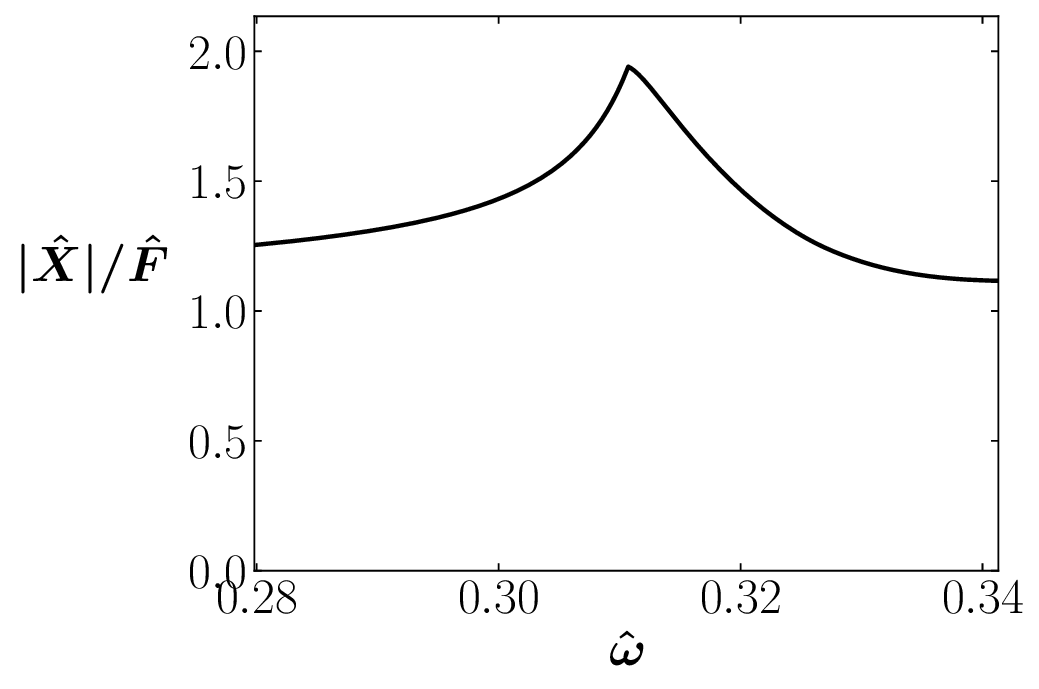}
		\caption{} 
	\end{subfigure}
	\\
		\begin{subfigure}{0.32\linewidth}
		\centering
		\includegraphics[width=\linewidth]{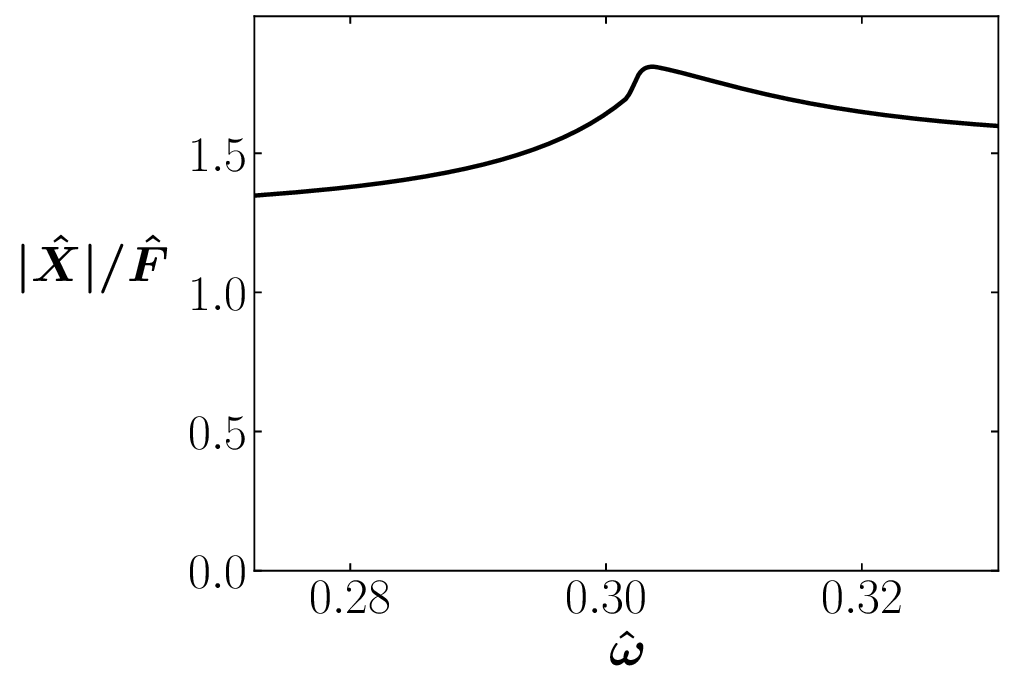}
		\caption{} 
	\end{subfigure}
	\hfill
	\begin{subfigure}{0.32\linewidth}
		\centering
		\includegraphics[width=\linewidth]{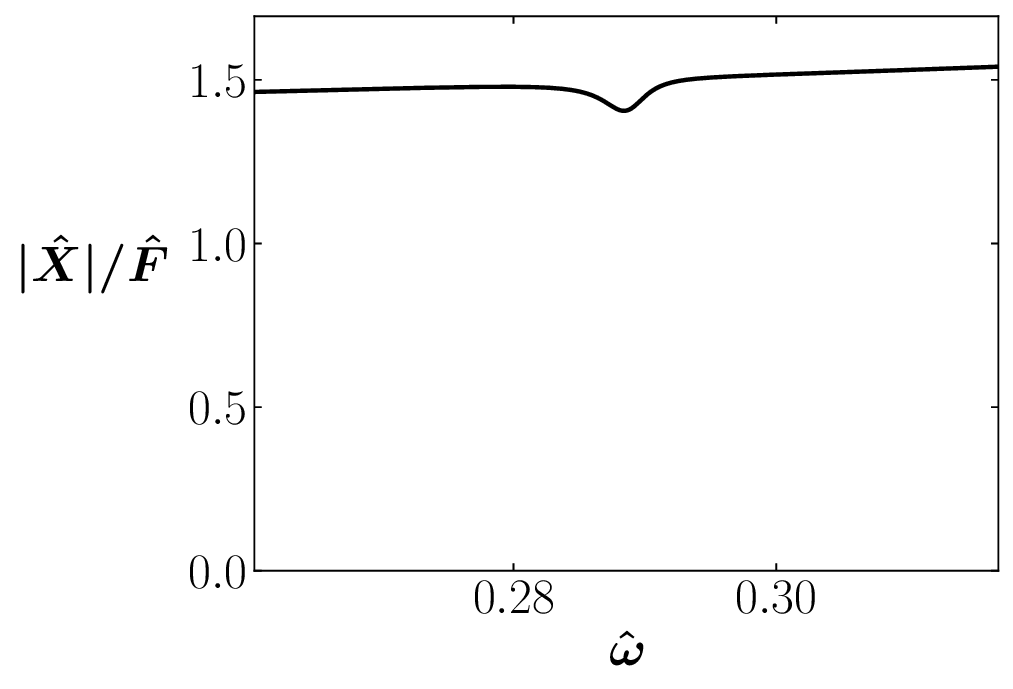}
		\caption{} 
	\end{subfigure}
	\hfill
	\begin{subfigure}{0.32\linewidth}
		\centering
		\includegraphics[width=\linewidth]{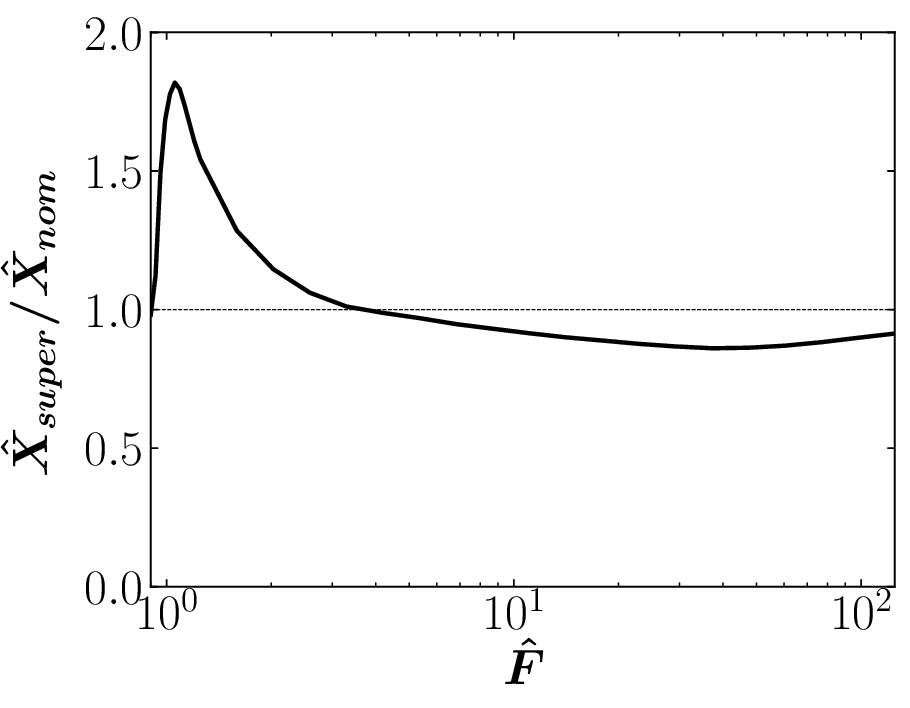}
		\caption{} 
	\end{subfigure}
		\caption{3:1 superharmonic resonances for Jenkins element case for nondimensional external force levels of (a) 0.9, (b) 1.02, (c) 1.13, (d) 2.59, (e) 125. Figure (f) shows the magnitude of superharmonic resonance relative to a nominal response for a range of force levels. The nominal response amplitude $\hat{X}_{nom}$ for (f) is the nondimensional response amplitude at 1.1 times the frequency of the maximum third harmonic response, which is also the maximum frequency plotted in (a)-(e). $\hat{X}_{super}$ for (f) is the total nondimensional response amplitude at the peak amplitude of the third harmonic.} 
		\label{fig:example_jenkins_superharmonic}
\end{figure}

The superharmonic resonances shown in \Cref{fig:example_FRC} and \Cref{sec:appendix_frcs} are only for a single force level, but the superharmonic resonances evolve over a range of force levels. Figures \ref{fig:example_duffing_superharmonic} and \ref{fig:example_jenkins_superharmonic} shows evolution for the 3:1 superharmonic resonance for the stiffening Duffing and Jenkins element cases respectively.  
For both cases, very low force levels produce nearly linear responses without notable superharmonic resonances. For the stiffening Duffing case, the superharmonic resonance is prominent for all forcing levels above a threshold (see \Cref{fig:example_duffing_superharmonic}(f)). Contrarily, the Jenkins element produces notable superharmonic resonances over a limited range of excitation amplitudes with the response at higher forcing amplitudes approaching a linear case again.
Additionally, the superharmonic resonance for the Jenkins element results in a local minimum at high forcing amplitudes because of the phase difference between the first and third harmonics.
These behaviors of superharmonic resonances are detailed further and the present curves are revisited in \Cref{sec:results}.

\FloatBarrier

\section{Modeling Superharmonic Resonances} \label{sec:modeling}

This section first summarizes an existing method \cite{volvertPhaseResonanceNonlinear2021, volvertResonantPhaseLags2022, abeloosControlbasedMethodsIdentification2022a} for tracking superharmonic resonances (\Cref{sec:previous_prnm}). 
However, the existing method cannot easily be applied to the most general cases of nonlinearities. Therefore, a new approach is derived by decomposing the nonlinear forces in the system in \Cref{sec:decomposeFnl}.
Then, \Cref{sec:apriori_phase} uses the decomposed nonlinear forces to predict the phase of the superharmonic resonances.
Finally, the new approach for tracking superharmonic resonances is formally defined in \Cref{sec:tracking_method}.

\subsection{Existing Phase Resonance Nonlinear Modes} \label{sec:previous_prnm}

Prior derivations of PRNM provided phase criteria for superharmonic, subharmonic, and ultra-subharmonic resonances for the Duffing oscillator \cite{volvertPhaseResonanceNonlinear2021, volvertResonantPhaseLags2022}. This phase criteria takes the form of a constant phase lag between the superharmonic response and either the external forcing or a response at the forcing frequency.
The phase lag allows for the tracking of superharmonic resonances without calculating full FRCs providing the opportunity to understand the system behavior over a range of force levels at significantly lower computational cost \cite{volvertPhaseResonanceNonlinear2021}.
For 3:1 superharmonic resonances, a phase lag of $ \pi/2 $ or $ 3 \pi/2 $ for the stiffening or softening Duffing oscillators respectively was derived \cite{volvertResonantPhaseLags2022}. In addition a phase lag of $ \pi/2 $ is reported for a 5:1 superharmonic resonance for the stiffening Duffing case \cite{volvertPhaseResonanceNonlinear2021, volvertResonantPhaseLags2022}. These studies predominately focused on SDOF systems and provided constant phase lags between the harmonic forcing of the fundamental frequency and the phase of the superharmonic of interest \cite{volvertPhaseResonanceNonlinear2021, volvertResonantPhaseLags2022}. 
The present work extends PRNM by considering amplitude dependent phase lags that are essential for characterizing hysteretic nonlinearities (e.g., see \Cref{sec:results_hyst}).

For an alternative formulation of PRNM, the phase criteria has also been derived using a perturbation method for odd and even nonlinearities while allowing for MDOF systems 
\cite{abeloosControlbasedMethodsIdentification2022a}. For a first order $ n$:1 superharmonic resonance (e.g., 3:1 for the Duffing oscillator), the phase criteria is \cite{abeloosControlbasedMethodsIdentification2022a}
\begin{equation}\label{eq:perturb_phase}
	\Delta_n = \phi_{1,n} - n \phi_{0,1}
\end{equation}
where $ \phi_{0,1} $ is the phase of the first harmonic from the order $ \varepsilon^0 $ solution and $ \phi_{1,n} $ is the phase of the $ n $th harmonic from the order $ \varepsilon^1 $ solution. Here $ \varepsilon \ll 1 $ is a small parameter that separates the scales of the weak nonlinearity (order $ \varepsilon^1 $) from the linear system (order $ \varepsilon^0 $) as described by the equation of motion \cite{abeloosControlbasedMethodsIdentification2022a}: 
\begin{equation}
	\ddot{x} + 2 \omega_0 \zeta_0 \dot{x} + \omega_0^2 x + \dfrac{\varepsilon}{m} f_{nl}(x, \dot{x}) = F \cos \omega t.
\end{equation}
The parameter $ \varepsilon $ serves as the basis of perturbation methods (e.g., \cite{nayfehNonlinearInteractionsAnalytical2000, nayfehNonlinearOscillations1995}). It is important to note that $ \phi_{0,1} $ represents the phase of the underlying linear system since the nonlinearity is order $ \varepsilon^1 $ and thus not included in the order $ \varepsilon^0 $ solution. More generally, it is possible that $ \phi_{0,1} $ is not equal to the phase of the first harmonic in the full perturbation solution. 
This means that the phase criteria of \eqref{eq:perturb_phase} cannot be easily applied to HBM type solutions or experiments in cases where the phase of the underlying linear system is not physically realized and is not known.\footnote{For some cases, it may be reasonable to utilize the physically realized phase of the first harmonic as $\phi_{0,1}$ such as the experiments demonstrated in \cite{abeloosControlbasedMethodsIdentification2022a}.}
The present extension of PRNM is formulated to allow for direct use with HBM without requiring any analytical solution steps and is more general than the formulation of PRNM  \cite{volvertPhaseResonanceNonlinear2021, volvertResonantPhaseLags2022, abeloosControlbasedMethodsIdentification2022a}. Furthermore, the present work enables the easy calculation of the phase criteria numerically rather than analytically via a perturbation method.

For 5:1 superharmonic resonances of a Duffing oscillator, a second order approach is necessary, providing a phase criteria of \cite{abeloosControlbasedMethodsIdentification2022a}
\begin{equation} \label{eq:secondorder_perturb}
	\Delta_5 = \phi_{2,5} - \phi_{1,3} - 2 \phi_{0,1}
\end{equation}
where $\phi_{2,5}$ is the phase of the fifth harmonic from the $\varepsilon^2$ solution and $\phi_{1,3}$ is the phase of the third harmonic from the $\varepsilon^1$ solution.
Based on the perturbation analysis, Section 5.2 of \cite{abeloosControlbasedMethodsIdentification2022a} conjectures that the phases described by \eqref{eq:perturb_phase} and \eqref{eq:secondorder_perturb} are $ -\pi/2 $ or 0 for the cases of odd or even nonlinearities respectively. For the case of the SDOF stiffening Duffing oscillator, this analysis is consistent with a phase lag of $ \pi/2 $ since the fundamental harmonic is away from resonance and thus responds near 0 phase \cite{volvertPhaseResonanceNonlinear2021, volvertResonantPhaseLags2022}. The logical arguments for the conjecture in \cite{abeloosControlbasedMethodsIdentification2022a} appear to hold for more general criteria of $ \pm \pi/2 $ or $ 0 $ and $ \pi $ for the cases of odd or even nonlinearities respectively. If the conjecture is generalized, the case of the softening Duffing oscillator from \cite{volvertResonantPhaseLags2022} is also consistent. 
This conjecture is analyzed for the examples in the present work, and the unilateral spring is found to violate the conjecture for even nonlinearities.

Previous derivations addressed phase lags for analytically described nonlinear forces, but the analysis techniques rely on perturbation methods that cannot be applied to numerically described nonlinearities (e.g., for hysteretic models). 
As an alternative, the present work provides a framework for applying a superharmonic phase criteria for a system with any arbitrary (including numerically described) nonlinearity (see \Cref{sec:tracking_method}). As shown in \Cref{sec:results}, the phase criteria is not required to be calculated a priori, but rather can be done during HBM type calculations allowing for a wider range of applications than can easily be analyzed with averaging and perturbation methods.
Furthermore, \Cref{sec:apriori_phase} illustrates that the phase criteria for hysteretic nonlinearities cannot be described as a constant value and thus prevents a simple analytically derived phase lag from being applied.

\subsection{Decomposing Nonlinear Forces} \label{sec:decomposeFnl}

The present work uses HBM to calculate the periodic steady-state responses to the nonlinear vibration problem and as a reference solution (see \Cref{sec:hbm})  \cite{krackHarmonicBalanceNonlinear2019}. 
For the present work with a SDOF system, the steady state motion is assumed to be
\begin{equation} \label{eq:x_approx}
	x(t) 
	= 
	X_0 + \sum_{k=1}^H  \left[ X_{kc} \cos(k \omega t) + X_{ks} \sin(k \omega t) \right]
	=
	X_0 + \sum_{k=1}^H X_{k} \cos(k \omega t - \phi_k)
\end{equation}
where $ H $ is the highest harmonic included in the approximation, $ \omega $ is the forcing frequency, $ X_0 $ is the zeroth harmonic displacement, and $ X_{kc} $ and $ X_{ks} $ are the harmonic displacements for cosine and sine respectively for the $ k $th harmonic. Alternatively, the total amplitude of harmonic $ k $ is $ X_k $ and has phase $ \phi_k $. For later convenience, the time series of displacements associated with the $ n $th harmonic is defined as
\begin{equation}
	x_n(t) = X_{nc} \cos(n \omega t) + X_{ns} \sin(n \omega t)
\end{equation}
and the time series associated with harmonics $ 0 $ through $ j $ is denoted as
\begin{equation}
	x_{0:j}(t) = X_0 + \sum_{k=1}^j X_{kc} \cos(k \omega t) + X_{ks} \sin(k \omega t) 
\end{equation}

The equations from HBM for the $ n $th harmonic are analyzed to understand occurrences of $ n$:1 superharmonic resonances. In general, the harmonic coefficients for the nonlinear forces acting on the $ n $th harmonic cosine and sine equations, $ F_{nl,nc} $ and $ F_{nl, ns} $ respectively, are calculated as
\begin{subequations} \label{eq:fourier_int}
\begin{equation}
	F_{nl,nc} 
	=
	\mathcal{F}_{nc} \{ f_{nl}[x(t)] \} 
	=
	\dfrac{\omega}{\pi} \int_0^{2\pi/\omega} f_{nl}[x(t)] \cos(n \omega t) dt
\end{equation}
\begin{equation}
	F_{nl,ns} 
	=
	\mathcal{F}_{ns} \{ f_{nl}[x(t)] \} 
	= 
	\dfrac{\omega}{\pi} \int_0^{2\pi/\omega} f_{nl}[x(t)] \sin(n \omega t) dt.
\end{equation}
\end{subequations}
Here, $ \mathcal{F}\{ \cdot \} $ denotes a Fourier transform with the subscripts denoting the harmonic number ($ n $) and cosine or sine ($ c $ or $ s $ respectively). Throughout this work, nonlinear force harmonic coefficients are evaluated with the Alternating Frequency-Time (AFT) method. AFT calculates the displacements in the time domain over a cycle and then evaluates the nonlinear forces in the time domain. Finally, the nonlinear forces are converted back into the frequency domain (see \Cref{sec:hbm} for more details and discussion about computational improvements to AFT for hysteretic models).

To decompose the nonlinear forces on the system, the broadband excitation of the $ n $th harmonic from the motion of the lower harmonics is defined as 
\begin{equation} \label{eq:fbroad}
	F_{nq,broad} = - \mathcal{F}_{nq} \{ f_{nl} [ x_{0:(n-1)}(t) ] \} ,
\end{equation}
where $ q=c,s $ denotes cosine or sine respectively.
This represents an excitation of the $ n $th harmonic from the presence of lower harmonics even if the $ n $th harmonic has no motion. This excitation is hypothesized to be a main cause of superharmonic resonances. 
Next, since the nonlinear forces do not satisfy superposition, a correction term for the violations of superposition is defined as 
\begin{equation} \label{eq:fsup}
	F_{kq,sup,n} = -\mathcal{F}_{kq} \{ f_{nl}[x(t)] - f_{nl}[x_n(t)] - f_{nl}[x_{0:(n-1)}(t)] \}.
\end{equation}
This force, $ F_{kq,sup,n} $, denotes the effect on the $ k $th harmonic of introducing the $ n $th harmonic to the solution.
In this section, only the first $ n $ harmonics are considered to analyze the $ n $th superharmonic resonance.\footnote{Note that only the first $n-1$ harmonics are used in $F_{nq,broad}$, but the $n$th harmonic is included in $F_{kq,sup,n}$.} 
Using these definitions, the nonlinear force on the $ n $th harmonic can be decomposed as\footnote{This decomposition is exact when $n$ harmonics are used in HBM. The decomposition is used to rearrange portions of the nonlinear force to better understand the dynamics and is exploited throughout the remainder of the paper.}
\begin{equation}
	F_{nl,nq} 
	=
	\mathcal{F}_{nq} \{ f_{nl}[x(t)] \} 
	= 
	\mathcal{F}_{nq} \{ f_{nl}[x_n(t)] \}
	- F_{nq,broad}
	- F_{nq,sup,n}.
\end{equation}
This leads to the equations of motion of the $ n $th harmonic from HBM of
\begin{equation} \label{eq:force_decomp_hbm}
\begin{split}
	(-n \omega^2 m + k) X_{nc} + (n \omega c) X_{ns} +  \mathcal{F}_{nc} \{ f_{nl}[x_n(t)] \}  
	&= 
	F_{nc,broad} + F_{nc,sup,n}
		\\
	(-n \omega^2 m + k) X_{ns} - (n \omega c) X_{nc} + \mathcal{F}_{ns} \{ f_{nl}[x_n(t)] \}  
	&= 
	F_{ns,broad} + F_{ns,sup,n} .
\end{split}
\end{equation}
Here, the left hand side of the equations is of the form of single harmonic motion for the nonlinear system, while the right hand side denotes what is treated as external forcing on the harmonic. 
Away from the superharmonic resonance, it is expected that the amplitude of the $ n $th harmonic will be small and the terms $ F_{nq,sup,n} $ will be small. More broadly, the present work assumes $F_{kq,sup,n}$ remains small for all harmonics $k$ and that the higher harmonic terms do not significantly influence the lower harmonic terms. The validity of this assumption is assessed empirically with the examples in \Cref{sec:results}. At frequencies below the superharmonic resonances, the nonlinear vibration is expected to occur in phase with the broadband excitation of $ F_{nq,broad} $. As the system passes the superharmonic resonance, the phase $ \phi_n $ is expected to increase by $ \pi $ to be out of phase with $ F_{nq,broad} $. Consistent with nonlinear phase resonance of the primary harmonic \cite{peetersDynamicTestingNonlinear2011, rensonForceAppropriationNonlinear2018, scheelPhaseResonanceApproach2018}, it is expected that the resonance of the $ n $th harmonic will occur when $ \phi_n $ has phase near
\begin{equation} \label{eq:simple_phase_res}
	\phi_{n} = arctan2(F_{ns,broad}, F_{nc,broad}) + \dfrac{\pi}{2}
\end{equation}
where $ arctan2 $ is the 4-quadrant arctangent operator.
The phase angle of the broadband excitation of a higher harmonic $ n $ is defined to be
\begin{equation} \label{eq:phibroad}
	\phi_{broad,n} = arctan2(F_{ns,broad}, F_{nc,broad}).
\end{equation}

\section{A Priori Phase Calculations}\label{sec:apriori_phase}

The present section provides a preliminary understanding of the superharmonic resonances for the considered nonlinear forces. A more in-depth exploration utilizing a new proposed method is reserved for \Cref{sec:results}.
Prior to simulating the dynamics of the system, the phase of the $ n $th superharmonic can be predicted for some simplified cases. Specifically, this section analyzes the expected phase \eqref{eq:simple_phase_res} of the superharmonic resonances that are excited by fundamental harmonic motion. Since the $ n$:1 superharmonic resonance occurs near the natural frequency of the system divided by $ n $, it is expected that the fundamental harmonic will respond nearly in phase with the fundamental forcing. This assumption can be violated by the fact that the $ n $th harmonic violates superposition resulting in the forces $ F_{1q,sup,n} $. 
This assumption also breaks down for systems with high damping due to the shift in phase of the fundamental response further from the primary resonance.
However, for the present analysis, these effects will be neglected. 
This analysis based on \eqref{eq:fourier_int}, \eqref{eq:fbroad}, and \eqref{eq:simple_phase_res}
is similar to the application of a perturbation method and gives phase criteria similar to those derived in \cite{abeloosControlbasedMethodsIdentification2022a} when the same nonlinear forces are considered.
The present section also provides insight into superharmonic responses of hysteretic systems based on a simple analysis of the interaction forces between harmonics described in \Cref{sec:decomposeFnl}.
Next, the primary cases of the lowest observed integers $ n $ for $ n $:1 superharmonic resonances are considered.\footnote{Here, the lowest observed superharmonic resonance is the lowest integer $n \geq 2$ that has nonzero $F_{nc,broad}$ or $F_{ns,broad}$ from \eqref{eq:fbroad} for a given system.}
Secondary superharmonic resonances corresponding to larger values of $n$ are similarly considered in \Cref{sec:apriori_secondary}.

\subsection{Primary Superharmonics}

For each of the nonlinear forces presented in \Cref{sec:nlforces}, the quantities $ F_{nc,broad} $ and $ F_{ns,broad} $ are calculated with \eqref{eq:fbroad} and either the analytical integral of \eqref{eq:fourier_int} or the AFT method (see \Cref{sec:aft}) and summarized in \Cref{tab:Fbroad_primary} under the assumption that the fundamental harmonic motion is
\begin{equation} \label{eq:fundamentalMotion}
	x(t) = X_{1} \cos(\omega t).
\end{equation}
To allow for analytical calculations and give conceptual insight, $X_0$ is not considered here.\footnote{The full solutions including $X_0$ are numerically calculated in \Cref{sec:results}.}
The present section intends to provide an analytical understanding of the nonlinear forces rather than exact predictions of the responses. The accuracy of the predicted phases ($\phi_n$ and $\phi_{broad,n}$, calculated with \eqref{eq:simple_phase_res} and \eqref{eq:phibroad} respectively) discussed here are analyzed against FRCs in \Cref{sec:results}.

\FloatBarrier

\begin{table}[h!]
	\centering
	\caption{Analytical calculations of $ F_{nc,broad} $ and $ F_{ns,broad} $ (see \eqref{eq:fbroad}) excitation of primary superharmonic resonances. The equations and values presented here are independent of the parameters chosen in \Cref{tab:nl_params}.}
	\label{tab:Fbroad_primary}
		\begin{tabular}{cccccc}
			\hline \hline
			Nonlinear Force & \multicolumn{1}{c}{\begin{tabular}[c]{@{}c@{}}Excited\\ Harmonic\end{tabular}} 
			& $ F_{nc,broad} $ [N] & $ F_{ns,broad} [N] $ & $ \phi_{broad,n} $ [rad]  & $ \phi_n $ [rad]
			\\ \hline
			
			Stiffening Duffing & 3 & $ -\alpha X_1^3/4 $ & 0 & $ -\pi $ & $ -\pi/2 $
			\\
			Quintic Stiffness & 3 & $ -5 \eta X_1^5/16 $ & 0 & $ -\pi $ & $ -\pi/2 $
			\\ \hline
			
			Softening Duffing & 3 & $ -\alpha X_1^3/4 $ & 0 & $ 0 $ & $ \pi/2 $
			\\
			Conservative Softening II & 3 & \multicolumn{2}{c}{\Cref{fig:fbroad_mag}} & 0 & $ \pi/2 $
			\\ \hline
			
			Unilateral Spring & 2 & $ -2 k_{nl} X_1 / 3 \pi $  & 0 & $ -\pi $ & $ -\pi/2 $
						\\ \hline

			Cubic Damping & 3 & 0 & $ -\gamma \omega^3 X_1^3 / 4 $ & $ -\pi/2 $ & $ 0 $ 
			\\
			Jenkins Element & 3 & \multicolumn{2}{c}{\Cref{fig:fbroad_mag}}  & \Cref{fig:fbroad_phase} & Variable
			\\
			Iwan Element & 3 & \multicolumn{2}{c}{\Cref{fig:fbroad_mag}} & \Cref{fig:fbroad_phase}  & Variable

			\\ \hline\hline
		\end{tabular}
\end{table}

\begin{figure}[!h]
	\centering
	\newcommand{\fbroadwidth}{0.49}  
	\begin{subfigure}{\fbroadwidth\linewidth}
		\centering
		\includegraphics[width=\linewidth]{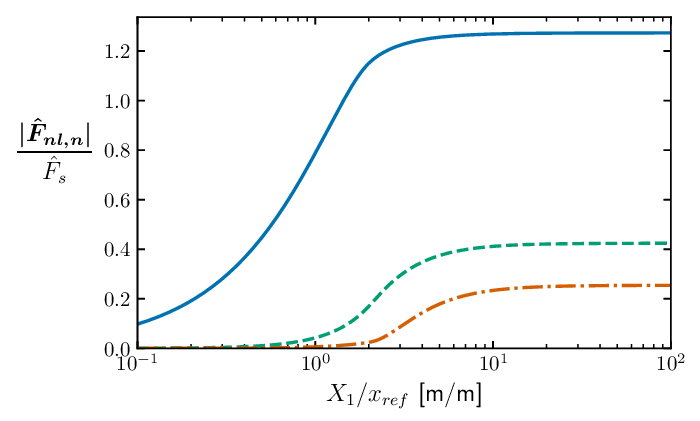}
		\caption{}
	\end{subfigure}
	\quad \quad
	\begin{subfigure}{0.2\linewidth}
		\centering
		\includegraphics[width=\linewidth]{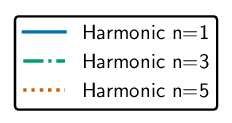}
		\vspace{0.75 in}
	\end{subfigure}
	\\
	\begin{subfigure}{\fbroadwidth\linewidth}
		\centering
		\includegraphics[width=\linewidth]{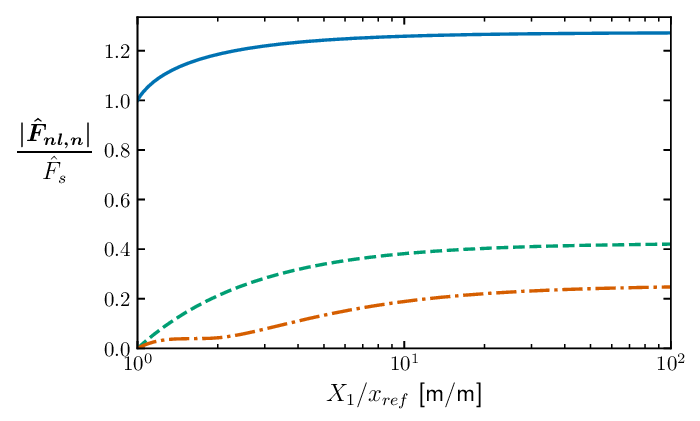}
		\caption{}
	\end{subfigure}
	\hfill
	\begin{subfigure}{\fbroadwidth\linewidth}
		\centering
		\includegraphics[width=\linewidth]{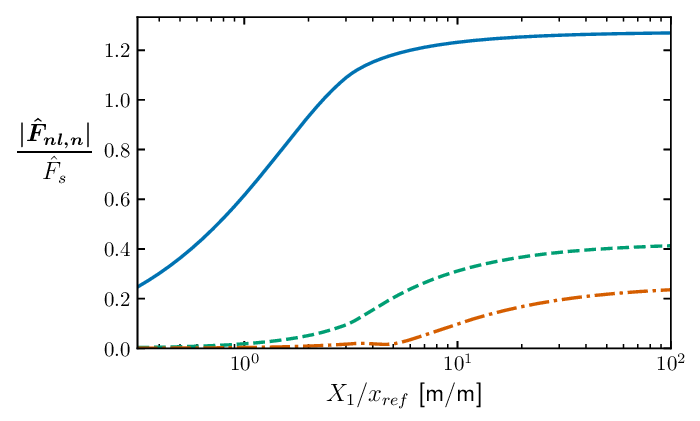}
		\caption{}
	\end{subfigure}
		\caption{Magnitudes of nonlinear force coefficients for fundamental motion of $ x(t) = X_1 \cos(\omega t ) $ for (a) II softening nonlinearity, (b) Jenkins element, and (c) Iwan element. Plots are normalized by $ x_{ref} $, the displacement where the force saturates, and $ \hat{F}_s $, the nondimensional saturated force value.
	Here the solutions are numerically calculated for the parameters in \Cref{tab:nl_params}. The normalization removes any dependency on the values of $k_t$ and $F_s$. For the II softening nonlinearity and the Iwan element, increasing $\beta$ causes the behavior to become more similar to that of the Jenkins element. For those nonlinearities, decreasing $\chi$ decreased the force values over the plotted region.
	The magnitudes of the nonlinear force coefficients are calculated as $|\boldsymbol{\hat{F}_{nl,n}}| = (F_{nl,nc}^2 + F_{nl,ns}^2) / (k_{lin} x_{ref})$.
	} 
						\label{fig:fbroad_mag}
\end{figure}

\begin{figure}[!h]
	\centering
	\newcommand{\fbroadwidth}{0.49}  
	\begin{subfigure}{\fbroadwidth\linewidth}
		\centering
		\includegraphics[width=\linewidth]{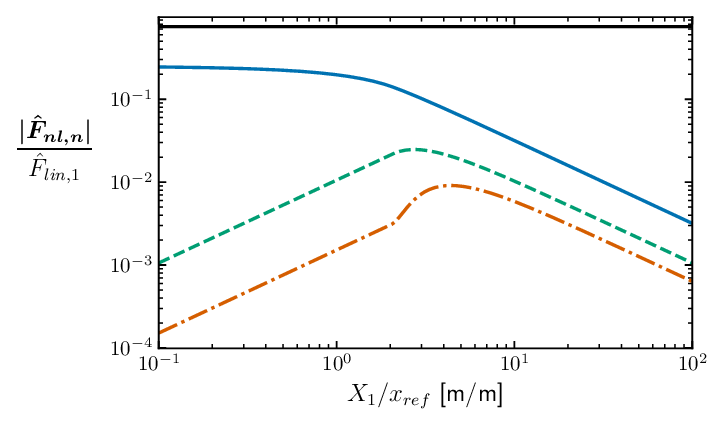}
		\caption{}
	\end{subfigure}
	\quad \quad
	\begin{subfigure}{0.2\linewidth}
		\centering
		\includegraphics[width=\linewidth]{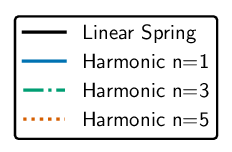}
		\vspace{0.75 in}
	\end{subfigure}
	\\
	\begin{subfigure}{\fbroadwidth\linewidth}
		\centering
		\includegraphics[width=\linewidth]{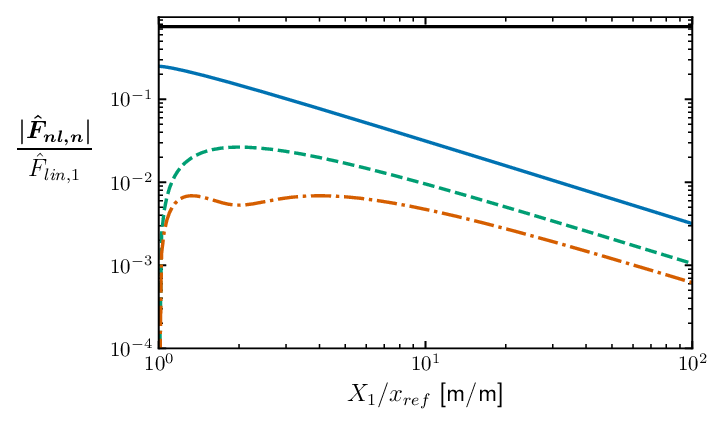}
		\caption{}
	\end{subfigure}
	\hfill
	\begin{subfigure}{\fbroadwidth\linewidth}
		\centering
		\includegraphics[width=\linewidth]{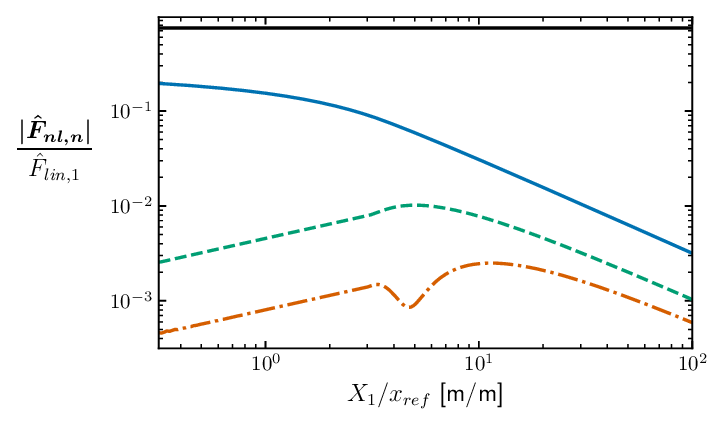}
		\caption{}
	\end{subfigure}
		\caption{Replotted nonlinear harmonic forces from \Cref{fig:fbroad_mag} normalized by $\hat{F}_{lin,1}$, the first harmonic of the nondimensional force produced by a spring with stiffness $k_{lin}$, for (a) II softening nonlinearity, (b) Jenkins element, and (c) Iwan element. The horizontal black lines in each plot represents the contribution of the first harmonic of the linear force due to the linear spring included in each system in \Cref{sec:nlforces}.}
	\label{fig:fbroad_mag_norm}
\end{figure}

For the Duffing, quintic, and cubic damping nonlinearities, the excitation of the third harmonic grows faster than linearly with increasing amplitude (see \Cref{tab:Fbroad_primary}) and thus the importance of the superharmonic resonances can be expected to increase with increasing amplitude. For the unilateral spring, the excitation of the second harmonic is linear in $ X_1 $, and thus the extent of the superharmonic resonance is expected to be constant with amplitude. On the other hand, the conservative softening II nonlinearity, the Jenkins element, and the Iwan element all show saturating excitation of the third harmonic (see \Cref{fig:fbroad_mag}). While these cases cannot be evaluated analytically, the calculations shown in \Cref{fig:fbroad_mag} do not require any solutions of systems of equations, but rather result from just nonlinear force evaluations.
For the three saturating nonlinearities, the range of force levels that result in prominent superharmonic resonances is of interest (see \Cref{sec:frcs_main} and \Cref{fig:example_jenkins_superharmonic}).
To better understand this range, the calculations form \Cref{fig:fbroad_mag} are normalized in \Cref{fig:fbroad_mag_norm} by the force produced by a linearized spring with stiffness $k_{lin}$ for a given displacement value.
For these cases, it is expected that at low amplitudes the superharmonic resonance will be small since the excitation of the third harmonic is low. 
Then at a moderate amplitude level, the superharmonic resonance will be clear. The exact amplitude range where the superharmonic resonance is prominent cannot be determined based on the simple analysis here. However, the superharmonic resonance is most likely to become appreciable for a range of displacements between
the displacement where the force saturates and 10 times that displacement level because the excitation of the higher harmonics grows most quickly there (see \Cref{fig:fbroad_mag}) and thus peaks relative to the contribution of the linear spring (see \Cref{fig:fbroad_mag_norm}). At higher amplitudes, the superharmonic resonance will likely decrease in prominence as the response of the first harmonic continues to grow, but the excitation force for the third harmonic saturates. 
This behavior is previewed in \Cref{sec:frcs_main} and confirmed in \Cref{sec:results}. This analysis illustrates the potential understanding of the system behavior that can be obtained through a simple calculations.

\FloatBarrier

The two polynomial stiffening nonlinearities (stiffening Duffing and Quintic stiffness) both result in excitation of the third harmonic at phase $ -\pi $ and thus are expected to have 3:1 superharmonic resonances at a phase of $ -\pi/2 $ for the third harmonic. These phase criteria are consistent with the analysis of \cite{volvertPhaseResonanceNonlinear2021, volvertResonantPhaseLags2022, abeloosControlbasedMethodsIdentification2022a}, though those studies did not directly consider a quintic stiffness. 
For the softening Duffing nonlinearity and the conservative softening II model, the third harmonic is excited at phase 0 and the expected phase of the 3:1 superharmonic resonance is $ \pi/2 $. This result for the softening Duffing nonlinearity is consistent with the analysis of \cite{volvertResonantPhaseLags2022}. Previous studies have not considered a nonlinearity of the form of the conservative softening II model. 

From \Cref{tab:Fbroad_primary}, the unilateral spring excites the second harmonic at phase $ -\pi $ resulting in an expected resonance phase of $ -\pi/2 $ for the second harmonic. This differs from the results of \cite{volvertPhaseResonanceNonlinear2021} that suggested a phase lag of $ 3\pi/4 $ for a 2:1 superharmonic resonance for the Duffing oscillator. However, that analysis only considered the Duffing nonlinearity, so may not be applicable in this case. In addition, the unilateral spring is an even nonlinearity (see \Cref{sec:unispring_force}), but the conjecture from \cite{abeloosControlbasedMethodsIdentification2022a} is not consistent with this expected phase lag. 
These discrepancies motivate the present work to investigate a range of nonlinear forces to better understand the phase lags of superharmonic resonances.

The cubic damping nonlinearity is expected to have a phase of 0 for the 3:1 superharmonic resonance. Cubic damping cannot be categorized as either an even or an odd function of displacement and thus a general form of the phase lag is not conjectured by \cite{abeloosControlbasedMethodsIdentification2022a}. However, the present approach of analyzing the Fourier coefficients provides a straight forward way of predicting the phase without a full dynamic analysis. The accuracy of the phases predicted here are discussed in \Cref{sec:results}.

\begin{figure}[!h]
	\centering
	\begin{subfigure}{0.4\linewidth}
		\centering
		\includegraphics[width=\linewidth]{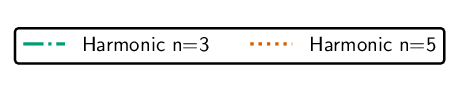}
	\end{subfigure}
	\\
	\newcommand{\fbroadwidth}{0.49}  
	\begin{subfigure}{\fbroadwidth\linewidth}
		\centering
		\includegraphics[width=\linewidth]{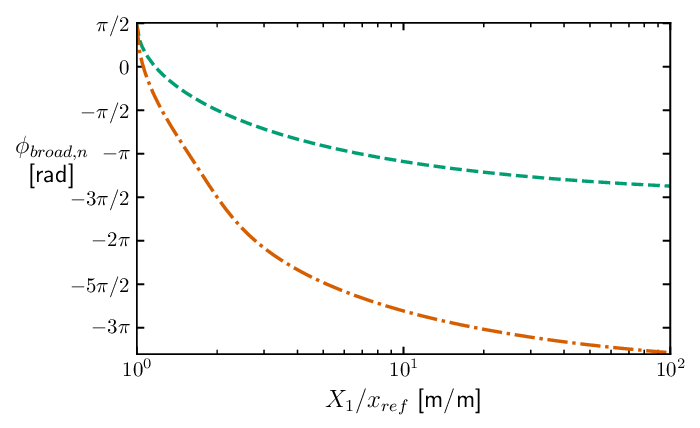}
		\caption{}
	\end{subfigure}
	\hfill
	\begin{subfigure}{\fbroadwidth\linewidth}
		\centering
		\includegraphics[width=\linewidth]{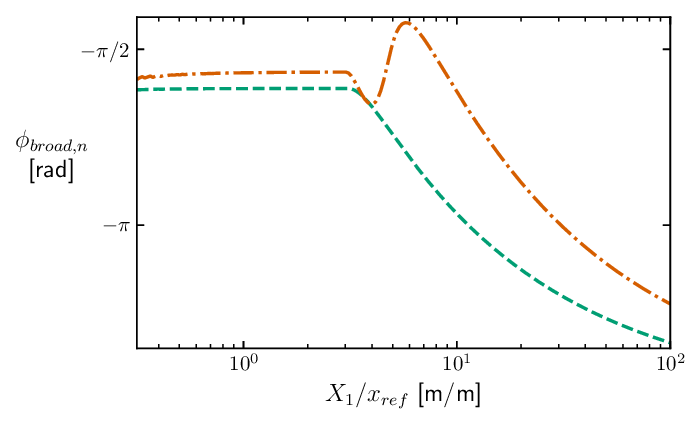}
		\caption{}
	\end{subfigure}
		\caption{Phase of broadband excitation of higher harmonics for fundamental motion of $ x(t) = X_1 \cos(\omega t ) $ for (a) Jenkins element and (b) Iwan element. Plots are normalized by $ x_{ref} $, the displacement where the force saturates. The normalization removes any dependency on the specific values of $k_t$ and $F_s$. Increasing $\beta$ for the Iwan element results in behavior similar to the Jenkins element. Increasing $\chi$ for the Iwan element increases the magnitude of the phase shift near $X_1/x_{ref} = 10 $.} 
	\label{fig:fbroad_phase}
\end{figure}

Finally, the phase characteristics of the two hysteretic models can be analyzed numerically. 
The phase of the excitation of the third harmonic is variable for the Jenkins and Iwan elements as presented in \Cref{fig:fbroad_phase}. 
Any constant phase criteria for the Jenkins and Iwan elements would result in significant errors since the phase varies dramatically over the region plotted in \Cref{fig:fbroad_phase}.
This motivates the tracking method presented in \Cref{sec:tracking_method} as a way to handle the most general nonlinear forces that cannot be analytically analyzed. However, as the amplitude approaches infinity, both hysteretic models converge to square waves of the velocity and towards a constant phase. As discussed in \Cref{sec:results}, the superharmonic resonances are most prominent over the displacement ranges shown in \Cref{fig:fbroad_phase}. Therefore, the limit of a square wave for the nonlinear forces is not informative.  
An approach that can address hysteretic models such as the Jenkins and Iwan models is necessary for assembled structures that are generally modeled with hysteretic nonlinear forces to capture the frictional interactions in bolted joints \cite{brakeMechanicsJointedStructures2017, mathisReviewDampingModels2020}.

The present analysis of the broadband excitation of higher harmonics suggests that the lowest superharmonic resonance for the Duffing oscillator is the 3:1 case. However, 2:1 superharmonic resonance have been documented for the Duffing oscillator as part of branches resulting from symmetry-breaking bifurcations \cite{volvertPhaseResonanceNonlinear2021}.
In addition, similar 2:1 superharmonic resonances may exist for hysteretic nonlinearities \cite{okuizumiMultipleTimeScale2004}.
It is hypothesized that these branches are a result of the forces that describe superposition violations of \eqref{eq:fsup} and not the broadband excitation forces of \eqref{eq:fbroad}. Therefore, these cases are not considered in the present work. 
However, it is important to note that the present work may not be able to track all possible superharmonic resonances. 
The cases of a priori phase calculations for secondary superharmonic resonances are discussed in \Cref{sec:apriori_secondary}.

\FloatBarrier

\section{Variable Phase Resonance Nonlinear Modes} \label{sec:tracking_method}

To more generally track superharmonic resonances, a new method termed variable phase resonance nonlinear modes (VPRNM) is proposed in this section.
From \Cref{sec:apriori_phase}, it is clear that fixing a constant phase criteria is impossible for the hysteretic nonlinear forces. Furthermore, for MDOF systems, the fundamental motion may not be at zero phase during the superharmonic resonance. Therefore, a more general approach is needed for tracking the superharmonic resonances than a constant phase criteria as is provided in \Cref{sec:apriori_phase} and the previous approach of PRNM \cite{volvertPhaseResonanceNonlinear2021, volvertResonantPhaseLags2022, abeloosControlbasedMethodsIdentification2022a}. 
The proposed method of PRNM in \cite{abeloosControlbasedMethodsIdentification2022a} allows for MDOF systems (see \Cref{sec:previous_prnm}); however, the present work differs in that it calculates a variable phase for the superharmonic resonance dependent only on the current response state from HBM. This allows for tracking superharmonic resonances for nonlinearities that produce harmonics at amplitude dependent phases such as the hysteretic nonlinearities presented in \Cref{fig:fbroad_phase}.

First a vector is defined to capture the phase and magnitude of the broadband excitation of the higher harmonic $ n $ as
\begin{equation} \label{eq:vecFbroad}
	\boldsymbol{F_{broad,n}} = 
	\begin{bmatrix}
		F_{nc,broad} \\ F_{ns,broad}
	\end{bmatrix}
\end{equation}
This force is a function of the displacements of all of the harmonics less than $ n $ (i.e., 0 and 1 through $ n-1 $). This calculation can be easily completed using the AFT algorithm with a simple update to the input arguments to eliminate the $ n $th and higher harmonic components. VPRNM treats $ \boldsymbol{F_{broad,n}} $ as an external forcing to the superharmonic, and thus maintains phase resonance of the response with respect to this forcing. This is similar to how a primary resonance can be tracked with phase resonance with respect to an external force \cite{peetersDynamicTestingNonlinear2011}.
It is emphasized that this phase criteria assumes that the $F_{kq,sup,n}$ forces defined in \Cref{sec:decomposeFnl} remain small and that the higher harmonics do not influence the response of the lower harmonics. If $F_{nq,sup,n}$ is not small then the excitation of the higher harmonic may not be well approximated by $\boldsymbol{F_{broad,n}}$. Similarly, if $F_{kq,sup,n}$ for $k < n$ is not small, then the phase of the lower harmonics may be influenced, distorting the calculated $\boldsymbol{F_{broad,n}}$.

As presented in \Cref{sec:hbm}, HBM includes unknowns for harmonic coefficients of each degree of freedom at a fixed frequency and forcing level. For VPRNM, the method seeks to track the superharmonic resonance and therefore the response at only a single frequency at a given force level. Therefore, frequency becomes an additional unknown, and an additional equation is required. VPRNM proposes to use the orthogonality of the response of the superharmonic $ n $ that is being tracked with the force vector $ \boldsymbol{F_{broad,n}} $ as
\begin{equation} \label{eq:sprnm_constraint}
	\dfrac{\boldsymbol{F_{broad,n}}^T }{||\boldsymbol{F_{broad,n}}||_2}
	\begin{bmatrix}
		X_{nc} \\ X_{ns}
	\end{bmatrix}
	=
	0.
\end{equation}
For the SDOF systems considered here, orthogonality in the complex plane as expressed by the inner product in \eqref{eq:sprnm_constraint} is equivalent to enforcing a $ \pi/2 $ phase between the forcing of the $ n $th harmonic ($ \boldsymbol{F_{broad,n}} $) and the response of the $ n $th harmonic.

Implementing VPRNM as additional constraint for an existing HBM routine only requires the addition of \eqref{eq:sprnm_constraint} and the use of a second AFT evaluation to determine $ \boldsymbol{F_{broad,n}} $. The gradients required for many numerical solvers\footnote{For instance, Newton-Raphson requires a linearization of the system around the current state to iteratively update the calculation of the roots of the system of nonlinear equations.} are straightforward to derive using the gradients from AFT that are generally required with HBM. Furthermore, the present approach is in a form that could be generalized to MDOF systems by replacing the scalars $ F_{nc,broad}, F_{ns,broad}, X_{nc} $, and $  X_{ns} $ with vectors in \eqref{eq:vecFbroad} and \eqref{eq:sprnm_constraint}. However, the accuracy of such a generalization is unclear and is left to future work. To improve conditioning, the numerical implementation divides \eqref{eq:sprnm_constraint} by the magnitude of $ \boldsymbol{F_{broad,n}} $ (i.e., the 2-norm $||\boldsymbol{F_{broad,n}}||_2$). The present work does not consider any external excitation of higher harmonics. It is hypothesized that if the $n$th harmonic is externally excited, then one could modify \eqref{eq:sprnm_constraint} to utilize the sum of $\boldsymbol{F_{broad,n}}$ and the external excitation of the $n$th harmonic.

Using the constraint of \eqref{eq:sprnm_constraint}, VPRNM now has $ N(2H+1) + 1 $ unknowns for a system with $ N $ degrees of freedom or $ 2H + 2 $ unknowns for the present work with a SDOF system. These unknowns correspond to the harmonic coefficients of the zeroth and first $ H $ harmonics and the frequency. HBM provides $ N(2H + 1) $ equations and \eqref{eq:sprnm_constraint} is the additional required equation.
Continuation is then used to calculate responses over a range of force levels.
Since the HBM equations are augmented with an additional equation without modification, standard approaches such as those discussed in \cite{krackHarmonicBalanceNonlinear2019} could be applied to analyze the stability of the HBM solutions along the VPRNM backbone.

An alternative tracking approach could be to calculate the point where the first derivative of the response amplitude (or the amplitude of a specific harmonic) is equal to zero. However, implementation of such an algorithm would prove challenging since the residual function now requires the first derivative of the response including the first derivative of the nonlinear forces. Therefore, the gradient based solvers used for the present work would require the evaluation of the second derivative of the nonlinear forces with respect to the displacements. 
Recent work has demonstrated such an approach for smooth nonlinearities but showed significantly increases in computation times compared to phase resonance based approaches \cite{razeTrackingAmplitudeExtrema2024}.
Furthermore, the piecewise linear representations of frictional nonlinearities considered here do not lend themselves to calculating useful second derivatives via an AFT procedure.\footnote{Since the nonlinearity is treated as piecewise linear, the second derivative is zero or undefined at all time points in AFT.} 
On the other hand, the proposed method can easily be added to existing harmonic balance codes since the additional equation only requires an AFT evaluation that is already required for HBM and an inner product that is trivial to implement.

The present implementation of VPRNM differs from the previously proposed PRNM \cite{volvertPhaseResonanceNonlinear2021} in that the phase criteria is formulated to be amplitude dependent and allows VPRNM to be applied to hysteretic nonlinearities. 
Additionally, the phase criteria is dependent on the HBM solution for the phase of the first harmonic response rather than the phase of the underlying linear system at the frequency of interest as is used in \cite{abeloosControlbasedMethodsIdentification2022a}. 
An additional difference exists in how the exact equations are formulated. The present work considers external forcing at a fixed phase. On the other hand, PRNM uses forcing based on velocity feedback and applies a constraint on the phase of the superharmonic of interest to  make the solution unique \cite{volvertPhaseResonanceNonlinear2021}. For the present work with SDOF systems, this difference is only in the formulation of the equations, but does not cause differences in the results. 
Rather, the results of VPRNM differ from PRNM because the amplitude dependent phase criteria is used. 
The present implementation of the phase criteria is used instead of velocity feedback because it is more straightforward to implement for the variable phase criteria.
Lastly, since \eqref{eq:sprnm_constraint} relies on knowledge of the exact nonlinearity the current form of VPRNM cannot be applied to experiments.

\section{Results} \label{sec:results}

\newcommand{\basecaption}[5]{{#1}:1 superharmonic resonances for {#2} using harmonics 0 and 1 through {#3}. Left: FRCs (top: amplitude and bottom: phase of {#4} harmonic). Right: evolution of response with varying force level. Axes are shared between the left/right and top/bottom plots. The shaded region on the right is the envelope of plots on left; dots on right indicate force levels used on left FRC plots. Blue is the truth solution from HBM at discrete force levels; orange is the present approximation found by continuation. The dashed gray lines on the bottom left indicate the phase of the $ \boldsymbol{F_{broad,{#1}}} $ excitation of the {#4} harmonic (see \eqref{eq:fbroad} and \eqref{eq:phibroad}). {#5}}

\newcommand{\forcearrow}[0]{The arrows indicate the approximate direction of evolution of the superharmonic resonance for increasing external force.}

\newcommand{\PropResp}[0]{The response is proportional to force amplitude, so all curves on the left directly overlay. Normalization also results in single points for the orange and blue responses on the left.}

This section presents FRCs for the eight different nonlinear forces to highlight superharmonic resonances. These plots highlight the utility of the proposed VPRNM tracking method for superharmonics presented in \Cref{sec:tracking_method} and how the calculations in \Cref{sec:apriori_phase} can be applied. Additionally, the limitations of VPRNM are discussed. 
Results are divided by the nonlinear force types with the stiffening nonlinearities in \Cref{sec:results_harden}, the conservative softening nonlinearities in \Cref{sec:results_conserve_soft}, the even nonlinearity in \Cref{sec:results_uni}, and the damping and hysteretic nonlinearities in \Cref{sec:results_hyst}.
For each type of nonlinearity, sections are further divided between primary superharmonic resonances (e.g., 3:1 for most nonlinearities or 2:1 for the unilateral spring) and secondary superharmonic resonances corresponding to higher harmonics.
Finally, \Cref{sec:timing} discusses the relative errors for different cases and the computation time for the proposed method.
Throughout, VPRNM results are compared to an amplitude resonance of the superharmonic and differences are defined as errors since it is expected that the amplitude resonance would be most important for design. Alternatively, differences in responses could be interpreted as differences in definitions rather than error.
Plots in this section focus on the superharmonic resonances rather than the full FRCs. The context of the full FRCs are provided in \Cref{sec:system} and \Cref{sec:appendix_frcs}.

Simulations with stiffening Duffing, quintic stiffness, and unilateral spring nonlinearities are all conducted with the zeroth and first 12 harmonics since these superharmonic resonances showed notable changes in behavior when increasing the highest considered harmonic from 3 to 12. For all other nonlinearities, only the minimum number of harmonics necessary for the superharmonic resonance are included (e.g., 0th and first $ n $ harmonics for an $ n$:1 superharmonic resonance) since the behavior does not notably change with the inclusion of additional harmonics. Including the minimal number of harmonics allows for more clear isolation of the individual superharmonic resonance being considered for the clarity of the plots.

\subsection{Conservative Stiffening Nonlinearities} \label{sec:results_harden}

\subsubsection{Primary Superharmonics} 

\FloatBarrier

\begin{figure}[!h]
	\centering
	\includegraphics[width=0.7\linewidth]{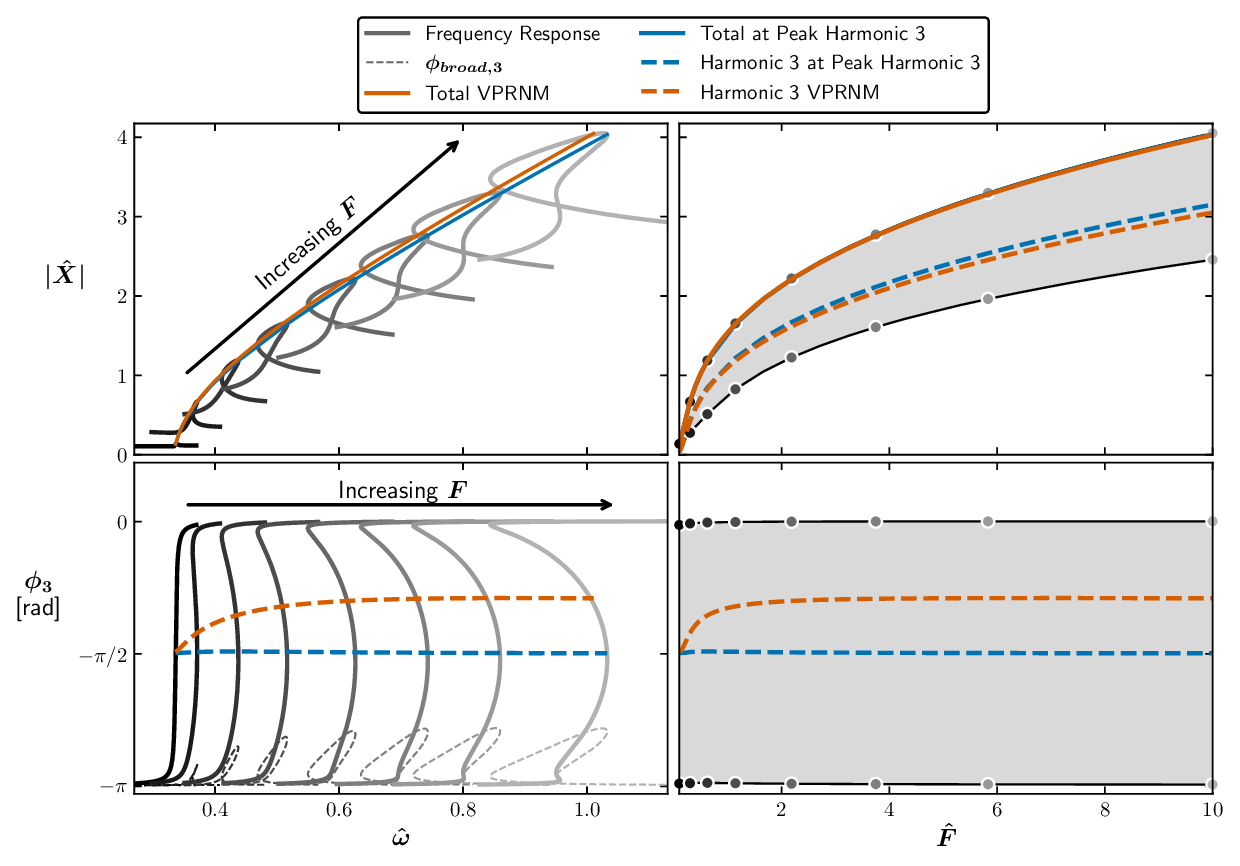}
	\caption{ \basecaption{3}{stiffening Duffing nonlinearity}{12}{third}{\forcearrow} } 
	\label{fig:duffing_FRC}
\end{figure}

The response of the stiffening Duffing oscillator near a 3:1 superharmonic resonance is shown in \Cref{fig:duffing_FRC} (see \Cref{fig:example_FRC} for wider context). On the top left of \Cref{fig:duffing_FRC}, the FRCs for various forcing amplitudes are shown in gray and exhibit loops near the superharmonic resonance.  The orange line in the top left of \Cref{fig:duffing_FRC} shows that the proposed VPRNM method of tracking the superharmonic resonance captures points near the peak amplitude (blue) of the 3:1 superharmonic resonance. 
The blue lines correspond to the details of the response at the frequency with the maximum contribution from the third harmonic. These quantities are calculated by post-processing full FRCs at discrete force levels and require running HBM for multiple force levels over a range of frequencies. 
Alternatively, orange lines represent the responses tracked via VPRNM, which does continuation with respect to force level and only calculates one response for each force level. Therefore, the VPRNM curves require significantly less computation than the curves for the peak third harmonic (see \Cref{sec:timing}).

The top right of \Cref{fig:duffing_FRC} shows the bounds of the FRCs (shaded region), magnitude of the third harmonic (dashed lines), and total response magnitude for VPRNM and the frequency with maximum third harmonic contribution (solid colored lines) as a function of external force level.
Here, the height of the shaded region corresponds to the additional amplitude caused by the superharmonic resonance.
On the top right, VPRNM directly overlays the total amplitude at the peak of the third harmonic, and consequently the solid blue line is not visible. 
Since the total amplitude is captured, the slight errors in frequency (see the top left) and magnitude of the third harmonic (dashed lines on top right) are not a significant limitation of VPRNM.

The bottom of \Cref{fig:duffing_FRC} illustrates how the phase of the third harmonic evolves for the FRCs shown on the top of the figure (with respect to frequency and external force level on left and right respectively). For each FRC, the phase of the third harmonic rises from $ -\pi $ to $ 0 $ as the 3:1 superharmonic resonance is crossed. 
The dashed orange lines on the bottom of \Cref{fig:duffing_FRC} show that the phase of the third harmonic calculated from VPRNM increases from near $ -\pi/2 $ at low amplitudes and then appears to saturate at higher phase near $ -\pi/3 $. This differs from the the phase of the third harmonic at the peak amplitude of the third harmonic (calculated from full FRCs), which remains at phase of $ -\pi/2 $.

The shift in the phase criteria for VPRNM is caused by the loops in the excitation phase of the third harmonic ($ \phi_{broad,3} $ from \eqref{eq:phibroad}), which is plotted as dashed gray lines on the bottom left of \Cref{fig:duffing_FRC}.
These loops correspond to a shift in the phase of the fundamental harmonic caused by the presence of the third harmonic resonance due to $ F_{1c,sup,3} $ and $ F_{1s,sup,3} $ in \eqref{eq:force_decomp_hbm}.
Here, the shift in the phase criteria for VPRNM explains the slight errors in capturing the amplitude of the superharmonic resonance with VPRNM.
The forces $ F_{1c,sup,3} $ and $ F_{1s,sup,3} $ correspond to how the assumptions related to superposition in VPRNM are violated (see \Cref{sec:decomposeFnl}). 
The results in the present section empirically illustrate when sufficient accuracy can be achieved while neglecting the superposition effects.\footnote{Superposition effects are only neglected in the additional constraint equation for VPRNM, which is used to set the frequency at a given force level. The full response at the selected frequency is calculated with HBM without assuming any form of superposition.} 
It is noted that the phase resonance criteria in \cite{abeloosControlbasedMethodsIdentification2022a} is partially based on the phase of the first harmonic calculated without any influence of the higher harmonics and that calculation is not influenced by these loops in the first harmonic phase (see \Cref{sec:previous_prnm}).

\FloatBarrier

\begin{figure}[!h]
	\centering
	\includegraphics[width=0.7\linewidth]{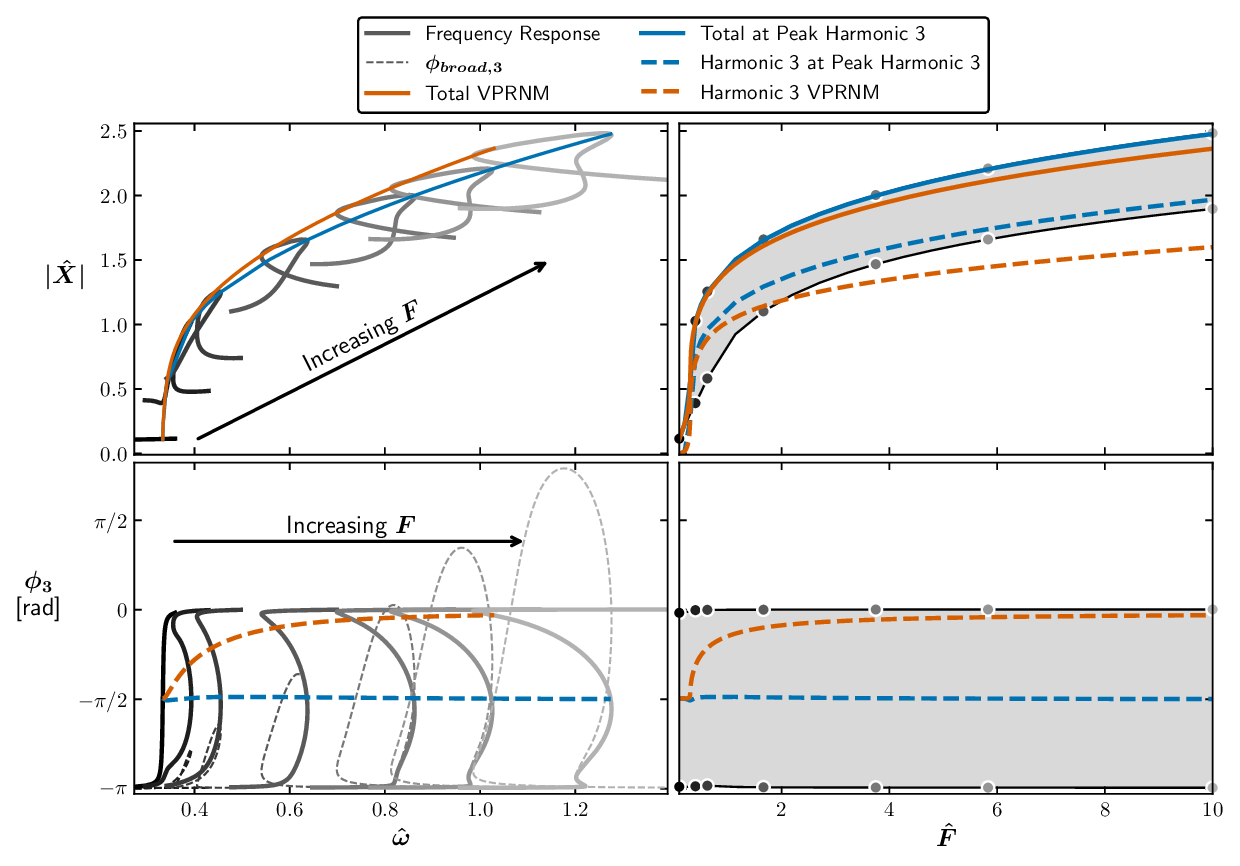}
	\caption{
		\basecaption{3}{quintic stiffness nonlinearity}{12}{third}{\forcearrow}
	} 
	\label{fig:quintic_FRC}
\end{figure}

\Cref{fig:quintic_FRC} shows similar behavior for the 3:1 superharmonic resonance of the quintic stiffness as is observed for the stiffening Duffing oscillator.
Specifically, $ \phi_{broad,3} $ develops loops in bottom left of \Cref{fig:quintic_FRC} resulting in VPRNM missing the peak amplitudes on the top of \Cref{fig:quintic_FRC}.\footnote{Some slight faceting is shown on the top right of \Cref{fig:quintic_FRC} for the bounds of the shaded gray region, which are based on FRCs at discrete force values. Instead, VPRNM does continuation with respect to the force level and more easily produces a smooth curve with lower computational cost.}
However, VPRNM still provides a useful understanding of the superharmonic resonance compared to the alternative of neglecting the superharmonic resonance as is often done with nonlinear modal methods \cite{schwarzValidationTurbineBlade2020}.
It is possible that the extent of variations in $ \phi_{broad,3} $ could be analyzed based on the nonlinear forces along the VPRNM curve to assess the errors in the calculation. Such calculations are outside the scope of the present work. 
Although one could adopt a constant phase criteria of $ -\pi/2 $ similar to the analysis of PRNM, such an approach cannot be applied generally to all of the nonlinearities considered here.

\FloatBarrier

\subsubsection{Secondary Superharmonics}

\FloatBarrier

\begin{figure}[!h]
	\centering
	\includegraphics[width=0.7\linewidth]{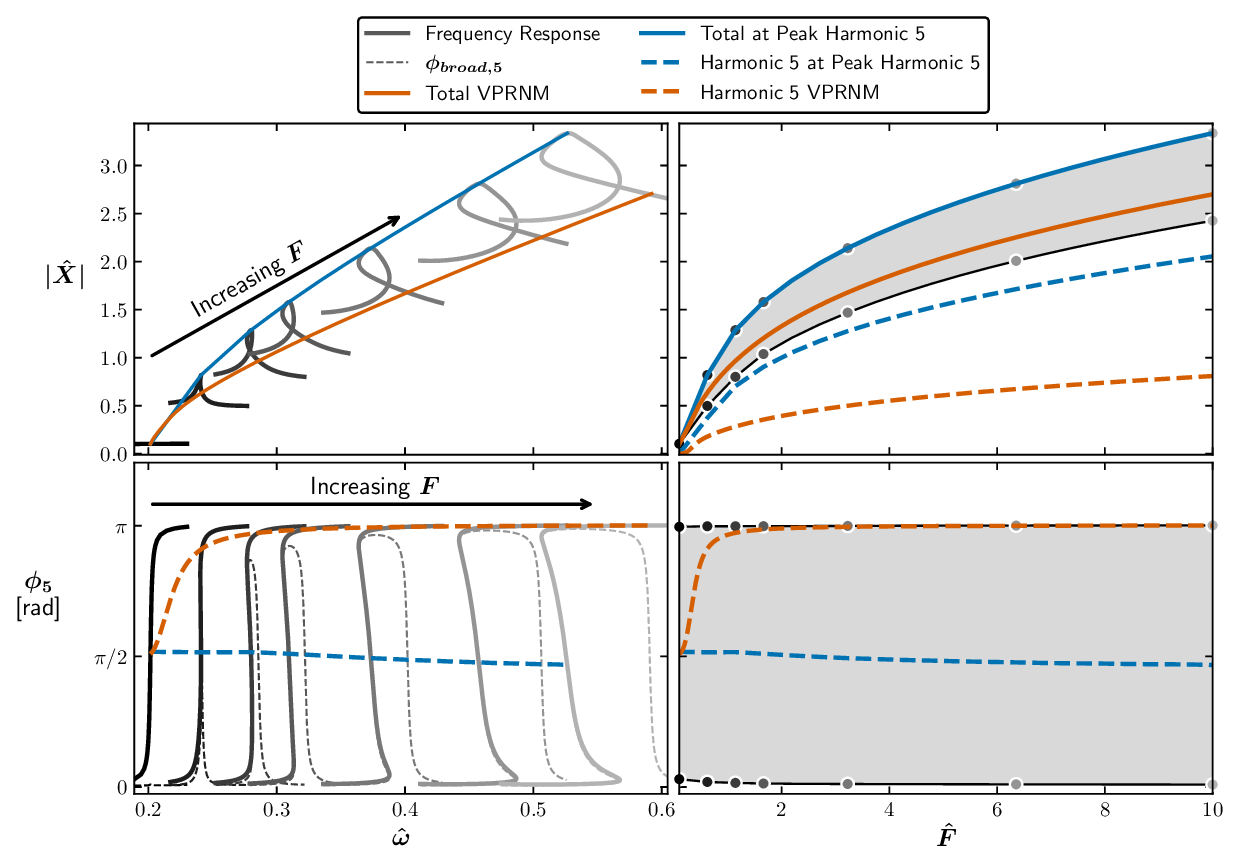}
	\caption{
		\basecaption{5}{stiffening Duffing nonlinearity}{12}{fifth}{\forcearrow}
	}
	\label{fig:duffing_5to1_FRC}
\end{figure}

For the 5:1 superharmonic resonance, the stiffening Duffing oscillator in \Cref{fig:duffing_5to1_FRC} shows more extreme errors related to similar phenomena as observed in the 3:1 case.  
In this case, the presence of the fifth harmonic results in significant phase shifts in the third harmonic that result in significant phase shifts in the predicted excitation phase of the fifth harmonic, $ \phi_{broad,5} $ from \eqref{eq:phibroad}. This phase shift of the third harmonic is attributed to the terms $F_{3q,sup,5}$ that are neglected in the VPRNM phase constraint.
Given the errors in VPRNM, the solutions could be used to initialize HBM close to the superharmonic resonance\footnote{VPRNM solutions exactly satisfy the HBM equations at a given frequency and could allow for a good HBM initialization compared to alternative methods that work best away from resonance. HBM could then be run from the VPRNM solution over a smaller range of frequencies to capture a superharmonic resonance that is not well characterized with VPRNM.} and potentially to determine which way to run continuation utilizing gradient information.

\FloatBarrier

\begin{figure}[!h]
	\centering
	\includegraphics[width=0.7\linewidth]{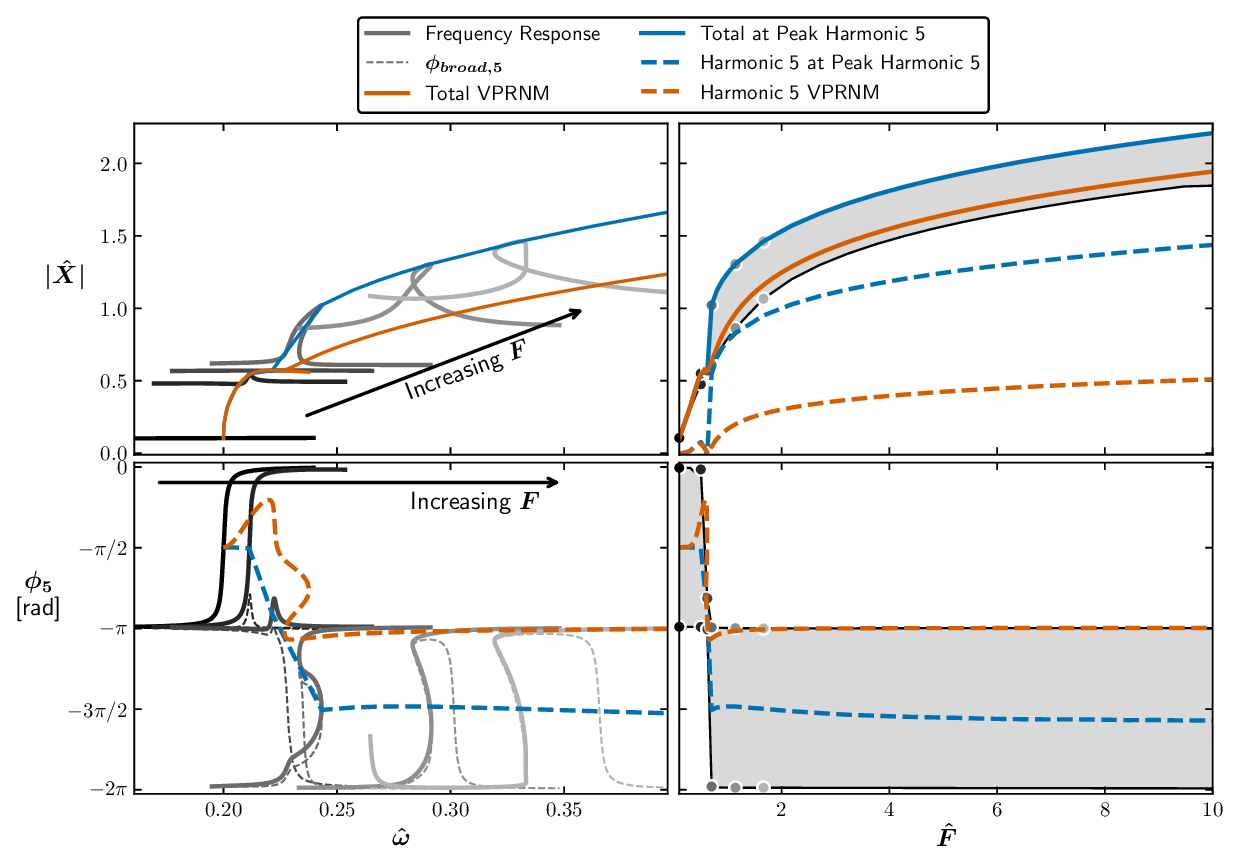}
	\caption{
		\basecaption{5}{quintic stiffness nonlinearity}{12}{fifth}{\forcearrow}
	}
	\label{fig:quintic_5to1_FRC}
\end{figure}

The 5:1 superharmonic responses for the quinitic stiffness (see \Cref{fig:quintic_5to1_FRC}) show significantly different responses than the 3:1 superharmonic case or the cases of the stiffening Duffing oscillator. 
For the quintic stiffness, nondimensional force levels of 0.10 and 0.62 do not produce a notable 5:1 superharmonic resonance (the lines with the lowest and third from lowest amplitude in the top left of \Cref{fig:quintic_5to1_FRC}). Nevertheless, a nondimensional force level of 0.4 plotted in the top left of \Cref{fig:quintic_5to1_FRC} (the line with the second lowest amplitude) does show a notable 5:1 superharmonic resonance. Then at higher force levels, behavior similar to that of the other stiffening nonlinearity cases is observed. 
This behavior at low amplitudes can be attributed to the relative magnitudes of the first and third harmonics and thus the relative importance of the terms $ -\eta X_1^5 $ and $ +20\eta X_1^4 X_3 $ from \Cref{sec:apriori_secondary} and \Cref{tab:Fbroad_second}. When $ X_3 $ grows sufficiently so that the latter term dominates, the phase of the excitation changes from $ -\pi $ to $ -2 \pi$ resulting in the shift seen on the bottom of \Cref{fig:quintic_5to1_FRC}. Therefore, the lack of a superharmonic response at a nondimensional force level of 0.62 is attributed to these terms approximately canceling out. 
Indeed, VPRNM is able to capture the forcing amplitude of this sudden shift in the phase criteria for the superharmonic resonance. On the other hand, 
constant phase criteria from perturbation analysis \cite{volvertPhaseResonanceNonlinear2021, volvertResonantPhaseLags2022, abeloosControlbasedMethodsIdentification2022a} would miss this behavior.

At higher force levels, VPRNM does a poor job of capturing the 5:1 superharmonic resonance for the quintic stiffness. As with the stiffening Duffing oscillator, this can be attributed to significant changes in the phase $ \phi_{broad,5} $ in the bottom left of \Cref{fig:quintic_5to1_FRC} near the superharmonic resonance. 
It is also noted that VPRNM performs better in the low amplitude regime where the excitation of the fifth harmonic is a direct result of the fundamental motion rather than an interaction between the first and third harmonics. Similarly, VPRNM performs poorly for the stiffening Duffing oscillator where the fifth harmonic is excited by interactions between the first and third harmonics.

\FloatBarrier

\subsection{Conservative Softening Nonlinearities} \label{sec:results_conserve_soft}

\subsubsection{Primary Superharmonics}

\FloatBarrier

\begin{figure}[!h]
	\centering
	\includegraphics[width=0.7\linewidth]{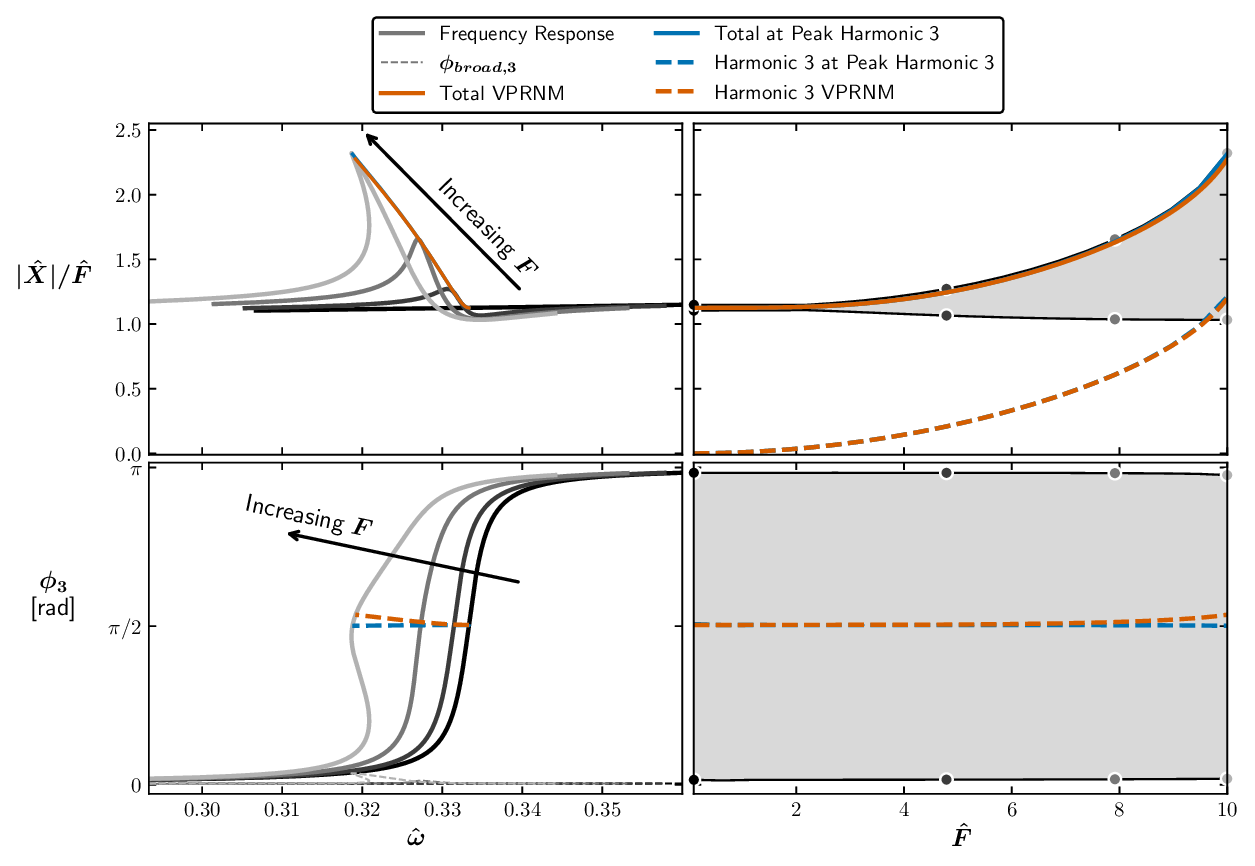}
	\caption{
		\basecaption{3}{softening Duffing nonlinearity}{3}{third}{\forcearrow}
	}
	\label{fig:softduffing_FRC}
\end{figure}

The softening Duffing nonlinearity produces 3:1 superharmonic resonances that are well captured by VPRNM in \Cref{fig:softduffing_FRC}. The higher accuracy of VPRNM here compared to the stiffening Duffing case is likely related to limiting the strength of the nonlinearity so that the total stiffness remains positive.
In addition, the phase of the third harmonic resonance remains close to $ \pi/2 $ as predicted in \Cref{sec:apriori_phase}. The change in phase criteria from $ -\pi/2 $ for the stiffening case is captured automatically with VPRNM rather than requiring new analytical analysis like previous methods \cite{volvertPhaseResonanceNonlinear2021, volvertResonantPhaseLags2022, abeloosControlbasedMethodsIdentification2022a}.

\FloatBarrier

\begin{figure}[!h]
	\centering
	\includegraphics[width=0.7\linewidth]{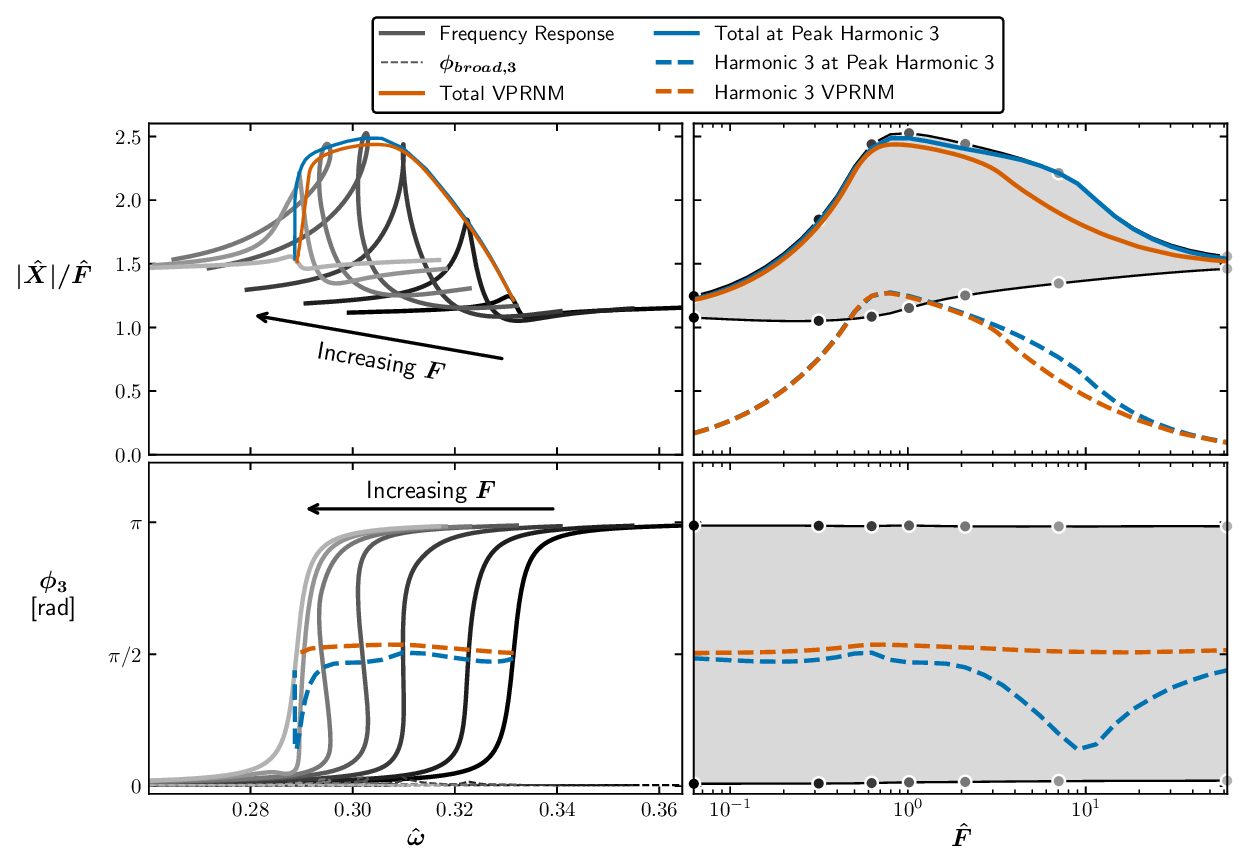}
	\caption{
		\basecaption{3}{conservative softening II nonlinearity}{3}{third}{\forcearrow}
	}
	\label{fig:consv_iwan_FRC}
\end{figure}

Due to its saturating nature,
the conservative softening II nonlinearity provides notably different characteristics (see \Cref{fig:consv_iwan_FRC}) than the previous nonlinearities. 
The top of \Cref{fig:consv_iwan_FRC} shows how the superharmonic resonance is most prominent at intermediate force amplitudes with small superharmonic responses at low and high force levels as predicted and discussed in \Cref{sec:apriori_phase}. 
Specifically, the height of the shaded region in the top right of \Cref{fig:consv_iwan_FRC} corresponds to the additional amplitude caused by the superharmonic resonance and is clearly smaller at high and low force values compared to a range of force values near the middle.
Additionally, the nondimensional linear frequency of the system using the linearized nonlinear force is 1 at low amplitudes, and the nondimensional frequency of the system without the nonlinear force is $ \sqrt{0.75} \approx 0.866 $. Here, the frequency of the 3:1 superharmonic resonances decreases from approximately one third of the linearized natural frequency at low amplitude to approximate one third of the natural frequency without the nonlinear force at high amplitudes. 
VPRNM accurately captures these behaviors because the phase $ \phi_{broad,3} $ remains near constant in the bottom left of \Cref{fig:consv_iwan_FRC}.
Analytical approaches \cite{volvertPhaseResonanceNonlinear2021, volvertResonantPhaseLags2022, abeloosControlbasedMethodsIdentification2022a} would be difficult to apply given the complexity of the nonlinear force.
Conversely, VPRNM is easily applied with the AFT procedure that is commonly used for nonlinear forces with HBM.

\FloatBarrier

\subsubsection{Secondary Superharmonics}

\FloatBarrier

\begin{figure}[!h]
	\centering
	\includegraphics[width=0.7\linewidth]{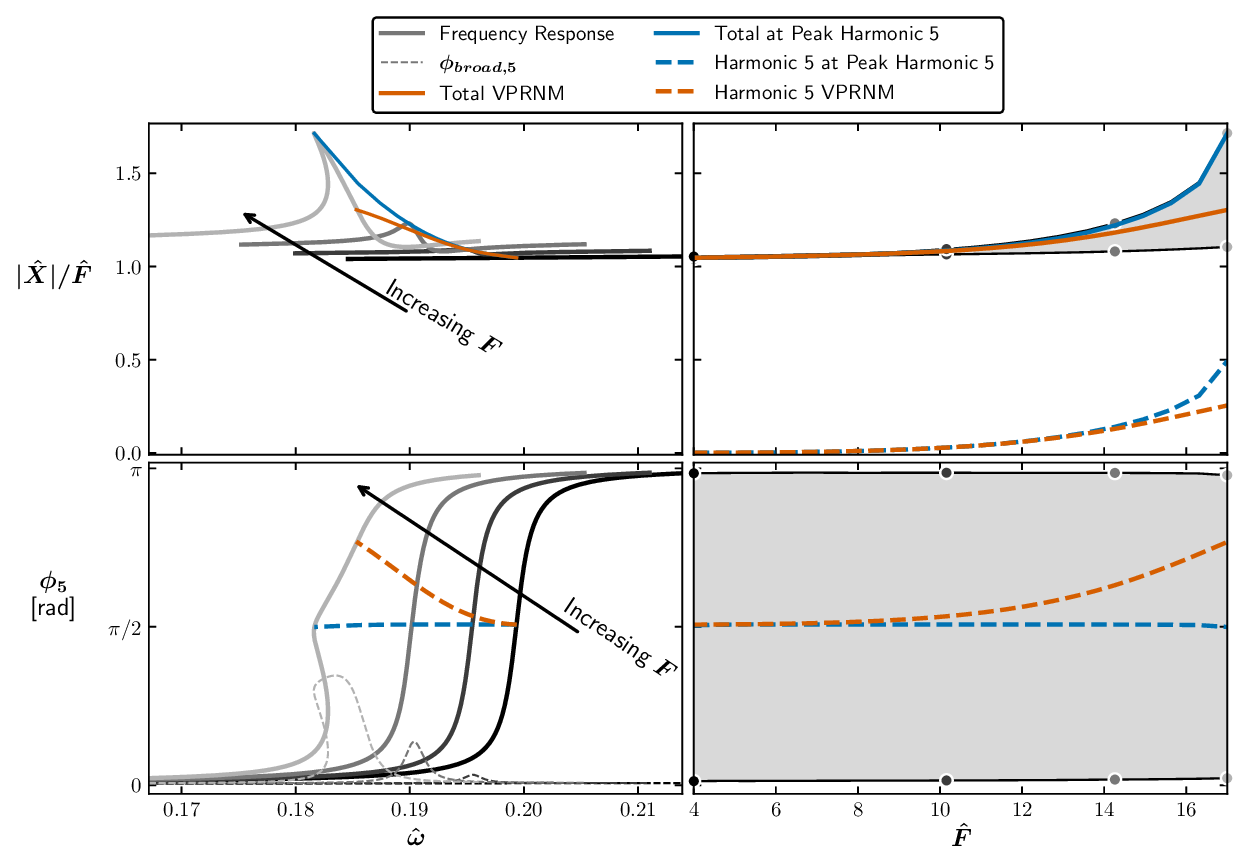}
	\caption{
		\basecaption{5}{softening Duffing nonlinearity}{5}{fifth}{\forcearrow}
	}
	\label{fig:softduffing_5to1_FRC}
\end{figure}

\begin{figure}[!h]
	\centering
	\includegraphics[width=0.7\linewidth]{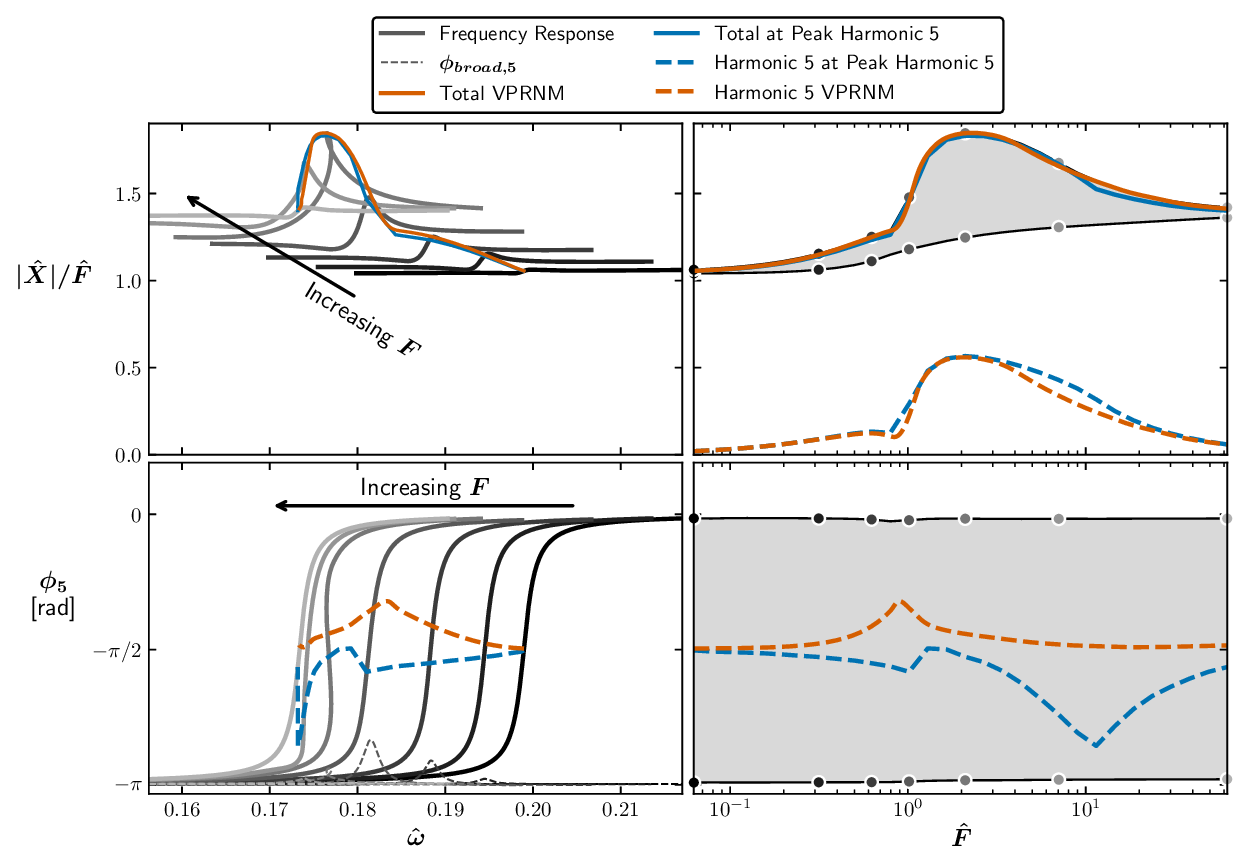}
	\caption{
		\basecaption{5}{conservative softening II nonlinearity}{5}{fifth}{\forcearrow}
	}
	\label{fig:consv_iwan_5to1_FRC}
\end{figure}

The softening Duffing nonlinearity produces 5:1 superharmonic resonances\footnote{Including more than 5 harmonics in \Cref{fig:softduffing_5to1_FRC} decreased the amplitude of the superharmonic response, but did not qualitatively alter the results. Therefore, only 5 harmonics are included to more clearly isolate the superharmonic resonance.}
as shown in \Cref{fig:softduffing_5to1_FRC}.
As with the stiffening nonlinearities, the softening Duffing nonlinearity shows some clear errors for VPRNM in predicting the peak superharmonic response as a result of shifts in the phase $ \phi_{broad,5} $ shown in the bottom left of \Cref{fig:softduffing_5to1_FRC}. 
As before, the VPRNM results could be used to initialize HBM solutions instead of simply using the results from VPRNM.
Contrarily, the 5:1 superharmonic resonance for the conservative softening II nonlinearity shows similar behavior to the 3:1 case and good accuracy for VPRNM (see \Cref{fig:consv_iwan_5to1_FRC}).
One possible reason for the good performance of VPRNM is that the fundamental harmonic motion excites the fifth harmonic directly for the conservative softening II model. This differs from the two Duffing nonlinearities where both the first and third harmonics are required to excite the fifth harmonic. Furthermore, VPRNM only breaks down for the quintic stiffness at higher amplitudes where the presence of the third harmonic is important for the excitation of the fifth superharmonic. 
These results suggest that VPRNM may be more effective for higher superharmonic resonances if the nonlinear force under fundamental motion directly excites the higher harmonic.

\FloatBarrier

\FloatBarrier

\subsection{Even Nonlinearity} \label{sec:results_uni}

\subsubsection{Primary Superharmonic}

\FloatBarrier

\begin{figure}[!h]
	\centering
	\includegraphics[width=0.7\linewidth]{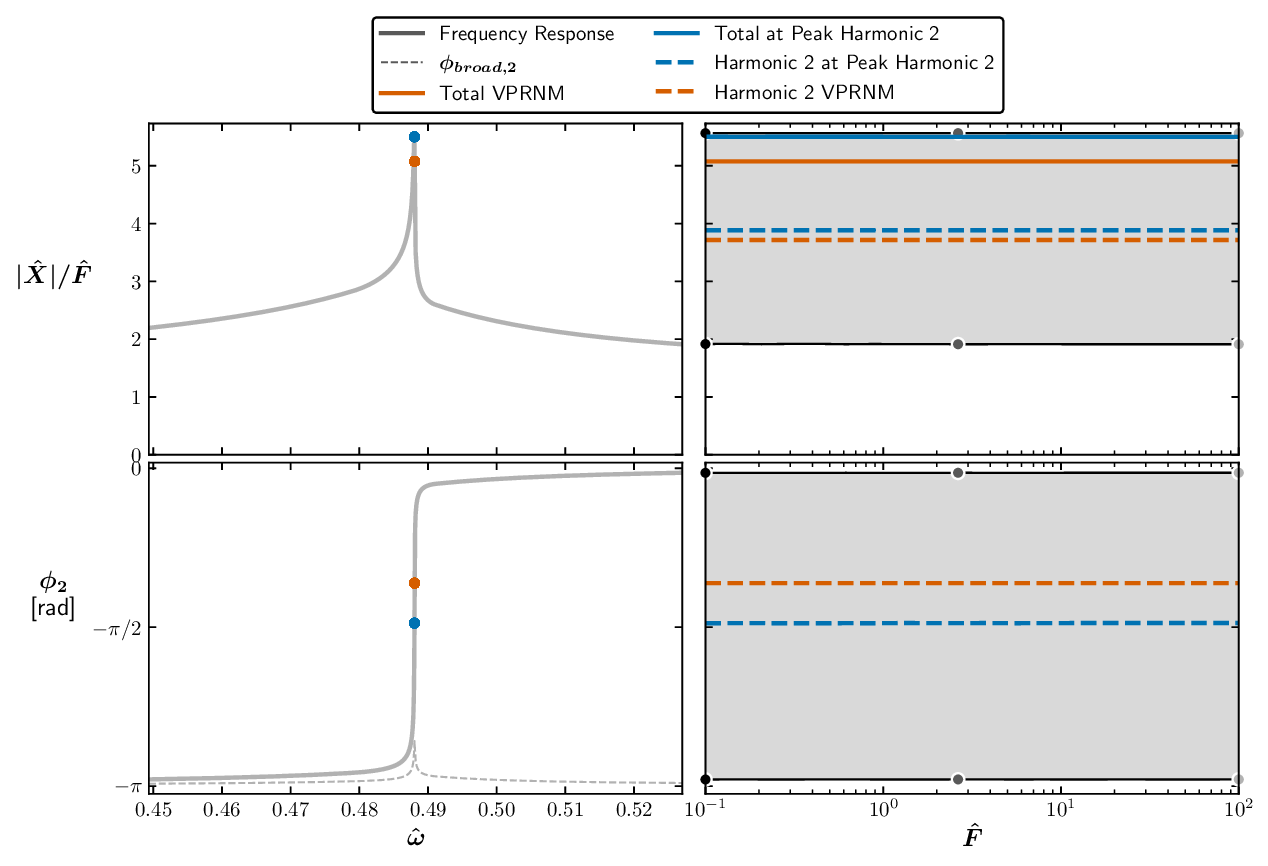}
	\caption{
		\basecaption{2}{unilateral spring nonlinearity}{12}{second}{\PropResp}
	}
	\label{fig:unispring_FRC}
\end{figure}

In addition to the conservative and odd functions of displacement previously discussed, even nonlinearities, such as unilateral contact, occur in real structures (see \Cref{sec:unispring_force} for formulated a unilateral spring as an even nonlinearity).
The case of an SDOF system with unilateral contact is shown in \Cref{fig:unispring_FRC}.
Because of the nature of the unilateral spring for this case, the response is proportional to the external force level as shown on the right of \Cref{fig:unispring_FRC} (and as was predicted in \Cref{sec:apriori_phase}). On the left of \Cref{fig:unispring_FRC}, the FRCs all directly overlay when normalized so only one is visible.
Here, VPRNM has notable errors in the prediction of the peak response amplitude, which are attributed to shifts in the fundamental harmonic and $ \phi_{broad,2} $ as was the case with the stiffening nonlinearities.

Previous efforts at tracking superharmonic resonances have not considered nonlinear forces of the form of a unilateral spring, so cannot be directly compared.
Nevertheless, as noted in \Cref{sec:apriori_phase}, these results are inconsistent with the conjecture from \cite{abeloosControlbasedMethodsIdentification2022a} that even nonlinearities should result in phase lags of 0.\footnote{Also see discussion on the conjecture from \cite{abeloosControlbasedMethodsIdentification2022a} in \Cref{sec:previous_prnm}.}
An alternative approach to tracking the superharmonic resonance here would be to fix the phase at a constant value of $ -\pi/2 $ as is derived in \Cref{sec:apriori_phase}. Such an approach would match the peak response of the system given that the blue line on the right of \Cref{fig:unispring_FRC} indicates that the peak response of the second harmonic occurs at a phase of $ -\pi/2 $.
Notwithstanding, a constant phase criteria is not selected in this work since it cannot be applied for the hysteretic forces.

\FloatBarrier

\subsubsection{Secondary Superharmonics}

\FloatBarrier

\begin{figure}[!h]
	\centering
	\includegraphics[width=0.7\linewidth]{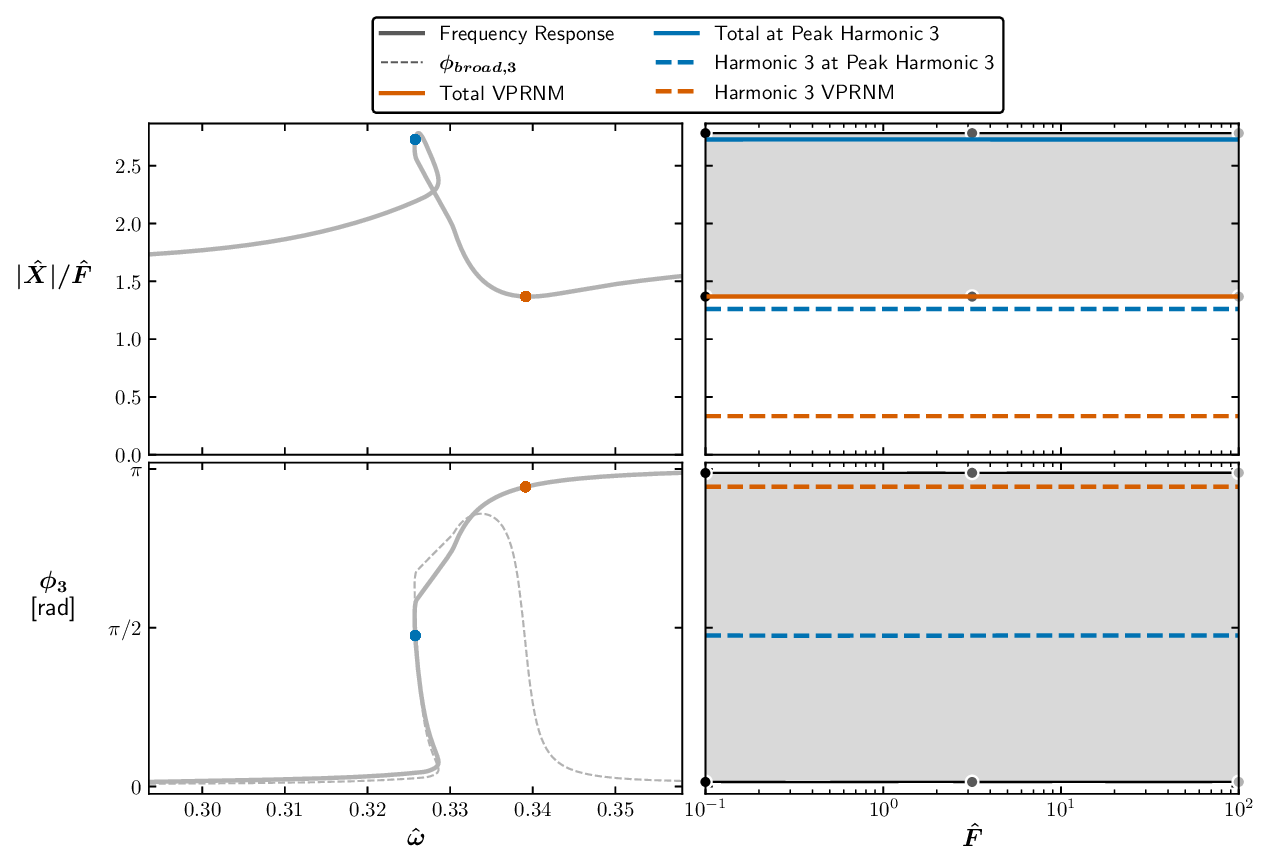}
	\caption{
		\basecaption{3}{unilateral spring nonlinearity}{12}{third}{\PropResp}
	}
	\label{fig:unispring_3to1_FRC}
\end{figure}

\begin{figure}[!h]
	\centering
	\includegraphics[width=0.7\linewidth]{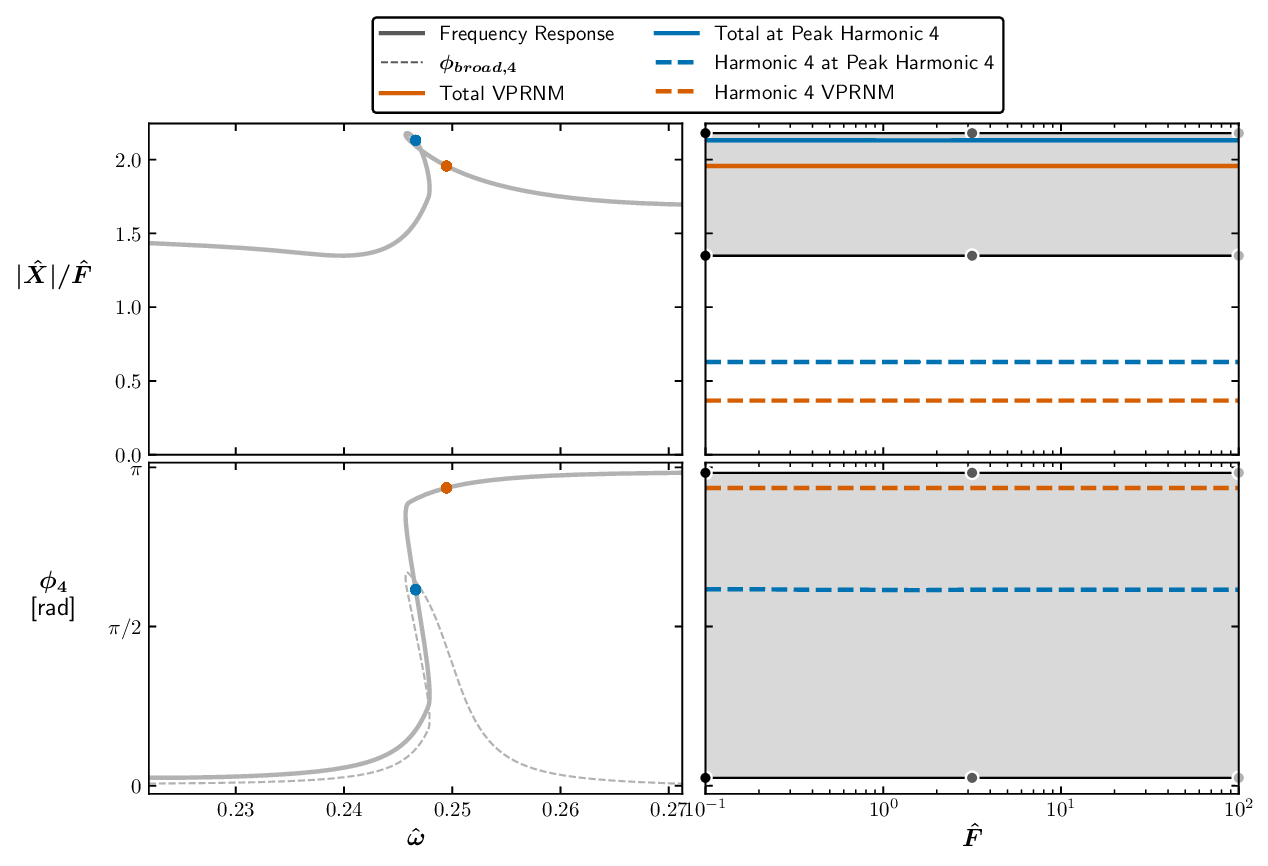}
	\caption{
		\basecaption{4}{unilateral spring nonlinearity}{12}{fourth}{\PropResp}
	}
	\label{fig:unispring_4to1_FRC}
\end{figure}

The even nonlinearity of the unilateral spring gives rise to 3:1 and 4:1 superhamonic resonances (see Figures \ref{fig:unispring_3to1_FRC} and \ref{fig:unispring_4to1_FRC} respectively).
As with the 2:1 case, the 3:1 and 4:1 cases show proportional responses and clear errors for VPRNM. 
As with the stiffening nonlinearities, the secondary superharmonic resonances for the unilateral spring show higher errors than the primary superharmonic resonance due to larger shifts in the phase of excitation (here $ \phi_{broad,3} $ and $ \phi_{broad,4} $).
For these cases, VPRNM may be used to initialize HBM close to the superharmonic resonance as opposed to predicting responses.

\FloatBarrier

\subsection{Damping and Hysteretic Nonlinearities} \label{sec:results_hyst}

\subsubsection{Primary Superharmonics}

\FloatBarrier

\begin{figure}[!h]
	\centering
	\includegraphics[width=0.7\linewidth]{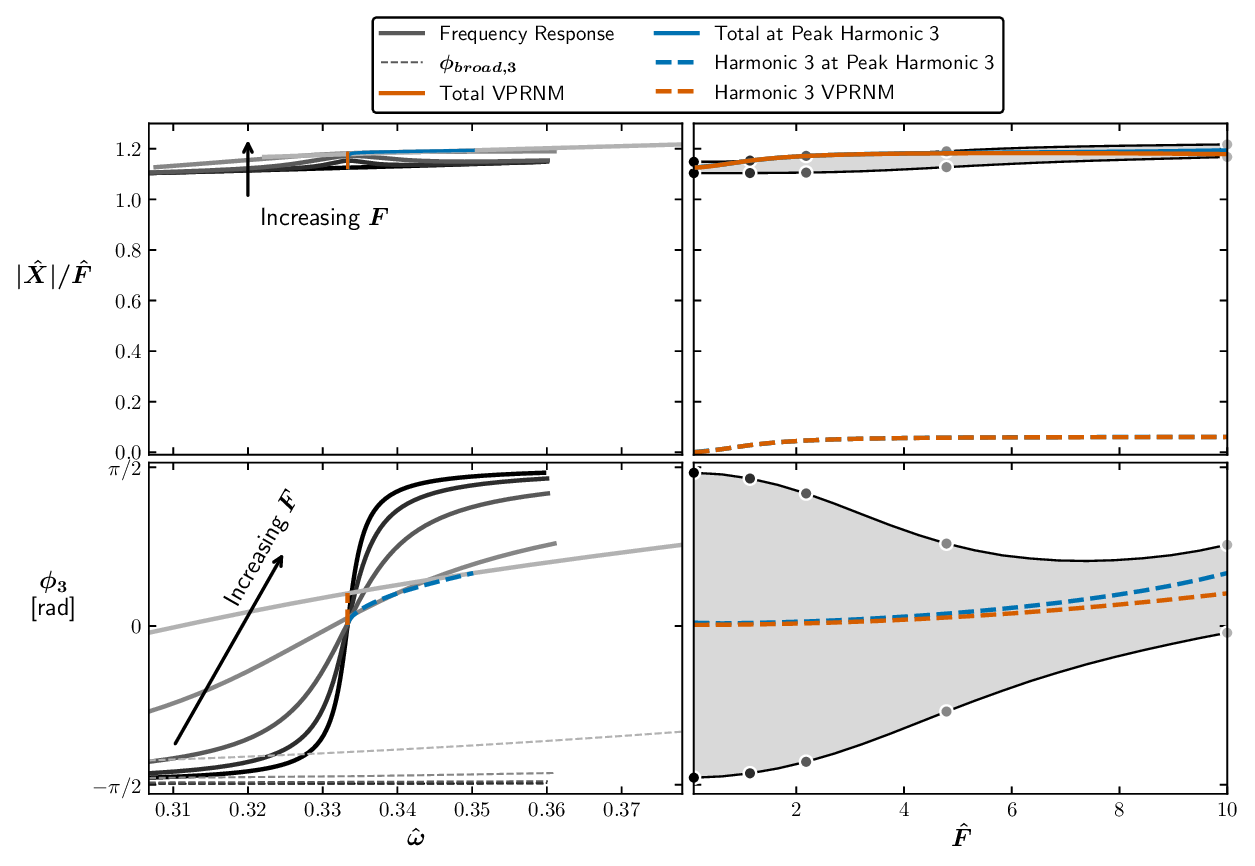}
	\caption{
		\basecaption{3}{cubic damping nonlinearity}{3}{third}{\forcearrow}
	}
	\label{fig:cubicdamp_FRC}
\end{figure}

A small 3:1 superharmonic resonance for the SDOF system with cubic damping is shown in \Cref{fig:cubicdamp_FRC}. The response is small because a strong enough nonlinear force to excite the third harmonic corresponds to a highly damped system.
VPRNM captures the shift in the phase of the fundamental harmonic related to the expanding primary resonance and adjusts the phase for the superharmonic resonance (see the bottom right of \Cref{fig:cubicdamp_FRC}).
This illustrates that the formulation of VPRNM is sufficiently general that it can be applied when the phase of the fundamental harmonic is not known a priori, and thus the method should be applicable to MDOF structures where fundamental and superharmonic resonances could interact.\footnote{Investigation of MDOF structures is left to future work.}

\begin{figure}[!h]
	\centering
	\includegraphics[width=0.7\linewidth]{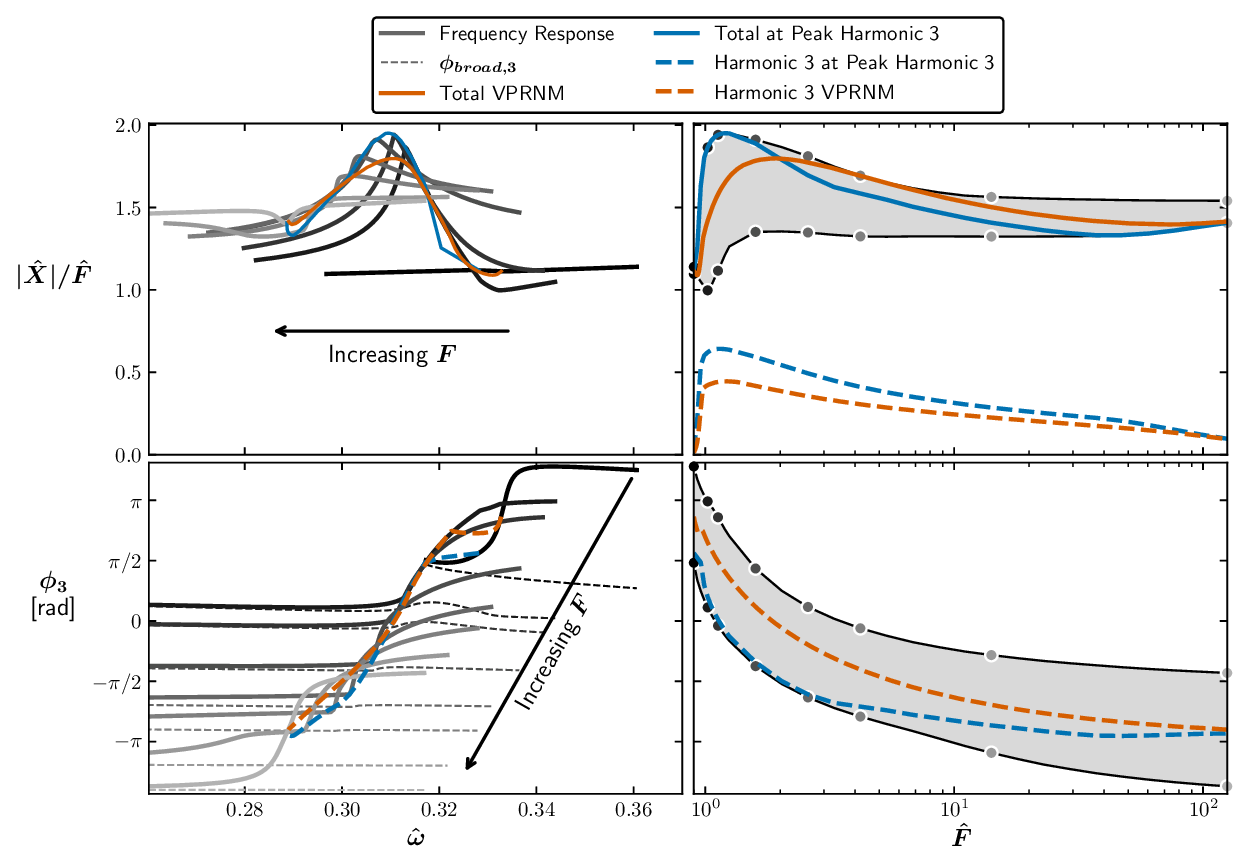}
	\caption{
		\basecaption{3}{Jenkins element}{3}{third}{\forcearrow} On the bottom left, the phase corresponding to the third harmonic for the lowest force level is only plotted for frequencies greater than approximately 0.32 rad/s since the system responds in the linear regime and does not excite the third harmonic until that frequency.  
	}
	\label{fig:jenkins_FRC}
\end{figure}

\Cref{fig:jenkins_FRC} shows the 3:1 superharmonic resonances for the hysteretic Jenkins element, which exhibits similar behavior in terms of amplitude as the conservative softening II nonlinearity due to the saturating nature.\footnote{Specifically, a large superharmonic resonance at intermediate amplitudes and small superharmonic resonances at the maximum and minimum external force values. Likewise, the frequency decreases with increasing external force levels shifting between approximately 1/3 of the low amplitude linearized frequency and 1/3 of the system frequency without the nonlinear force.}
The amplitude in the top left of \Cref{fig:jenkins_FRC} away from the superharmonic resonance generally increases with increasing force level as the primary resonance decreases in frequency, so the plotted region is closer to the primary resonance. 
VPRNM shows some clear errors in \Cref{fig:jenkins_FRC} compared to HBM, yet VPRNM requires significantly less computational time (see \Cref{sec:timing}) and alternative tracking methods are not available given the variable phase behavior for the Jenkins nonlinearity.
Additionally, VPRNM identifies that the most prominent superharmonic resonances occur near $\hat{F} = 2.5$, and this could be used to focus HBM computations near the area of interest.
Indeed, without VPRNM, one could mistakenly choose HBM forcing amplitudes that do not adequately cover the range of the most prominent superharmonic resonance (as was done in initial analysis for the present work).

At high external force amplitudes, the superharmonic resonance corresponds to a local minima (see the lightest gray FRC curve on the top left of \Cref{fig:jenkins_FRC} and the transition of the solid lines from the top to the bottom of the gray region on the top right).
The shifting phase of the superharmonic resonance is critical to understand since it indicates whether the superharmonic resonance corresponds to a local maximum or local minimum in the total response amplitude. 
VPRNM captures this qualitative transition from local maximum to local minimum of the superharmonic resonance with increasing force levels.

\FloatBarrier

\begin{figure}[!h]
	\centering
	\includegraphics[width=0.7\linewidth]{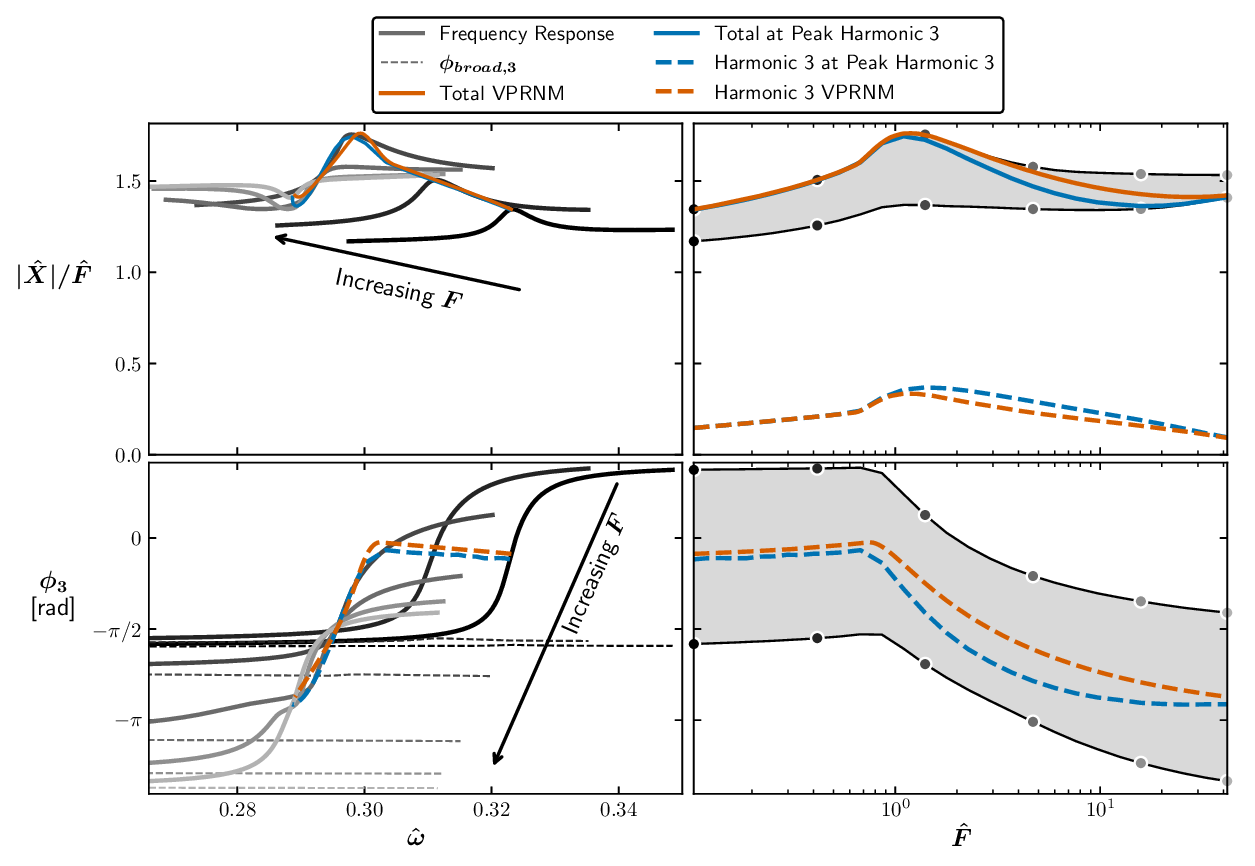}
	\caption{
		\basecaption{3}{Iwan element}{3}{third}{\forcearrow}
	}
	\label{fig:iwan_FRC}
\end{figure}

The Iwan element, which is frequently used to model bolted connections, shows similar trends to the Jenkins element (see \Cref{fig:iwan_FRC}).
The Iwan element is smoother than the Jenkins element; hence, the superharmonic resonance is visible over a wider range of force levels.
Furthermore, the superharmonic resonances are more accurately captured by VPRNM for the Iwan element (see the top right of \Cref{fig:iwan_FRC}) than for the Jenkins element. This is promising for the application of VPRNM to jointed structures, which exhibit smooth behavior more similar to that of the Iwan element than that of a single Jenkins element.

In the literature, it has been observed that hysteretic nonlinearities generally result in weak internal resonances for systems with multiple modes \cite{woiwodeEffectModalInteractions2020a}. Here, it is hypothesized that internal resonances for hysteretic systems not only require commensurate frequencies, but also require that those commensurate frequencies occur within a specific amplitude range where the excitation of the higher harmonics is the largest. The combination of these two requirements would make observations of internal resonances in MDOF systems with hysteresis more difficult than in systems with other nonlinearities (e.g., stiffening Duffing nonlinearity). VPRNM provides the opportunity for a computational efficient approach to characterize if the conditions for internal resonances are met.

\FloatBarrier

\subsubsection{Secondary Superharmonics}

The cubic damping nonlinearity did not produce any notable superharmonic resonances other than the 3:1 superharmonic resonance previously discussed, and therefore it is not included in the present section. It is hypothesized that this occurs because the cubic damping nonlinearity results in very high damping in the nonlinear regime, and consequently there is not a regime where the nonlinear forces are strong enough to excite the fifth harmonic without being strongly damped.

\FloatBarrier

\begin{figure}[!h]
	\centering
	\includegraphics[width=0.7\linewidth]{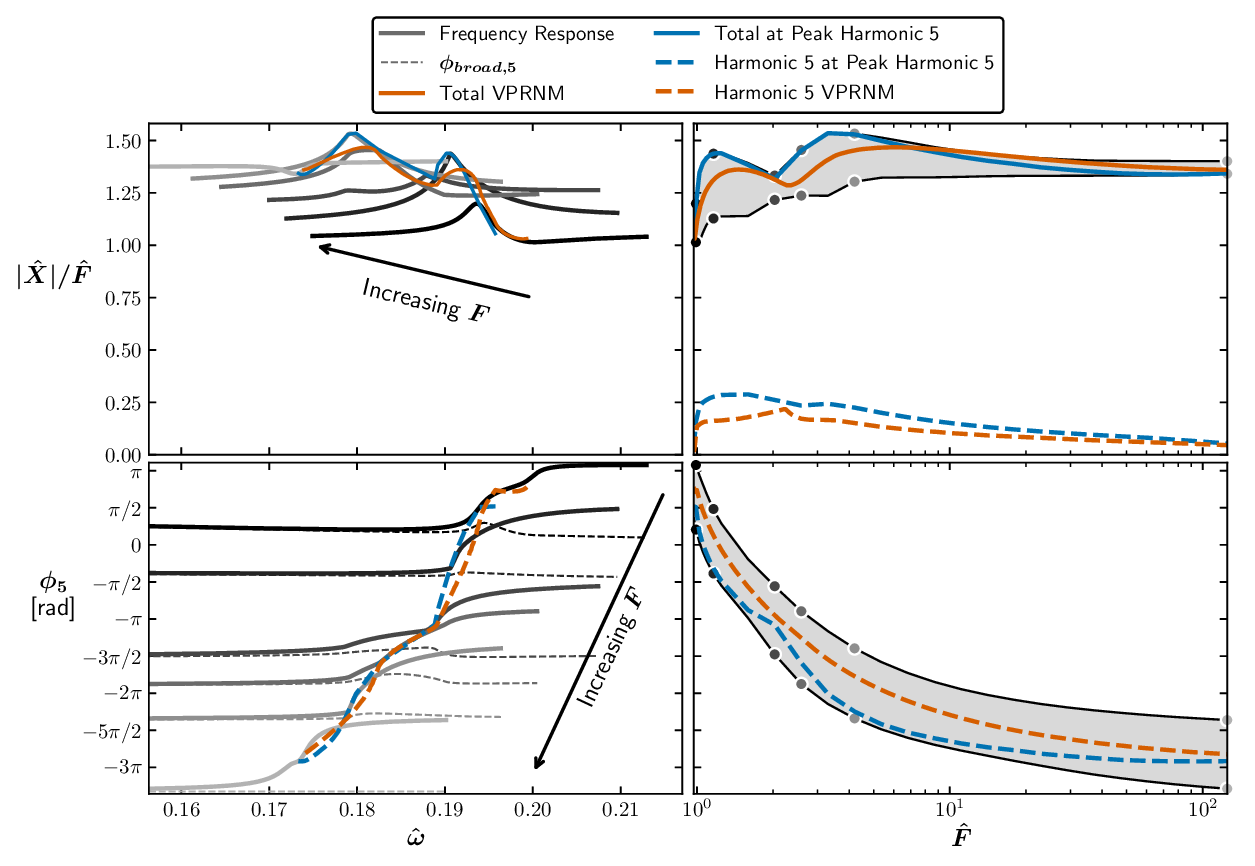}
	\caption{
		\basecaption{5}{Jenkins element}{5}{fifth}{\forcearrow}
	}
	\label{fig:jenkins_5to1_FRC}
\end{figure}

The two hysteretic nonlinearities both produce 5:1 superharmonic resonances. First, \Cref{fig:jenkins_5to1_FRC} shows the 5:1 superharmonic resonance for the Jenkins element with two external force amplitude regions where the fifth harmonic produces a clear superharmonic resonance (around nondimensional forces of 1.16 and 3.25). These regions can be clearly seen on the top right of \Cref{fig:jenkins_5to1_FRC} and can be approximately identified with VPRNM. 
Other behaviors of the 5:1 superharmonic resonance are similar to the 3:1 case including the frequency shift, a significant shift in $ \phi_{broad,5} $, and a transition to a local minimum at high amplitudes. These are qualitatively captured by VPRNM illustrating the value of the method in identifying regions of interest for superharmonic resonances.
The relative accuracy of VPRNM for the Jenkins element compared to some other 5:1 superharmonic resonance cases can be attributed to the relatively small shifts in $\phi_{broad,5}$, shown as mostly straight dashed lines in the bottom left of \Cref{fig:jenkins_5to1_FRC}.

\FloatBarrier

\begin{figure}[!h]
	\centering
	\includegraphics[width=0.7\linewidth]{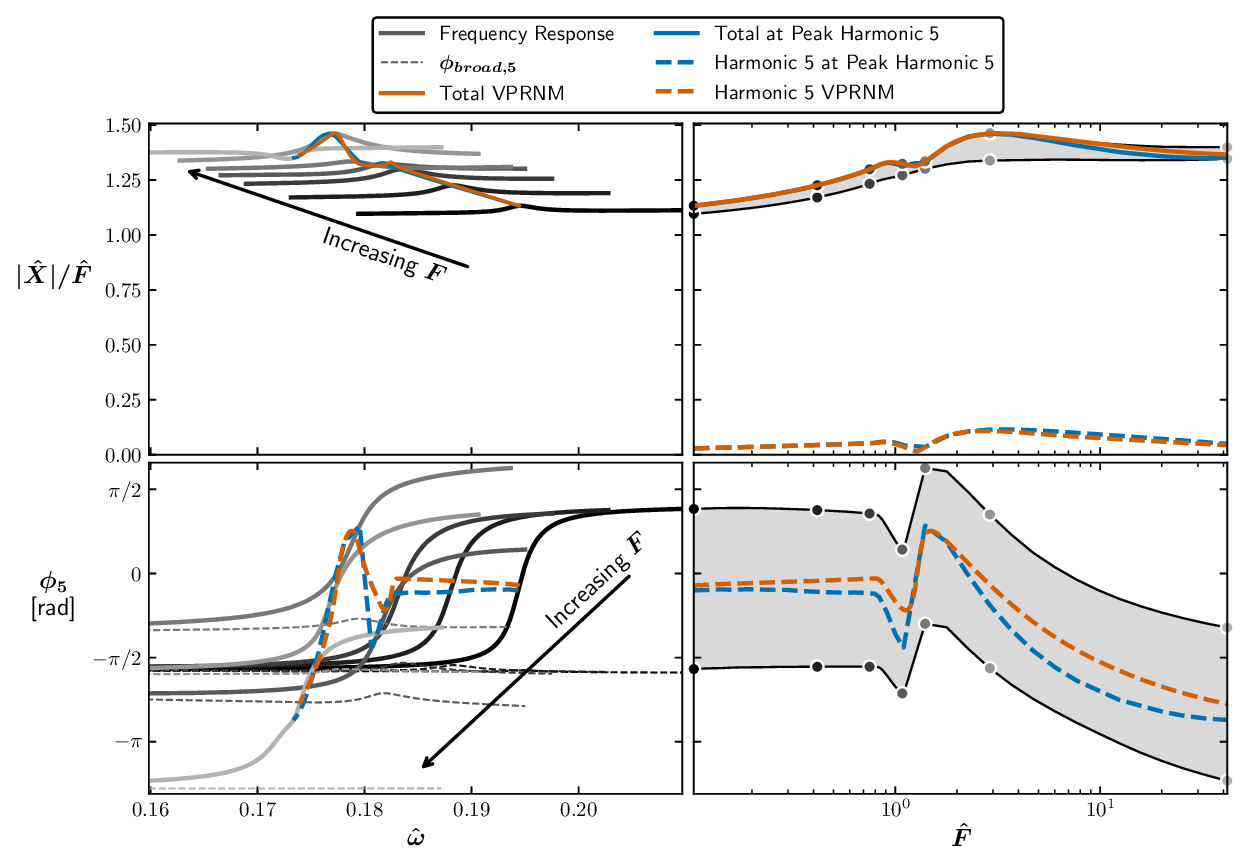}
	\caption{
		\basecaption{5}{Iwan element}{5}{fifth}{\forcearrow}
	}
	\label{fig:iwan_5to1_FRC}
\end{figure}

As a final case, the 5:1 superharmonic resonance for the Iwan element is analyzed in \Cref{fig:iwan_5to1_FRC}.
As was the case for the Jenkins element, the Iwan element shows local maxima with respect to the force at two force levels, as seen in the top right of \Cref{fig:iwan_5to1_FRC}.
VPRNM is able to accurately capture the peak amplitudes of the superharmonic resonance as shown in the top right of \Cref{fig:iwan_5to1_FRC} despite the significant shift in the phase of the fifth harmonic that would break alternative methods. 
The double peak behavior for both the Iwan and Jenkins models likely occurs as a result of the slight leveling off or local peak in the excitation magnitude of $F_{broad,5}$ in \Cref{fig:fbroad_mag} for both models.\footnote{Noting that the results in Figures \ref{fig:jenkins_5to1_FRC} and \ref{fig:iwan_5to1_FRC} are both normalized by external force level.}
The results for the hysteretic nonlinearities illustrate how the present framework of analyzing the higher harmonic components of the nonlinear forces can provide insights into the superharmonic resonance behavior of systems with different nonlinear forces.

\FloatBarrier

\subsection{Outlook and Computation Time} \label{sec:timing}

The present paper proposes VPRNM as a framework for understanding and tracking superharmonic resonances through analyzing the higher harmonic components of nonlinear forces. 
VPRNM finds a single solution at each force level and traces a one dimensional curve as opposed to HBM, which must vary frequency and external force level to produce several FRCs. This allows for a significant reduction in computation time for VPRNM compared to HBM.
Similar computational benefits are observed for other tracking methods such as the EPMC \cite{krackNonlinearModalAnalysis2015} and PRNM \cite{volvertPhaseResonanceNonlinear2021, volvertResonantPhaseLags2022, abeloosControlbasedMethodsIdentification2022a}. 
However, EPMC breaks down in the presence of modal interactions or internal resonances. 
Additionally, the present work has generalized PRNM to be applicable to a range of nonlinear forces that cannot be well characterized by a constant phase criteria (see \Cref{sec:previous_prnm} for additional discussion on the limitations of PRNM). \Cref{tab:summ_accuracy} summarizes the percent differences between VPRNM and the amplitude resonance of the superharmonic for the considered cases. Here, many cases even with large differences of up to 35\% still produce useful information about the characteristics of the superharmonic resonance with significantly reduced computational cost compared to HBM, as discussed next.

\begin{table}[h!]
	\centering
	\caption{Percent difference in the amplitude shift caused by the superharmonic resonance between VPRNM and the HBM truth solution. Using the top right subplot of each case, the error is calculated as the area between the ``Total VPRNM'' line (solid orange) and ``Total at Peak Harmonic n'' for n=2, 3, 4, or 5 line (solid blue) divided by the area between the maximum and minimum values of the plotted HBM solutions (shaded gray area). The area is calculated using linear or log scale force matching how it is plotted for each case. The area is also calculated using the same response amplitude, which is normalized by the force in some cases.} \label{tab:summ_accuracy}
	\resizebox{\textwidth}{!}{ 	\begin{tabular}{ccc}
		\hline \hline
		Nonlinear Force & Primary Case (n:1) & Secondary Case (n:1) \\ \hline
		Stiffening Duffing & 0.4\% (3:1) & 68.9\% (5:1) \\ 
		Quintic Stiffness  & 13.9\% (3:1) & 3.1\% (5:1 with $\hat{F} \leq 0.62 $) or 77.3\% (5:1 with $\hat{F} \geq 0.62 $) \\ \hline
		
		Softening Duffing & 1.1\% (3:1) & 35.9\% (5:1) \\ 
		
		Conservative Softening II & 12.4\% (3:1) & 6.1\% (5:1) \\  \hline
		
		Unilateral Spring & 11.6\% (2:1) & 96.2\% (3:1) and 20.9\% (4:1) \\  \hline
		
		Cubic Damping & 5.4\% (3:1) & No Appreciable Resonance \\ \hline
		
		Jenkins Element & 31.2\% (3:1) & 32.4\% (5:1) \\ 
		Iwan Element & 14.4\% (3:1)  & 12.6\% (5:1) \\ 
		\hline \hline
	\end{tabular}}
\end{table}

Shown in \Cref{tab:timing}, HBM requires 12 to 248 times more computation time than VPRNM.\footnote{When comparing cases with the same continuation step size.} 
Here, the stiffening Duffing, Jenkins, and Iwan nonlinearities are considered to show the range of computational benefits (see \Cref{sec:full_timing} for computation times of all nonlinearities).
\Cref{tab:timing} considers the maximum step size used for all of the plots in the present paper
and 10 times that maximum step size.\footnote{Increasing the step size for continuation decreases the computation time while reducing resolution. In general, using the same step size for HBM and VPRNM can be considered as giving a similar accuracy.} 
For the larger step sizes, only the HBM solutions for the Jenkins and Iwan elements show some faceting. For these cases, the step size for VPRNM could be further increased achieving larger speedups compared to HBM while maintaining a similar accuracy to HBM.
Here, the stiffening Duffing nonlinearity shows the largest speedup for VPRNM because of the large number of harmonics used and the large frequency range used for HBM.\footnote{The large frequency range was required because the superharmonic resonance shifts in frequency and to ensure that the initial point was away from any superharmonic resonances (e.g., 5:1, 7:1 etc.) to reliably converge.} 
On the other hand, the Jenkins and Iwan models show smaller decreases in computation time because fewer harmonics are used over a smaller frequency range for HBM and the hysteretic models require more computation for the additional AFT evaluation for VPRNM.
Improvements to the AFT evaluation for the hysteretic nonlinearities discussed in \Cref{sec:aft} contribute to limiting the additional computational cost of the second AFT evaluation for VPRNM.
In all cases, the timing results in \Cref{tab:timing} and \Cref{sec:full_timing} show a clear improvement for using VPRNM rather than HBM at discrete force levels.
The computation times for HBM could be decreased by better initializing the HBM continuation closer to the superharmonic resonance. VPRNM solutions provide one such way to find points to initialize HBM.

\FloatBarrier
\begin{table}[h!]
	\centering
	\caption{Timing comparison between HBM and VPRNM method for 3:1 superharmonic resonances using solver settings identical to the previous sections. Additional setting details for HBM and VPRNM are provided for reference. Simulations are timed on a desktop computer (Intel i7-10710U CPU, 1.10 GHz processor with 6 cores and 32 GB of RAM). HBM computations are for a single run and VPRNM is averaged over 20 runs.} 
	\label{tab:timing}
	
	\resizebox{\textwidth}{!}{ 	\begin{tabular}{cccc}
		
		\hline \hline
		& Stiffening Duffing & Jenkins Element & Iwan Element 
		\\
		\hline
		HBM Time (sec) & 1140 & 237 & 967 \\ 
		HBM with 10$ \times $Plotted Step Size Time (sec) & 160 & 34.1 & 122 \\ 
		\hline
		VPRNM Time (sec) & 4.6 & 13 & 11 \\ 
		VPRNM with 10$ \times $Plotted Step Size Time (sec) & 0.73 & 2.8 & 1.6 \\
		\hline
		
		Harmonics (HBM and VPRNM) & 0, 1-12 & 0,1-3 & 0, 1-3 \\
		Number of HBM Force Levels & 25 & 30 & 30 \\
		HBM Frequency Range (rad/s) & 0.25-1.25 & 0.2-0.4 & 0.2-0.4  \\
		
		\hline \hline
		
	\end{tabular}}
	
\end{table}

The present paper has only considered the case of SDOF systems; future work will consider VPRNM for MDOF systems. The exact behavior of VPRNM for MDOF systems remains unknown, however the formulation in \Cref{sec:tracking_method} has been developed so that it can be applied to MDOF systems with minimal modifications. Furthermore, analysis of the SDOF system has provided potential insights into superharmonic and internal resonances for eight different nonlinear forces (e.g., transitions from local maximum to local minima for hysteretic models). 
Future development for VPRNM could also seek to analyze where the approximation breaks down based on the forces related to superposition ($ F_{kc,sup,n} $ and $ F_{ks,sup,n} $ from \eqref{eq:fsup}). 
Lastly, some of the errors in the VPRNM solutions presented here are attributed to considering more strongly nonlinear regimes than previous studies based on perturbation methods \cite{volvertPhaseResonanceNonlinear2021, volvertResonantPhaseLags2022, abeloosControlbasedMethodsIdentification2022a}.

\section{Conclusions}
\label{sec:concl}

The present work analyzes superharmonic resonances for eight different nonlinear forces showing a range of characteristics applied to a single degree of freedom system. 
A new method termed variable phase resonance nonlinear modes (VPRNM) is proposed for tracking superharmonic resonances and extends the concept of phase resonance nonlinear modes (PRNM). VPRNM utilizes a phase difference between the internal forces exciting the superharmonic resonance and the superharmonic response to identify and track the evolution of superharmonic resonances over varying external force levels. VPRNM generalizes PRNM to allow for easier application with arbitrary and numerically described nonlinear forces including hysteretic forces. 
The present work evaluates VPRNM for stiffening, softening, even, damping, and hysteretic nonlinearities.
The major conclusions are summarized as follows:
\begin{itemize}
	
	\item VPRNM reduces computation time compared to the harmonic balance method (HBM) by up to a factor of 248 while identifying behavior of superharmonic resonances such as phase shifts and transitions from local maxima to local minima that could easily be missed with HBM.
	
	\item VPRNM may fail to exactly track the superharmonic resonance in some cases (e.g., 5:1 superharmonics for stiffening nonlinearities). However, VPRNM solutions still provide points that could be used to initialize HBM near the superharmonic resonance, reducing computational costs.
	
	\item VPRNM accurately tracks superharmonic resonances for the considered hysteretic models and thus is a promising method for applications to jointed structures.
	
	\item Superharmonic resonances occur for hysteretic and saturating nonlinear forces over a limited range of external forcing amplitudes. Therefore, it is hypothesized that internal resonances for multiple degree of freedom (MDOF) systems with hysteresis require commensurate frequencies to occur at an appropriate amplitude. 
	This results in less frequently observed internal resonances than for other nonlinearities (e.g., Duffing) where commensurate frequencies need only occur above a certain threshold to be prominent.
	
	\item VPRNM performs well when the fifth harmonic is directly excited by the fundamental motion. Conversely, VPRNM has higher errors for cases where an interaction between the fundamental and third harmonic motions excites the fifth harmonic (e.g., Duffing nonlinearities).
	
	\item While the present work has only considered SDOF systems, VPRNM has been formulated to be easily generalized to MDOF systems as preliminary work has demonstrated \cite{porterTrackingSuperharmonicInternal2024}. Additionally, the case of cubic damping has illustrated that VPRNM can handle shifts in the phase of the fundamental motion near fundamental resonances that are expected for modal interactions in MDOF systems. Considerations of the accuracy of a generalized VPRNM for MDOF systems is left to future work.

\end{itemize}
The proposed VPRNM method provides a useful tool for characterizing superharmonic resonances at significantly reduced computational cost and has been applied to provide insights into superharmonic resonance behavior for a range of nonlinear forces.

\section*{Acknowledgments}
\label{sec:acknowledgments}

Funding: 
This material is based upon work
	supported by the U.S. Department of Energy, Office of
	Science, Office of Advanced Scientific Computing
	Research, Department of Energy Computational Science
	Graduate Fellowship under Award Number(s) DE-SC0021110.
The authors are thankful for the support of the National
Science Foundation under Grant Number 1847130.

This report was prepared as an account of
work sponsored by an agency of the United States
Government. Neither the United States Government nor
any agency thereof, nor any of their employees, makes
any warranty, express or implied, or assumes any legal
liability or responsibility for the accuracy, completeness,
or usefulness of any information, apparatus, product, or
process disclosed, or represents that its use would not
infringe privately owned rights. Reference herein to any
specific commercial product, process, or service by trade
name, trademark, manufacturer, or otherwise does not
necessarily constitute or imply its endorsement,
recommendation, or favoring by the United States
Government or any agency thereof. The views and
opinions of authors expressed herein do not necessarily
state or reflect those of the United States Government or
any agency thereof.

\printbibliography

\appendix
\section{Additional Example Frequency Response Curves} \label{sec:appendix_frcs}

\begin{figure}[!h]
	\centering
	\newcommand{\exwidth}{0.85}  
			\begin{subfigure}{\exwidth\linewidth}
				\centering
				\includegraphics[width=\linewidth]{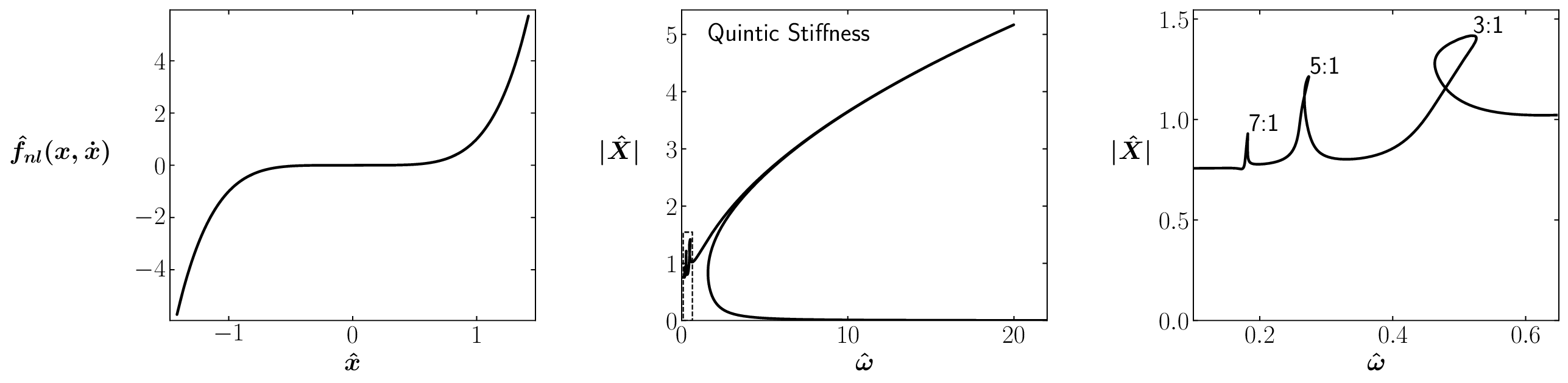}
				\caption{}
			\end{subfigure}
			\begin{subfigure}{\exwidth\linewidth}
	\centering
	\includegraphics[width=\linewidth]{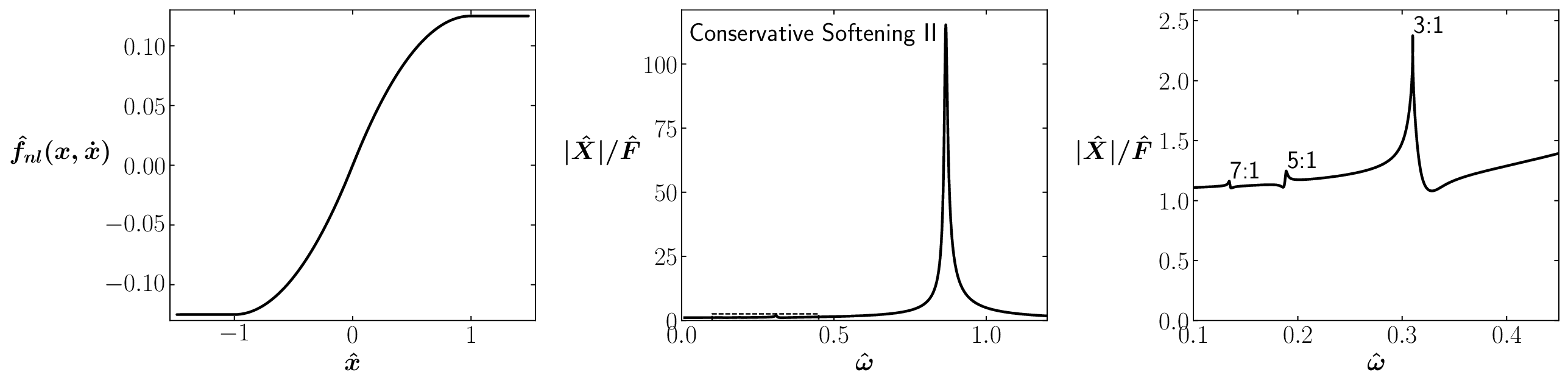}
	\caption{}
	\end{subfigure}
		\begin{subfigure}{\exwidth\linewidth}
	\centering
	\includegraphics[width=\linewidth]{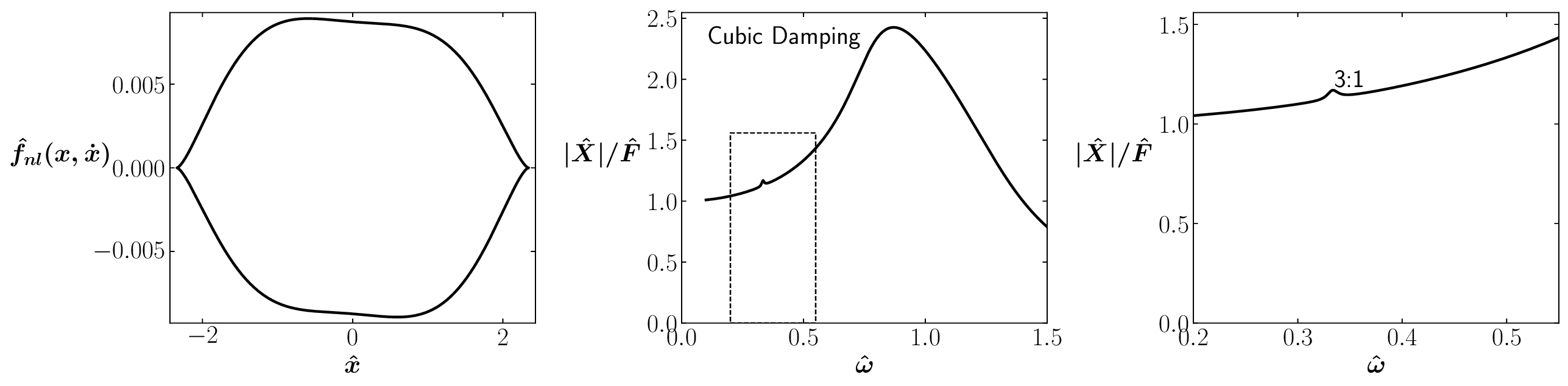}
	\caption{}
	\end{subfigure}
		\begin{subfigure}{\exwidth\linewidth}
	\centering
	\includegraphics[width=\linewidth]{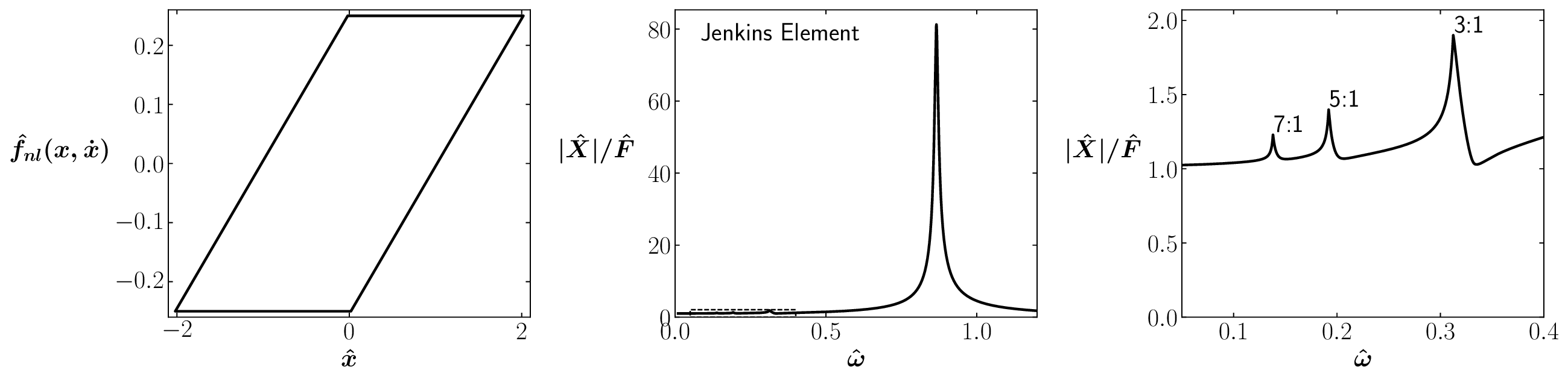}
	\caption{}
	\end{subfigure}
		\caption{Examples of FRCs illustrating superharmonic behavior calculated with harmonics 0 and 1 through 8. The left is the force displacement relationship for the 3:1 superharmonic resonance, the middle is the FRC over the full frequency range, and the right is a zoomed in version for the dashed box on the middle plot. The systems are 
			(a) quintic stiffness with $ \hat{F} = 1.0 $,
			(b) conservative softening II with $ \hat{F} = 0.625 $,
			(c) cubic damping with $ \hat{F} = 2.0 $, 
			and
			(d) Jenkins element with $ \hat{F} = 1.0625 $.
	} 
	\label{fig:example_FRC_appendix}
\end{figure}

\Cref{fig:example_FRC_appendix} shows FRCs with superharmonic responses for the quintic stiffness, conservative softening II, cubic damping, and Jenkins nonlinearities. Examples of FRCs for the other nonlinear forces are shown in \Cref{fig:example_FRC}.

\FloatBarrier

\FloatBarrier

\section{Harmonic Balance Method} \label{sec:hbm}

The harmonic balance method (HBM) assumes that the nonlinear vibration system responds in a periodic motion that is a sum of sine and cosine terms at integer multiples of the forcing frequency \cite{krackHarmonicBalanceNonlinear2019}. Combining \eqref{eq:x_approx} with the equation of motion \eqref{eq:system} yields the discrete set of equations
\begin{equation}
\begin{split}
	k X_0 + F_{nl, 0} - F_{ext, 0} &= 0
	\\
	(-n \omega^2 m + k) X_{nc} + (n \omega c) X_{ns} + F_{nl, nc} - F_{ext, nc} &= 0 \ \ 
	\forall n \in \{1, \dots, H\}
	\\
	(-n \omega^2 m + k) X_{ns} - (n \omega c) X_{nc} + F_{nl, ns} - F_{ext, ns} &= 0 \ \ 
	\forall n \in \{1, \dots, H\}
\end{split}
\end{equation}
Here, $ F_{ext, 0} $ is static external forcing. Additionally, $ F_{ext, nc} $ and $ F_{ext, ns} $ represent external forcing applied to the $ n $th harmonic cosine or the $ n $th harmonic sine respectively. The terms $ F_{nl, 0} $, $ F_{nl, nc} $, $ F_{nl, ns} $ are calculated from the harmonic components of the nonlinear forces (see \Cref{sec:aft}). 

HBM solutions are tracked for the FRCs utilizing continuation. Continuation augments a set of equations with an additional constraint and finds solutions for a range of an additional unknown parameter. For the case of HBM, the additional unknown parameter becomes the frequency of the response. The present work uses a tangent predictor and an orthogonal corrector similar to that described in \cite{peetersNonlinearNormalModes2009, rensonNumericalComputationNonlinear2016}. Continuation is also applied with VPRNM, in which the additional unknown parameter is the external forcing level.
The code for the present implementation of harmonic balance and continuation is made available \cite{porterTMDSimPy}.

\subsection{Alternating Frequency-Time Method}\label{sec:aft}

The nonlinear force terms in HBM are calculated with the Alternating Frequency-Time (AFT) method. For AFT, the current displacement time history $ x(t) $ is calculated.\footnote{All simulations use 1024 time steps of size $ \Delta t $ to discretize time in the interval of $ [ 0, 2 \pi / \omega - \Delta t ] $ unless otherwise noted.} Then, the nonlinear forces are calculated for the time series (repeating twice for hysteretic models to reach steady state \cite{krackHarmonicBalanceNonlinear2019}). Finally, a Fourier transform of the resulting steady-state time series of forces is evaluated to obtain the nonlinear forces in the frequency domain.

To improve computation times, an alternative version of the Jenkins and Iwan force evaluations is implemented for AFT. 
The alternative procedure still uses the same discrete history of $ x(t) $ as the normal AFT procedure. 
However, the algorithm then identifies critical time instants, $ t_{crit} $, where
\begin{equation} \label{eq:tcrit}
	\text{sgn}\left[ x(t_{crit} + \Delta t) - x(t_{crit})\right] 
	\neq 
	\text{sgn}\left[ x(t_{crit}) - x(t_{crit} - \Delta t)\right]
\end{equation}
for the discrete time $ \Delta t $ between entries in $ x(t) $. These time instants correspond to the local maximum and minimum displacements at velocity reversals. 
Then, two repeats of the hysteresis loops are evaluated for only the ordered points at times $ t_{crit} $ for the given friction model. 
Lastly, the frictional force is calculated at each time instant using the previous state as the displacement and force at the value of $ t_{crit} $ that occurs most closely before the time of the given instant. 
This gives an identical time series of nonlinear forces to the standard AFT algorithm because the nonlinear friction evaluations are equivalent when using the previous instant or using the previous velocity reversal point.
The standard steps from the AFT procedure are then applied to convert the time series of forces back into the frequency domain.

The present algorithm sped up the full AFT friction evaluations by a factor of 54 for the Jenkins model and 17 for the Iwan model compared to the standard approach of serially evaluating two hysteresis loops at all time instants.\footnote{The test used 1024 time instants and harmonics 0, 1, 2, and 3 on a desktop computer (Intel i7-10710U CPU, 1.10 GHz processor with 6 cores and 32 GB of RAM). The time was averaged over 1000 evaluations for Jenkins and 200 for Iwan, both corresponding to about 1.5 second of evaluations for the present algorithm.}
The speedup can be attributed to using vectorized operations for both \eqref{eq:tcrit} and the evaluation of the frictional forces at time instants between values of $ t_{crit} $ since these calculations can be trivially parallelized. The benefit of vectorization is especially significant in the present work since Python and NumPy are utilized for the computations \cite{harris2020array}. The Iwan model shows less speedup because significant vectorization is already utilized in the serial version for the 101 slider evaluations per frictional evaluation. 
However, implementation of the algorithm to use multiple processors or GPUs could allow for greater speedup.

Here, the critical path of the computation is independent of the number of times steps used in the AFT procedure.\footnote{Assuming sufficient number of discrete times are included to capture all velocity reversals.} The identification of $ t_{crit} $ in the limit of infinite time steps is equivalent to determining the zeros of a Fourier series (the velocity) with $ H $ terms. This Fourier series has a maximum of $ 2 H $ unique roots \cite{boydComputingZerosMaxima2007}. Therefore, a maximum of $ 4H + 1 $ friction evaluations are required along the critical path corresponding to evaluating each of $ 2H $ roots twice to reach steady state plus a single evaluation for all time instants that can be trivially parallelized. Since a maximum of 5 harmonics is used in this paper for the hysteretic nonlinearities, the critical path has a maximum of 21 friction evaluations\footnote{A total of 1034 friction evaluations are required corresponding to the 10 critical points evaluated twice and 1014 other points. However, the last 1014 evaluations can be fully parallelized.} compared to 2048 corresponding to two full evaluations of the time series.
Additionally, the proposed algorithm only requires 1034 friction evaluations in the present case.
The evaluation of \Cref{eq:tcrit} can easily be parallelized or vectorized, so the additional computation is not significant compared to the drastic reduction in friction evaluations and the length of the critical path. It is noted that a similar approach could be applied to speed up AFT evaluations of the elastic dry friction model that uses a normal load dependent slip limit. Such an algorithm would require additional critical points based on normal load conditions such as those constructed in \cite{brakeMasingManifoldsReconciling2023}. However, such nonlinearities are not applicable to the present work with SDOF systems so are not further considered here.

\section{A Priori Phase Calculations for Secondary Superharmonics} \label{sec:apriori_secondary}

This section analyzes a second set of superharmonic resonances corresponding to larger integer ratios.
The cases of primary secondary superharmonics are analyzed in \Cref{sec:apriori_phase}.
Some of these are directly excited by fundamental harmonic motion as presented in \Cref{tab:Fbroad_second}. However, the 5:1 superharmonic resonance for the Duffing oscillator, the 3:1 superharmonic resonance for the unilateral spring, and the 5:1 superharmonic resonance for cubic damping are not directly excited by fundamental harmonic motion. For the Duffing oscillator, quintic stiffness, and cubic damping the excitation of the fifth superharmonic is analyzed for motion of 
\begin{equation}
	x(t) = 
	X_1 \cos(\omega t) 
	+ 
	X_3 \cos(3 \omega t - \phi_{broad,3})
\end{equation}
This is consistent with the analysis in \Cref{sec:decomposeFnl} in that it is assumed that the third harmonic oscillates in phase with the broadband excitation until the 3:1 superharmonic resonances, which occurs at a higher frequency than the 5:1 superharmonic resonance. For simplicity, only the fundamental motion of \eqref{eq:fundamentalMotion} is considered for the conservative softening II, unilateral spring, Jenkins, and Iwan nonlinearities. 
The observation that the fundamental motion of some nonlinear forces excites the fifth harmonic while other nonlinear forces do not provides an interesting difference that may be useful to future attempts to understand superharmonic resonances.

\FloatBarrier

\begin{table}[h!]
	\centering
	\caption{Analytical calculations of $ F_{nq,broad} $ (see \eqref{eq:fbroad}) excitation of secondary superharmonic resonances. The equations and values presented here are independent of the parameters chosen in \Cref{tab:nl_params}.}
	\label{tab:Fbroad_second}
	\resizebox{\textwidth}{!}{ 		\begin{tabular}{cccccc}
			\hline \hline
			Nonlinear Force & \multicolumn{1}{c}{\begin{tabular}[c]{@{}c@{}}Excited\\ Harmonic\end{tabular}} 
			& $ F_{nc,broad} $ [N] & $ F_{ns,broad} [N] $ & $ \phi_{broad,n} $ [rad] & $ \phi_n $ [rad]
			\\ \hline
			
			Stiffening Duffing & 5 & $ 3 \alpha (X_1^2 X_3 - X_1 X_3^2) / 4 $ & 0 & 0 or $ \pi $ & $ \pm  \pi/2 $
			\\ \hline
			Quintic Stiffness & 5 & 
						\multicolumn{1}{r}{\begin{tabular}[r]{@{}r@{}} $ -\eta (X_1^5 - 20 X_1^4 X_3 + 30 X_1^3 X_3^2 $\\ $ - 30 X_1^2 X_3^3 + 20 X_1 X_3^4)/16 $ \end{tabular}} 
			& 0 & 0 or $ \pi $ & $ \pm  \pi/2 $
			\\ \hline
			
			Softening Duffing & 5 & $ -3 \alpha (X_1 X_3^2 + X_1^2 X_3 ) / 4 $ & 0 & 0 & $ \pi/2 $
			\\
			Conservative Softening II & 5 & \multicolumn{2}{c}{\Cref{fig:fbroad_mag}} & $ -\pi $ & $ -\pi/2 $
			\\ \hline
			
			Unilateral Spring & 4 & $ 2 k_{nl} X_1 / 15 \pi $  & 0 & 0 & $ \pi/2 $
						\\ \hline

			Cubic Damping & 5 & $ -9 \gamma \omega^3 X_1^2 X_3 /4$ & $ -27 \gamma \omega^3 X_1 X_3^2 / 4$ & Variable & Variable
			\\
			Jenkins Element & 5 & \multicolumn{2}{c}{\Cref{fig:fbroad_mag}} & \Cref{fig:fbroad_phase}  & Variable
			\\
			Iwan Element & 5 & \multicolumn{2}{c}{\Cref{fig:fbroad_mag}} & \Cref{fig:fbroad_phase}  & Variable

			\\ \hline\hline
		\end{tabular}
	}
\end{table}

As was done for the case of the primary superharmonic resonances, the magnitude of the broadband excitation of the higher harmonic can inform the expected behavior of the superharmonic responses. Here, the stiffening Duffing, quintic stiffness, softening Duffing, and cubic damping nonlinearities could all be expected to show increasing superharmonic resonances with increasing amplitudes. On the other hand, the conservative softening II nonlinearity, the Jenkins element, and the Iwan element will likely show a peak at intermediate amplitudes because of the saturating nature of the nonlinear force. As before, the unilateral spring results in an excitation that is proportional to the amplitude and may result in superharmonic resonances at all force levels.

Similar to the case of the primary superharmonic resonances, expected phase criteria for each of the nonlinear forces can be obtained by inspecting the harmonic components of the nonlinear force. For the stiffening Duffing oscillator, a phase of 0 or $ \pi $ for $ \phi_{broad,n} $ is possible depending on the relative magnitudes of $ X_1 $ and $ X_3 $. However, one could expect that $ X_3 < X_1 $ since neither harmonic is in resonance and the excitation of the third harmonic is less than the nonlinear restoring force for the fundamental harmonic of the Duffing oscillator. In that case, the phase resonance criteria would be $ \phi_n = \pi/2 $, which is consistent with the analysis of \cite{abeloosControlbasedMethodsIdentification2022a}.
For the quintic stiffness, it would again be expected that $ X_3 < X_1 $. However, given the large coefficients on the polynomial terms, it is possible that the sign may change. Therefore, the phase criteria for superharmonic resonance could be $ \pm \pi/2 $ (this is confirmed in \Cref{sec:results}). 
Finally, the cubic damping and hysteretic nonlinearities all show variable phases for the excitation of the fifth harmonic. This illustrates the need for a general tracking method rather than a constant phase criteria.

\FloatBarrier

\section{Computation Time} \label{sec:full_timing}

\Cref{tab:full_timing} shows computation times for HBM and VPRNM for all eight nonlinear forces for the primary superharmonic resonances presented in \Cref{sec:results}. HBM computations are timed for a single run and VPRNM computations are averaged over 20 repeated runs.
HBM requires between 18 and 248 times longer than VPRNM (for Jenkins and stiffening Duffing respectively).

\FloatBarrier
\begin{table}[h!]
	\centering
	\caption{Timing comparison between HBM and VPRNM method for 2:1 or 3:1 superharmonic resonances using solver settings identical to \Cref{sec:results}. Additional setting details for HBM and VPRNM are provided for reference. Simulations are timed on a desktop computer (Intel i7-10710U CPU, 1.10 GHz processor with 6 cores and 32 GB of RAM).} 
	\label{tab:full_timing}
	
	\resizebox{\textwidth}{!}{ 	\begin{tabular}{cccccc}
		
		\hline \hline
		Model & 
		\begin{tabular}[c]{@{}c@{}} HBM Time\\ (sec) \end{tabular} & 
		\begin{tabular}[c]{@{}c@{}} VPRNM Time\\ (sec) \end{tabular} & 
		\begin{tabular}[c]{@{}c@{}} Harmonics \\ (HBM and VPRNM) \end{tabular} & 
		\begin{tabular}[c]{@{}c@{}} HBM Frequency \\ Range (rad/s) \end{tabular} & 
		\begin{tabular}[c]{@{}c@{}} Number of HBM \\  Force Levels \end{tabular} 
		\\
		\hline
		Stiffening Duffing & 1140 & 4.6 & 0,1-12 & 0.01-1.25 & 25 \\
		Quintic Stiffness & 1190 & 10 & 0,1-12 & 0.01-2.0 & 25\\ \hline
		Softening Duffing & 125 & 4.6 & 0,1-3 & 0.1-0.4 & 20\\
		Conservative Softening II & 145 & 3.4 & 0,1-3 & 0.2-0.4 & 30\\ \hline
		Unilateral Spring & 379 & 3.3 & 0,1-12 & 0.35-0.65 & 20  		\\ \hline
		Cubic Damping & 465 & 3.5 & 0,1-3 & 0.27-0.4 & 20\\
		Jenkins Element & 237 & 13 & 0,1-3 & 0.2-0.4 & 30\\
		Iwan Element & 967 & 11 & 0,1-3 & 0.2-0.4 & 30\\
		\hline \hline
		
	\end{tabular} }
	
\end{table}

\end{document}
\typeout{get arXiv to do 4 passes: Label(s) may have changed. Rerun}

%% file: preamble.tex
\Large
\noindent \textbf{Tracking Superharmonic Resonances for Nonlinear Vibration of Conservative and Hysteretic Single Degree of Freedom Systems}

\vspace{24 pt}

\large

\noindent Justin H. Porter and Matthew R. W. Brake

\noindent Department of Mechanical Engineering, Rice University, Houston,
TX 77005

\vspace{1 in}

\normalsize

\noindent This is a preprint of the journal paper 
``Tracking Superharmonic Resonances for Nonlinear Vibration of Conservative and Hysteretic Single Degree of Freedom Systems.'' Please cite the final version available here: \url{doi.org/10.1016/j.ymssp.2024.111410}
\\
\\\\

\noindent The full final citation is:
\\
\\
Porter, J. H., Brake, M. R. W., 2024. Tracking superharmonic resonances for nonlinear vibration of conservative and hysteretic single degree of freedom systems. Mechanical Systems and Signal Processing 215, 111410. 

\thispagestyle{empty}